\documentclass{article}

\usepackage{arxiv}
\usepackage{latexsym,amsmath,amssymb,amsthm}
\usepackage[utf8]{inputenc} 
\usepackage[T1]{fontenc}    
\usepackage{hyperref}       
\usepackage{url}            
\usepackage{booktabs}       
\usepackage{amsfonts}       
\usepackage{nicefrac}       
\usepackage{microtype}      
\usepackage{natbib}
\usepackage{doi}

\usepackage{graphicx, caption}
\usepackage{subcaption}
\usepackage{graphics}
\usepackage{tikz}
\usetikzlibrary{positioning, calc, arrows.meta, shapes.geometric, fit,bending,matrix}
\graphicspath{ {./fig/} }
\DeclareMathOperator*{\argmin}{\arg\min}

\title{Impilict Runge-Kutta based sparse identification of governing equations in biologically motivated systems}

\author{ Mehrdad Anvari \\
	Department of applied mathematics\\
	University of Tabriz\\
	Tabriz, 51666-16471 \\
	\texttt{anvari@tabrizu.ac.ir} \\
	\And
	Hamidreza Marasi \\
	Department of applied mathematics\\
	University of Tabriz\\
	Tabriz, 51666-16471 \\
	\texttt{marasi@tabrizu.ac.ir} \\
	 \AND
	 Hossein Kheiri \\
	Department of applied mathematics\\
	University of Tabriz\\
	Tabriz, 51666-16471 \\
	\texttt{h-kheiri@tabrizu.ac.ir} \\
}

\date{}

\renewcommand{\shorttitle}{Impilict Runge-Kutta based sparse identification of governing equations}

\begin{document}
\maketitle

\begin{abstract}
Identifying governing equations in physical and biological systems from datasets remains a long-standing challenge across various scientific disciplines, providing mechanistic insights into complex system evolution. Common methods like sparse identification of nonlinear dynamics (SINDy) often rely on precise derivative estimations, making them vulnerable to data scarcity and noise. This study presents a novel data-driven framework by integrating high order implicit Runge-Kutta methods (IRKs) with the sparse identification, termed IRK-SINDy. The framework exhibits remarkable robustness to data scarcity and noise by leveraging the lower stepsize constraint of IRKs. Two methods for incorporating IRKs into sparse regression are introduced: one employs iterative schemes for numerically solving nonlinear algebraic system of equations, while the other utilizes deep neural networks to predict stage values of IRKs. The performance of IRK-SINDy is demonstrated through numerical experiments on benchmark problems with varied dynamical behaviors, including linear and nonlinear oscillators, the Lorenz system, and biologically relevant models like predator-prey dynamics, logistic growth, and the FitzHugh-Nagumo model. Results indicate that IRK-SINDy outperforms conventional SINDy and the RK4-SINDy framework, particularly under conditions of extreme data scarcity and noise, yielding interpretable and generalizable models.
\end{abstract}

\keywords{Model discovery \and machine learning \and deep learning \and  system identification \and  sparse regression  \and dynamical systems}

\section{introduction}

The discovery of equations governing physical and biological systems facilitates an increasing and comprehensive mechanistic understanding of the evolution and progression of complex systems, concurrently propelling research across numerous scientific disciplines alongside the expansion of datasets\cite{Brunton2025, Lorenzo2023, Metzcar2024, Rai2020}. The derivation of these equations from a foundational set of laws in mathematics, physics, and biology, provides interpretable and generalizable frameworks for describing and forecasting a variety of phenomena, but the direct derivation of these models from empirical data is challenging. By combining fundamental modeling principles with empirical data, numerous mathematical models have been established across the extensive domains of biology\cite{Jabbari2016}, epidemiology \cite{Dickson2024}, and pharmacology \cite{Chelliah2021} to investigate diverse problems, including the characterization of interaction networks among cells and proteins\cite{Glass2021}, metabolic networks\cite{Li2023a}, dynamics of tumor growth\cite{Kazerouni2020}, complex tumor-immune interactions \cite{Kirschner1998}, population dynamics \cite{Newman2014}, the propagation of diseases \cite{Sood2023}, as well as pharmacokinetic-pharmacodynamic models \cite{Kalaria2020}. Notable examples encompass predator-prey models \cite{Anjos2023, Gutierrez2023} and competition models \cite{Tang2021, Kareva2015}, which represent the foundational elements of systems biology of cancer \cite{Hamilton2022}. Moreover, within the exploration of diverse biological phenomena, oscillatory dynamics are frequently encountered: the presence of negative feedback, nonlinear behaviors, and the complex interactions inherent in biological processes pose significant obstacles to traditional frameworks for discovering main differential equations \cite{Prokop2024}. Even in scenarios where there exists a partial knowledge of the phenomenon under investigation, it is impractical to rely exclusively on these fundamental principles \cite{Lejarza2022}.

Over the past decades, considerable efforts have been directed toward the development of data-driven approaches for modeling and prediction within systems biology \cite{Cortés2022, Keating2020}, achieving notable success, particularly in the areas of genomic and medical image analysis\cite{Liu2021, Wong2024, Zhang2024, Ma2024}. Although these methods are capable of capturing the underlying structures of a biological system \cite{Cortés2022, Zhao2021, Galindez2023}, the predominant categories of these approaches lack interpretability in resulting models through neural networks \cite{Zhao2022}, random forest \cite{Jiao2021}, spectral methods \cite{Park2023}, and etc, as well as provide limited insights into the underlying mechanisms governing the evolution of biological processes and interactions \cite{AlQuraishi2021} that lead to e.g., oscillatory behavior. In order to address this challenge, innovative hybrid methodologies such as physics-informed machine learning \cite{Karniadakis2021, Raissi2019} and systems biology-informed deep learning \cite{Yazdani2020, Lagergren2020} have been developed, which not only facilitate the discovery of data patterns but also significantly enhance the interpretability of the findings. A subclass of these methodologies, refered to as universal differential equations \cite{Rackauckas2020}, delineates the simplest aspects of the model as knowledge-based terms and the more complex components as unknown terms by leveraging the universal approximation capabilities of deep neural networks.

The direct derivation of governing differential equations from the collected data leads to a better understanding and interpretation of the unknown mechanisms underlying the processes \cite{Brunton2016a, Brunton2016b}. Recently, such methodologies inspired by the foundational concepts of regression, such as symbolic regression \cite{Schmidt2009, Udrescu2020, Orzechowski2018} and sparse identification of nonlinear dynamics, SINDy \cite{Brunton2025, Brunton2016a}, have been effectively proposed. Although symbolic regression can be used to identify the nonlinear differential equations that describe the dynamic behavior of the underlying system, this approach is computationally prohibitive \cite{Virgolin2022}, restricting their applicability to real-world problems. Additionally, symbolic regression approaches can be prone to overfitting.

The conventional SINDy \cite{Brunton2016a} and its various extensions , including Implicit-SINDy \cite{Mangan2016}, as well as SINDy-PI \cite{Kaheman2020} to account identifying rational functions, SINDYc \cite{Brunton2016b, Kaiser2018} to account for control input, reactive SINDy \cite{Hoffmann2019} and Ensemble-SINDy \cite{Fasel2022}, leverage the principle of parsimony and the sparsity inherent in the governing equations within the high-dimensional space of nonlinear functions \cite{Mangan2016}, thereby yielding models that are more interpretable, and generalizable. The implementation of most of these extensions is included in the open source module PySINDy \cite{Kaptanoglu2021}. These methodologies are witnessing an escalation in their significance, particularly within the disciplines of physics \cite{Alves2022}, chemistry \cite{Hoffmann2019}, biology \cite{Mangan2016}, and engineering \cite{Li2023}, and have demonstrated remarkable efficacy when applied to synthetic datasets across a multitude of benchmark problems addressing both ordinary differential equations (ODEs)\cite{Brunton2016a} and partial differential equations (PDEs) \cite{Rudy2017}, including the discernment of wide range of dynamical behaviour such as linear and nonlinear systems with limit cycles and oscillations \cite{Prokop2024}, periodic and chaotic systems \cite{Kaptanoglu2023, Bramburger2020}, bifurcations \cite{Mangan2016}, multistability \cite{Qian2023}, as well as empirical datasets such as gene expression data \cite{Sandoz2023}, data pertaining to predator-prey dynamics \cite{Anjos2023, Gutierrez2023}, and in vitro datasets associated with chimeric antigen receptor (CAR) T-cell functionality \cite{Brummer2023}.

The main challenge of model discovery is learning dynamic  in scenarios of data scarcity and noisy data \cite{Fasel2022}. All developed versions of SINDy require accurate derivative approximations, which impose substantial constraints on the practical sampling time stepsizes of the data. Furthermore, the presence of noise within the datasets leads to considerable disturbances in derivative approximations \cite{Xu2021}. Schaeffer and McCalla \cite{Schaeffer2017} introduced integral formulation of SINDy to address this challenge. Messenger and Bortz \cite{Messenger2021a, Messenger2021b} by proposing Weak SINDy (WSINDy) extended this weak formulation to provide better robustness to noise.  Goyal and Benner \cite{Goyal2022}, by conceptualizing the integration of sparse identification with the classical fourth-order Runge-Kutta method \cite{Hairer1993}— termed  RK4-SINDy—have alleviated  the requirement for derivative approximation, similar to the study conducted on linear multistep methods by Chen \cite{Chen2023}, thereby prposing a novel framework for the sparse identification of dynamical systems. The enhanced efficacy of RK4-SINDy in comparison to conventional SINDy is attributable to the superior order of convergence associated with RK4. Nonetheless, explicit methods employed in the numerical solution of differential equations, such as RK4, are constrained by a limited stability region, inevitably encountering stepsize limitations \cite{Wanner1996}. Consequently, A-stable methods are necessary for effectively addressing many nonlinear problems, encompassing oscillators and stiff systems \cite{Wanner1996, Butcher2016}. Conversely, in accordance with the second Dahlquist barrier \cite{Dahlquist1963}, there are no A-stable linear multistep methods that exceed order $2$.

For any specified order, there exists an A-stable implicit Runge-Kutta method, also referred to as IRK, corresponding to that order \cite{Butcher2016, Butcher1964}. The reduced stepsize limitations and high accuracy of IRKs render these methods particularly suitable for the identification of dynamical systems, especially those of biological relevance \cite{Fasel2022, Xu2021}. Implicit Gauss methods constitute the A-stable class of fully implicit Runge-Kutta methods exhibiting the highest accuracy \cite{Butcher2016, Sato2023}, thus establishing them as an ideal choice for addressing stiff problems. In this manuscript, we introduce and evaluate two potential approaches for integrating IRKs with sparse identification to facilitate data-driven discovery of governing differential equations through the utilization of integral forms of these equations. The first approach is based on computation of stage values of IRKs  by solving the nonlinear system of algebraic equations. The second approach is founded on the prediction of stage values of IRKs by leveraging the universal approximation capabilities of deep neural networks.

The subsequent sections of the manuscript are organized as follows: Methods section provides two approaches based on IRKs for the purpose of learning governing equations through sparse regression techniques. Additionally, Results section evaluates the performance of the proposed frameworks utilizing synthetic datasets through a series of numerical experiments, while also comparing the results with existing approaches. Finally, conclusions and discussion are presented in last section, accompanied by a brief summary of prospective research directions.

\section{Methods}
In this section, employing implicit Range-Kutta methods, we propose our methodology to discover governing equations from sparsely sampled observations. The collected data may exhibit significant levels of noise. At the core of this approach is the combination of high-order numerical integration techniques, known for their robust numerical stability properties, with the intrinsic sparsity of the governing equations in nonlinear functions space. We consider a collection of candidate functions as a nonlinear feature library, from which the governing equations can be expressed as a linear combination of them. Subsequently, they are identified utilizing noisy data by specifying active terms in the library through a sparse regression technique. This strategy does not rely explicitly on derivative computations, making it more robust to data scarcity and noise, in contrast to the conventional approach \cite{Brunton2016a} and RK4-SINDy \cite{Goyal2022}. In this work, we employ Gauss methods as implicit Range-Kutta methods and assess the principal advantage of integrating them in the sparse identification of dynamical systems within benchmark problems.

\subsection{Implicit Runge-Kutta methods}
Runge-Kutta methods are extensively utilized to solve initial value problems due to the capability to construct them from any specified order \cite{Wanner1996, Butcher2016, Butcher1964}. Implicit Runge-Kutta methods exhibit A-stability properties and are widely recognized as a highly suitable candidate for addressing issues associated with stiffness \cite{Wanner1996, Butcher2016}. IRKs of higher orders impose fewer restrictions on the stepsize, thereby can play a crucial  role in the sparse identification of dynamical systems with constrained data availability. Consider the subsequent initial value problem:
\begin{equation}\label{eqdy}
\begin{cases}
\dot{x}(t) = f(x(t)),\\
x(t_{0}) = x_{0},
\end{cases}
\end{equation}
where the vector $x(t) = \begin{bmatrix}
x_{1}(t) & \dots & x_{d}(t)
\end{bmatrix}^{T}$ indicates the state variables and $f: \mathbb{R}^{d} \to \mathbb{R}^{d}$ represents the unidentified vector field that requires to be determined from data corrupted by noise. To identify the nonlinear dynamics of eq.\eqref{eqdy}, i.e. $f(.)$, it is assumed that the data $\{ x(t_{k}) \}_{k=0}^{m}$  is sampled at times $\{ t_{0}, t_{1}, \dots, t_{m}\}$ with stepsizes $h_{k} = t_{k+1}-t_{k}$. Our approach is based on comparing the data $x(t_{k+1})$ and the predicted value of $x$ at time $t_{k+1}$ based on data $x(t_{k})$ \cite{Goyal2022, Chen2023} using an IRK method. Hence, let us utilize the general form of Runge-Kutta methods with s stages \cite{Butcher2016} for eq.\eqref{eqdy}:
\begin{subequations}
\begin{equation}
x(t_{k+1}) \approx x(t_{k}) + h_{k}\sum_{j=1}^{s}b_{j} f(\chi_{j}^{f}(t_{k})),\label{eqrka}
\end{equation}
\begin{equation}
\chi_{i}^{f}(t_{k}) =  x(t_{k}) + h_{k} \sum_{j=1}^{s}a_{ij}f(\chi_{j}^{f}(t_{k})), \quad i=1, \dots, s,\label{eqrkb}
\end{equation}
\end{subequations}
where $\chi_{i}^{f}(t_{k}) \approx x(t_{k} + c_{i}h_{k})$, for all $i \in \{1, \dots s\}$ with respect to function $f$ in eq.\eqref{eqdy}, that can also be expressed in vectorized form as in eq.\eqref{vector}, with respect to notations \eqref{notation}:
\begin{subequations}
\begin{equation}\label{notation}
\chi =
\begin{bmatrix}
\chi_{1}^{f}\\
\vdots\\
\chi_{s}^{f}
\end{bmatrix}, \quad
F(\chi) = \begin{bmatrix}
f(\chi_{1}^{f})\\
\vdots\\
f(\chi_{d}^{f})
\end{bmatrix},
\end{equation}
\begin{equation}\label{vector}
\begin{bmatrix}
\begin{array}
{c}
\chi(t_{k})\\
\hline
x(t_{k+1})
\end{array}
\end{bmatrix}
\approx
\begin{bmatrix}
\begin{array}
{c|c}
A\otimes \mathbb{I}_{d} & \textbf{1}_{s} \otimes \mathbb{I}_{d}\\
\hline
b^{T}\otimes\mathbb{I}_{d} &  1
\end{array}
\end{bmatrix}
\begin{bmatrix}
\begin{array}
{c}
h_{k} F(\chi(t_{k}))\\
\hline
x(t_{k})
\end{array}
\end{bmatrix},
\end{equation}
\end{subequations}
where $\mathbb{I}_{d}$ denotes $d\times d$ identity matrix, $\textbf{1}_{s}$ is $s$-element unit vector, and $\otimes$ represents Kronecker produc. The versatility of this general formulation allows for both implicit and explicit time-stepping schemes, depending on the selection of $A$. An implicit method with $s$ stages can achieve at most order $2s$ for Gauss methods \cite{Wanner1996, Butcher2016, Butcher1964}. Our methodology for discovering governing equations of the system focuses on utilizing high-order IRKs. The key idea of the methodology is improving the identified vector field through minimizing the loss function used in the training process to quantify satisfaction of the governing ODE model through the reconstruction of data based on the stages values of IRKs.

\subsection{Discovering nonlinear differential equations with IRKs}
The local error of the $s$-stage Gauss method is $O(h_{k}^{2 s})$ \cite{Butcher2016}; Therefore, for sufficiently small stepsizes, $h_{k}$, it has high accuracy in calculating $x(t_{k+1})$ from $x(t_{k})$. If we consider the right hand side of eq.\eqref{eqrka} as a function of $f$, $x(t_{k})$ and $h_{k}$, denote by $\mathcal{F}_{irk}$; the correctness of $x(t_{k+1}) \approx \mathcal{F}_{irk}(f, x(t_{k}), h_{k})$ follows directly:
\begin{equation}\label{Firk}
\mathcal{F}_{irk}(f, x(t_{k}), h_{k}) := x(t_{k}) + h_{k}\sum_{j=1}^{s}b_{j}f(\chi_{j}^{f}(t_{k})),
\end{equation}
where, it is possible to predict $x(t_{k+1})$ and $x(t_{k-1})$ values with high accuracy  using IRKs through $x(t_{k})$. We aim to present the IRK-based approach for the sparse identification of nonlinear dynamics. Specifically, our goal is to derive the most parsimonious representation of the vector field $f(x(t))$ by leveraging the time series data of $x(t)$ at the time instances $\{ t_{0}, \dots, t_{m}\}$ that was previously assumed. For this purpose, we consider the data matrices in the following manner:
\begin{equation}\label{xlr}
X^{L} : =
\begin{bmatrix}
x^{T}(t_{0})\\
x^{T}(t_{1})\\
\vdots\\
x^{T}(t_{m-1})
\end{bmatrix},\quad
X^{R} : =
\begin{bmatrix}
x^{T}(t_{1})\\
x^{T}(t_{2})\\
\vdots\\
x^{T}(t_{m})
\end{bmatrix},
\end{equation}
and represent predictions as eq.\eqref{xlrf}:
\begin{equation}\label{xlrf}
X_{\mathcal{F}}^{R}(f) =
\begin{bmatrix}
\mathcal{F}_{irk}(f, x(t_{0}), h_{0})\\
\mathcal{F}_{irk}(f, x(t_{1}), h_{1})\\
\vdots\\
\mathcal{F}_{irk}(f, x(t_{m-1}), h_{m-1})
\end{bmatrix}, \quad
X_{\mathcal{F}}^{L}(f) =
\begin{bmatrix}
\mathcal{F}_{irk}(f, x(t_{0}), -h_{0})\\
\mathcal{F}_{irk}(f, x(t_{1}), -h_{1})\\
\vdots\\
\mathcal{F}_{irk}(f, x(t_{m-1}), -h_{m-1})
\end{bmatrix}.
\end{equation}
The subsequent essential step in sparse identification involves the generation of a large library of nonlinear features, denoted as $\Phi$, which encompasses potential nonlinear functions that can play a role in the dynamics. It is assumed that the function $f(.)$ can be expressed as a linear combination of a few selected terms derived from the library \cite{Brunton2016a}. Illustratively, one could opt for a collection of polynomials, exponential functions, as well as trigonometric functions within the library. Upon considering the vector $x=[x_{1}, \dots, x_{d}]^{T}$, a library may be given as:
\begin{equation}\label{dict}
\Phi(x) := [1, x, x^{\mathcal{P}_{2}}, \dots, x^{\mathcal{P}_{d}}, \dots, \sin(x), \cos(x), \dots, \exp(-x), \exp(-2x), \dots],
\end{equation}
where $x^{\mathcal{P}_{i}}$ represents polynomials of the degree $i$. To exemplify, in the scenario where $d=2$, $x^{\mathcal{P}_{2}}$ is given as follows:
\begin{equation*}
x^{\mathcal{P}_{2}} = [x_{1}^{2}(t), x_{1}(t)x_{2}(t), x_{2}^{2}(t)].
\end{equation*}
Each element within the $\Phi$ library stands as a suitable candidate for representing $f$. Moreover, depending on the specific context, a collection of meaningful features can be systematically or empirically devised for inclusion within the library.

Despite the possession of an extensive library, numerous choices for candidates will inevitably arise. The primary objective, however, revolves around identifying the minimal feasible candidate subset for the nonlinear representation of the function $f$ \cite{Brunton2016a, Goyal2022}. As a result, a regularized optimization problem is formulated to select a limited number of candidate functions from the library. In this context, it is imperative to note that if $\xi^{*}$ represents the sparse coefficient matrix derived from the optimization problem and $f_{k}: \mathbb{R}^{d} \to \mathbb{R}$ denotes the $k$th component of the function $f$, and $\xi_{k}^{*}$ $k$th collumn of $\xi^{*}$, the following must be satisfied:
\begin{equation}\label{eq1}
f_{k}(x) \approx \sum_{j} \Phi_{j}(x)\xi^{*}_{j,k} = \Phi(x)\xi^{*}_{k},
\end{equation}
Alternatively, in a more comprehensive manner:
\begin{equation}\label{eq1extend}
f(x) \approx \Phi(x) \xi^{*},
\end{equation}
Consequently, appropriately selecting candidate functions from the library determines the governing equations. Special attention must be given to eq.\eqref{eq2} while formulating the appropriate optimization problem:
\begin{equation}\label{eq2}
X^{i} = X^{i}_{\mathcal{F}}(f),
\end{equation}
for $i=L, R$. Acording to eqns.\eqref{eq1extend} and \eqref{eq2}, we can define loss function \eqref{loss} with respect to $\xi$ matrix:
\begin{equation}\label{loss}
\mathcal{L}(\xi) = \alpha \| X^{L} - X^{L}_{\mathcal{F}}(\Phi(.)\xi) \|_{2}^{2} + (1 - \alpha) \| X^{R} - X^{R}_{\mathcal{F}}(\Phi(.)\xi) \|_{2}^{2},
\end{equation}
where $0 \leq \alpha \leq 1$. The corresponding regularized optimization problem \eqref{optprob} can be formulated as \cite{Goyal2022, Chen2023, Messenger2021a, Xu2021}:
\begin{equation}\label{optprob}
\xi^{*} 
= 
\argmin_{\xi} 
\{ \mathcal{L}(\xi) +  \mathcal{R}(\xi, \lambda)\},
\end{equation}
where, $\mathcal{R}(\xi, \lambda)$ represents the regularization term with thresholding parameter $0 \leq \lambda \leq 1$. A choice is the utilization of $\mathit{l}_{1}$-regularization \cite{Sun2020, Cortiella2021}, defined as $\mathcal{R}(\xi, \lambda) = \lambda \| \xi \|_{1}$.
\begin{figure}
\flushleft
\begin{tikzpicture}[
node distance=0.0cm and 0.0cm,
neuron/.style={circle, draw, minimum size=0.6cm},
block/.style={rectangle, draw, rounded corners, fill=teal!30, minimum width=2cm, minimum height=1cm},
>=Stealth,
every node/.style={align=center} 
]

\begin{scope}[local bounding box=NN, rounded corners=5pt, minimum width=5cm, minimum height=0cm]

\node[block, fill=teal!5, rounded corners=5pt, minimum width=0.25\textwidth, minimum height=8.4cm, dashed] (datablock) {};

\node[block, fill=teal!0, below right=of datablock.north west, yshift=-0.2cm, xshift = 0.4cm] (image) {\includegraphics[width=0.18\textwidth]{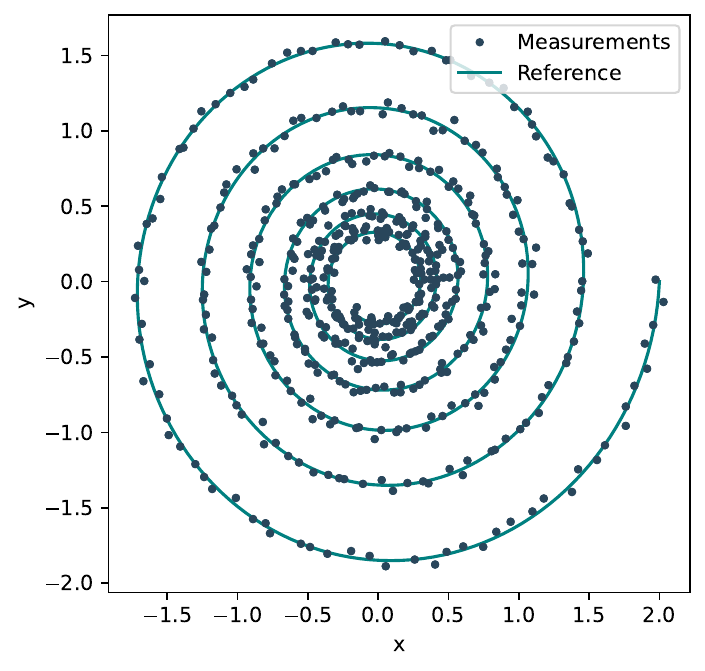}};

\node[block, fill=teal!10, rounded corners=5pt, minimum width=2.5cm, minimum height=0.4cm, below=of image, yshift=-0.35cm, text width=3.0cm, ,align=left] (data) {{\tiny $\begin{bmatrix}x_{1}(t_{0})& \dots &x_{d}(t_{0})\\x_{1}(t_{1})& \dots &x_{d}(t_{1})\\\vdots&  &\vdots\\x_{1}(t_{m})& \dots &x_{d}(t_{m})\end{bmatrix}$}};

\node[block, fill=teal!0, above right=of datablock.south west, yshift=0.2cm, xshift = 0.4cm] (image2) {\includegraphics[width=0.18\textwidth]{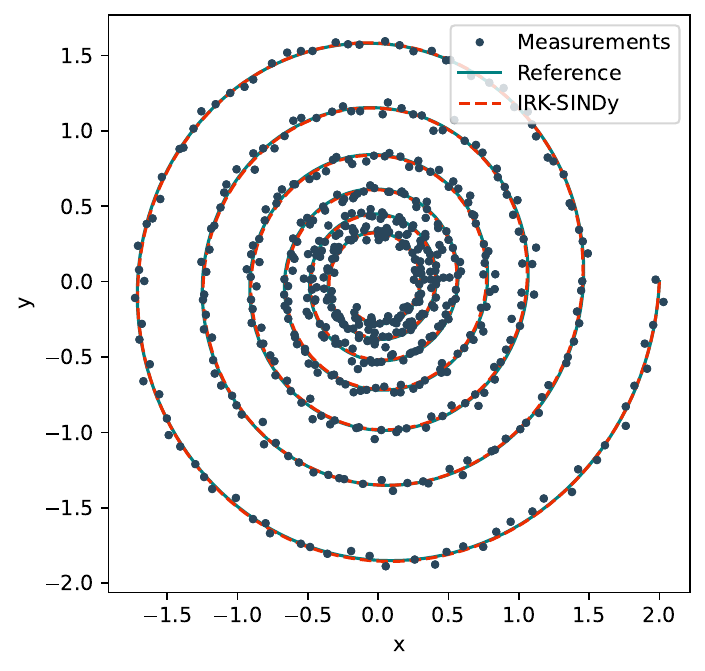}};

\node[above=of datablock.north]{\textbf{(a)}};

\node[block, fill=teal!5, rounded corners=5pt, minimum width=0.75\textwidth, minimum height=3.8cm, dashed, below right =of datablock.north east, xshift=0.6cm] (stageblock) {};

\node[above=of stageblock.north]{\textbf{(b)}};

\node[block, fill=teal!20, rounded corners=5pt, minimum width=1.1cm, minimum height=1.1cm, below right=of stageblock.north west, xshift=0.2cm, yshift=-0.2cm] (stage) {$\chi_{\xi}^{[0]}(t_{k})$};

\node[block, fill=blue!20, rounded corners=5pt, minimum width=1.3cm, minimum height=0.8cm, right=of stage, xshift=1.2cm, text width=1.3cm] (iter) {{\small Iterative schemes}};

\node[block, fill=teal!20, rounded corners=5pt, minimum width=6.8cm, minimum height=0.4cm, right=of iter, xshift=0.9cm, yshift=-0.0cm, text width=6.8cm, ,align=left] (fp) {{\small $\chi^{[i+1]} = (e\otimes I_{d})x(t_{k}) + h_{k}(A\otimes I_{d})\Phi(\chi^{[i]})\xi, \quad$ (I)}};

\node[block, fill=teal!20, rounded corners=5pt, minimum width=2.5cm, minimum height=0.8cm, below=of fp, yshift=-0.1cm, xshift=0.9cm, text width=5.0cm, ,align=left] (func) {{\tiny$\Psi(\chi) = \chi - (e\otimes I_{d})x(t_{k}) - h_{k}(A\otimes I_{d})\Phi(\chi)\xi$}};

\node[block, fill=teal!20, rounded corners=5pt, minimum width=4.9cm, minimum height=0.8cm, below=of func, yshift=-0.1cm, xshift=0.1cm, text width=4.9cm, ,align=left] (newton) {{\small $\begin{cases}J[\Psi](\chi^{[i]})\Delta = -\Psi(\chi^{[i]}), \\ \chi^{[i+1]} = \chi^{[i]} + \Delta.\end{cases}\quad$ (II)}};

\node[block,fill=teal!40,rounded corners=5pt, minimum width=4.8cm, minimum height=1.4cm, above right=of stageblock.south west,  yshift=0.2cm, xshift=1.6cm] (predstage) {};

\node[block, fill=teal!80, rounded corners=5pt, minimum width=1cm, minimum height=1cm, right=0.12cm of predstage.west] (stage1) {$\chi^{\xi}_{1}(t_{k})$};

\node[block, fill=teal!80, rounded corners=5pt, minimum width=1cm, minimum height=1cm, right=of stage1, xshift=0.3cm] (stage2) {$\chi^{\xi}_{2}(t_{k})$};

\node[right=of stage2, minimum width=0.1cm] (dots) {$\dots$};

\node[block, fill=teal!80, rounded corners=5pt, minimum width=1cm, minimum height=1cm, right=of dots] (stages) {$\chi^{\xi}_{s}(t_{k})$};

\node[block, below=of stageblock.south, fill=teal!5, rounded corners=5pt, minimum width=0.75\textwidth, minimum height=4.0cm, dashed, yshift=-0.6cm] (irkblock) {};

\node[above=of irkblock.north]{\textbf{(c)}};

\node[block, fill=teal!30, rounded corners=5pt, minimum width=1.1cm, minimum height=1.1cm, right=of irkblock.west, xshift=0.2cm, yshift=1.3cm] (datak) {$X(t_{k})$};

\node[right=of datak, minimum width=0.1cm, xshift = -0.2cm, yshift=-1.0cm] (sum1) {$\bigoplus$};

\node[block, fill=teal!30, rounded corners=5pt, minimum width=0.65cm, minimum height=2.5cm, right=of sum1, xshift=1.1cm] (s1) {\rotatebox[origin=c]{90}{$\Phi(\chi^{\xi}_{1}(t_{k})) \xi$}};

\node[right=of s1, minimum width=0.1cm, xshift=1.1cm] (t1) {$\bigoplus$};

\node[block, fill=teal!30, rounded corners=5pt, minimum width=0.65cm, minimum height=2.5cm, right=of t1, xshift=0.7cm] (s2) {\rotatebox[origin=c]{90}{$\Phi(\chi^{\xi}_{2}(t_{k}))\xi$}};

\node[right=of s2, minimum width=0.1cm, xshift=1.1cm] (t2) {$\bigoplus\dots\bigoplus$};

\node[block, fill=teal!30, rounded corners=5pt, minimum width=0.65cm, minimum height=2.5cm, right=of t2, xshift=0.7cm] (ss) {\rotatebox[origin=c]{90}{$\Phi(\chi^{\xi}_{s}(t_{k}))\xi$}};

\node[minimum size=0.0cm, draw, right=of ss, xshift=0.9cm] (t3) {$\approx$};

\node[block, fill=teal!30, rounded corners=5pt, minimum width=1.2cm, minimum height=1.2cm, above right=of irkblock.south west, xshift=0.2cm, yshift=0.2cm] (datakp1) {$X(t_{k+1})$};

\node[above=of datak.north, xshift=-0.2cm, yshift=-0.2cm] (fakep) {};
\draw[-] (stageblock.east) -- ++(1.0, 0.0);
\draw[->] (stageblock.east) ++(1.0, 0.0) |- (irkblock.east);

\draw[-] (iter.east) -- (fp);
\draw[->] (stage) -- (iter);
\draw[->] (iter.south west) -| (stage1.north);
\draw[->] (iter.south)++(0.3cm, 0.0cm) to (stage2.north);
\draw[->] (iter.south east) -| (stages.north);

\draw[->] (image.south) -- (data.north);
\draw[->] (data.east)++(0.0, 0.7cm) -| (fakep);

\draw[->] (datak) -- (stage);
\draw[->] (datak.south west) |- node[below, xshift=0.4cm]{$\times h_{k}$} (sum1);
\draw[->] (sum1) --  (s1);
\draw[->] (s1) -- node[above, xshift=0.1cm]{$\times h_{k}b_{1}$} (t1);
\draw[->] (t1) -- (s2);
\draw[->] (s2) -- node[above, xshift=0.1cm]{$\times h_{k}b_{2}$} (t2);
\draw[->] (t2) -- (ss);
\draw[->] (ss) -- node[above]{$\times h_{k}b_{s}$} (t3);
\draw[->] (t3.south) |- (datakp1);
\end{scope}

\begin{scope}[local bounding box=NN, rounded corners=5pt, minimum width=5cm, minimum height=0cm]

\node[block, fill=teal!5, rounded corners=5pt, minimum width=\textwidth + 0.6cm, minimum height=5.9cm, dashed, below right=of datablock.south west, yshift = -0.6cm] (sindyblock) {};

\node[above=of sindyblock.north]{\textbf{(d)}};

\node[block, fill=teal!10, rounded corners=5pt, minimum width=3.8cm, below left=of sindyblock.north east, yshift=-0.2cm, xshift=-0.2cm, text width=3.8cm, ,align=left] (linear) {\begin{equation*}\begin{split}&\dot{x}(t) = -0.1 x(t) + 2.0 y(t),\\ &\dot{y}(t) = -2.0 x(t) + 0.1 y(t).\end{split}\end{equation*}};

\node[below=of linear.north] (rmtitle) {\textbf{Reference model}};


\node[block, fill=teal!50, draw=teal!50, rounded corners=5pt, minimum width=2.8cm, minimum height=1.4cm, solid, below=of sindyblock.north, yshift=-0.2cm, xshift= 0.6cm] (xdata) {};

\node[block, fill=teal!30, rounded corners=5pt, minimum width=1.2cm, minimum height=1.2cm, below right=0.15 cm of xdata.north west] (xl) {$X^{L}$};

\node[block, fill=teal!30, rounded corners=5pt, minimum width=1.2cm, minimum height=1.2cm, right=of xl, xshift=0.2cm] (xl) {$X^{R}$};

\node[block, fill=orange!50!red, rounded corners=5pt, minimum width=1.6cm, above=of sindyblock.south, yshift=1.1cm, xshift= 0.6cm, text width=1.6cm, ,align=center] (irk) {{\small $IRK\big(\Phi(.)\xi\big)$}};



\node[block, fill=teal!50, draw=teal!50, rounded corners=5pt, minimum width=1.4cm, minimum height=2.8cm, solid, above left=0.3cm of sindyblock.south east] (xpred) {};

\node[block, fill=teal!30, rounded corners=5pt, minimum width=1.2cm, minimum height=1.2cm , above left=0.15cm of xpred.south east] (xrf) {$X^{R}_{\mathcal{F}}$};

\node[block, fill=teal!30, rounded corners=5pt, minimum width=1.2cm, minimum height=1.2cm, above=of xrf, yshift=0.2cm] (xlf) {$X^{L}_{\mathcal{F}}$};

\node[ellipse, draw, above=of xlf, xshift=-3.0cm, yshift=0.15cm, minimum width=1.1cm, minimum height=0.7cm] (loss) {Loss $\mathcal{L}$};

\node[diamond, minimum height=0.4cm, minimum width=1.5cm, xshift=-0.6cm, yshift=0.0cm, draw, left=of loss, shape aspect=2] (small) {{\small Small?}};

\node[ellipse, draw, below=of small.south, xshift=2.0cm, yshift=-0.1cm, minimum width=1.4cm, minimum height=0.5cm] (dloss) {${\partial\mathcal{L} \over \partial\xi}$};

\node[diamond, minimum height=0.4cm, minimum width=1.5cm, xshift=-1.9cm, yshift=-1.0cm, draw, left=of small, shape aspect=2] (sparse) {{\small Sparse?}};

\node[block, fill=teal!10, rounded corners=5pt, minimum width=3.5cm, above right=of sindyblock.south west, xshift=0.2cm, yshift = 0.2cm, text width=3.0cm, ,align=left] (coefs) {{\tiny
\begin{equation*}
\begin{matrix}
" &
\begin{matrix}\xi_{x} &\xi_{y}
\end{matrix}\\
\begin{matrix}
1\\
x\\
y\\
\vdots\\
.\\
y^{3}
\end{matrix} &
\begin{bmatrix}
[0.00] & [0.00]\\
[-0.10] & [-2.00] \\
[2.00] & [-0.10] \\
\vdots & \vdots\\
[0.00] & [0.00]\\
[0.00] & [0.00]
\end{bmatrix}
\end{matrix}
\end{equation*}}};

\node[below=of coefs.north] (imtitle) {{\small \textbf{Identified coefficients}}};

\node[block, fill=teal!10, rounded corners=5pt, minimum width=4.8cm, below right=0.3cm of sindyblock.north west, text width=4.8cm, ,align=left] (identified) {\begin{equation*}\begin{split}&\dot{x}(t) = -0.100 x(t) + 2.000 y(t),\\ &\dot{y}(t) = -2.000 x(t) + 0.100 y(t).\end{split}\end{equation*}};

\node[block, fill=teal!10, rounded corners=5pt, minimum width=8.3cm, minimum height = 0.0cm, above=0.15cm of sindyblock.south, text width=8.3cm, ,align=center, xshift = 1.1cm] (library) {\begin{tiny}
\begin{equation*}\Phi(x(t), y(t)) = \begin{bmatrix}
1 & x(t) & y(t) & x^{2}(t) & x(t) y(t) & y^{2}(t) & \dots, & y^{3}(t)
\end{bmatrix} \end{equation*}
\end{tiny}};

\node[below=of identified.north] (imtitle) {\textbf{Identified model}};

\draw[->] (linear.west) -- (xdata);
\draw[->] (xdata.south) -- (irk);
\draw[->] (irk.east) -- (xpred);
\draw[->] (xpred.north) |- (loss);
\draw[->] (loss.west) -- (small);
\draw[->] (small.south) |- node[above,  xshift=0.5cm]{{\small No}} (dloss.north west);
\draw[->] (dloss.west) -| node[below, xshift=1.1cm]{{\small Update $\xi$}} (irk.north east);
\draw[->] (small.west) node[above left, xshift=2.2cm]{{\small Yes}} -|  (sparse.north);
\draw[->] (sparse.east) node[above, xshift=0.2cm]{{\small No}} -- ++(1.0cm, 0.0cm);
\draw[->] (sparse.south) node[below left, yshift=0.1cm, xshift=2.1cm]{{\small Yes}}  |- (coefs.east);
\draw[-] (coefs.north)  -> ++(-0.0cm, 1.5cm);

\draw[-] (identified.west) -- ++(-1.0cm, 0.0cm);
\draw[->] (identified.west) ++(-1.0cm, 0.0cm) |- (image2.west) ;
\end{scope}

\end{tikzpicture}
\caption{Overview of the IRK-SINDy framework: a. For each benchmark problem, we perform measurements that incorporate noise and, thereafter form a dataset. Our objective is to construct a model that is parsimonious, interpretable, and possesses generalizability, capable of accurately forecasting reference dynamics. b. Given an appropriate initial guess (e.g., $X(t_{k})$), the stage values of the IRKs are approximated by solving the system of nonlinear equations \eqref{eqrkb} through iterative schemes. In this context, we employ two iterative approaches: (i) fixed point iteration and (ii) Newton's method. c. With the stage values estabilished the subsequent step values are computed according to eq.\eqref{Firk}. This computational process is depicted as the systematic IRK network. d. Within this structured representation of IRK-SINDy, the dataset is classified into two categories: forward and backward, followed by the formation of a symbolic features library comprising candidate nonlinear functions. To solve a nonlinear sparse regression problem using the forward and backward predictions illustrated in (b) and (c), an IRK step is applied, and the loss function is minimized by choosing a suitable optimizer. Following a certain number of epochs, a sparsity-promoting algorithm is employed. Finally, every non-zero element in the coefficient matrix $\xi^{*}$ signifies an active term within the feature library, thereby representing the resultant discovered model.}\label{fig1}
\end{figure}
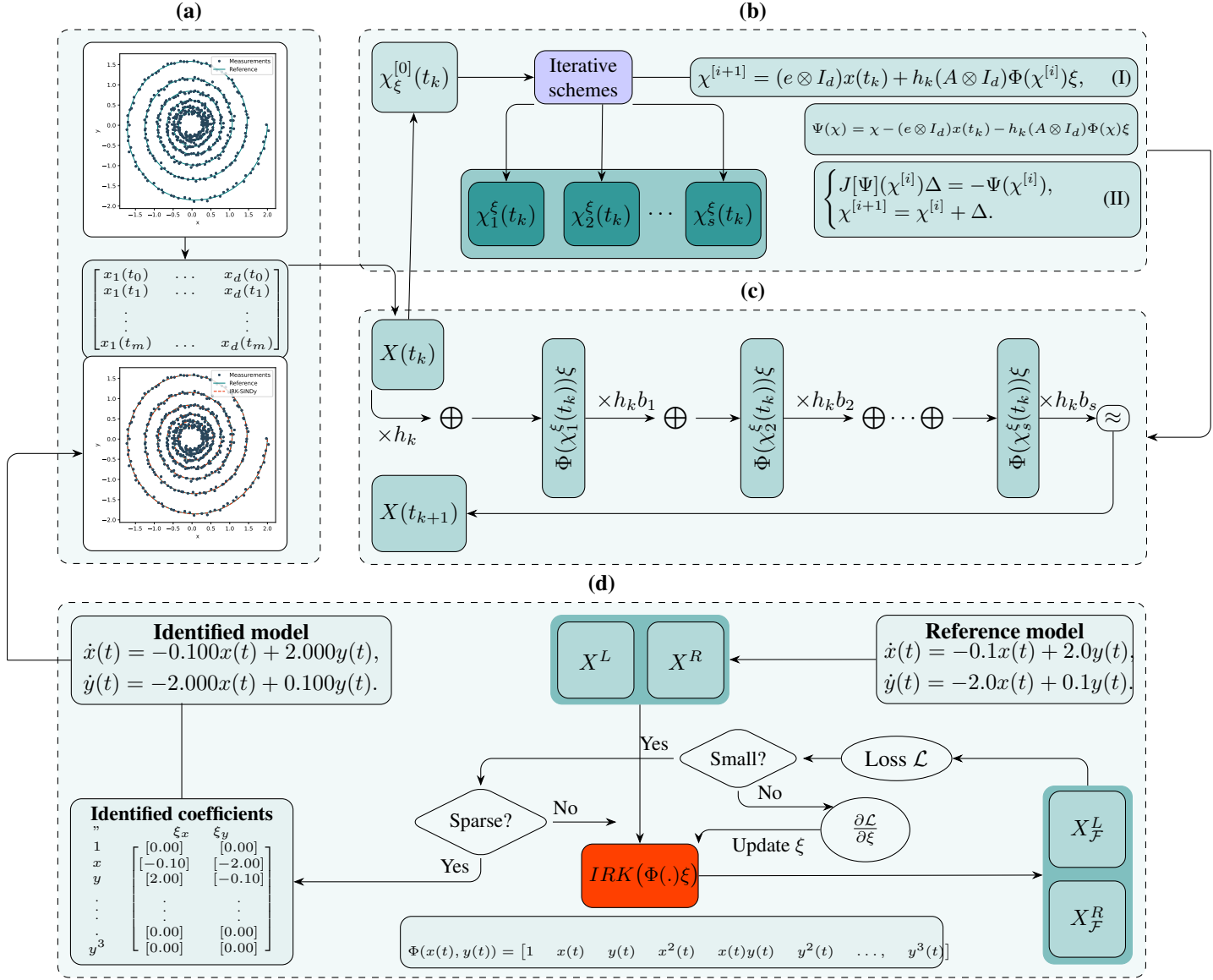

To implement the algorithm, it is crucial to acquire the stage values, $\chi_{i}(t_{k}), i=1, \dots, s$,  within the context of IRKs. In the classical implementation of IRKs, these values are computed by solving the nonlinear
 system of algebraic equations \eqref{eqrkb} employing iterative schemes such as fixed-point iteration and Newton's method \cite{Butcher2016}. We note that, in the implementation of Newton's method, automatic differentiation tool \cite{Baydin2018} exploited to calculate Jacobean matrix. In light of this foundational framework, as illustrated in Figure\ref{fig1}, we propose a novel sparse identification process of differential equations inspired by IRKs. Within this approach, we perform predictive analyses of the quantities $X^{L}$ and $X^{R}$ to address the optimization problem in \eqref{optprob}, which is achieved by obtaining the $sd$ stage values through the aforementioned iterative techniques and subsequently substituting them into eq.\eqref{eqrka}. It is crucial to emphasize that in the context of fixed point methods, the convergence condition for calculating the solution of the system represented by $g(x) = 0$ is depended on Lipschitz constant associated to $g$; conversely, for an initial guess in proximity to the stage values, Newton's method exhibits a rapid convergence to the solution of the system \cite{Butcher2016, Wanner1996}. Therefore, given the limitations of fixed-point methods in solving nonlinear and stiff problems \cite{Butcher1964, Sato2023}, they are inefficient for the sparse identification of nonlinear dynamical systems.

Despite the efficiency of this approach, the necessity of solving the system of nonlinear equations at each epoch significantly slows down the overall optimization process \cite{Jay2000}. Moreover, to enhance the accuracy of the predictions, it is essential to employ higher-order IRKs, and therefore an increased number of stage values. As $s$ increases, the corresponding computational cost increases exponentially associated with these calculations. Therefore, similiar to Raissi et al.\cite{Raissi2019}, we address this computational challenge using an auxiliary deep neural network, DNN, that ensures a linear increase in computational cost in relation to the number of stage values employed.

\subsection{Discovering nonlinear differential equations through combining DNNs and IRKs}
Here we make use of an auxiliary DNN that is parameterized as a nonlinear mapping from time $t_{k}$, along with the corresponding state variable $x$ evaluated at $t_{k}$, i.e. $x(t_{k})$, to the stage values ofIRKs in the approximation of $x(t_{k+1})$with stepsize $h_{k}$. Therefore, we represent the DNN as $\chi^{\theta} = [\chi^{\theta}_{1}, \dots, \chi^{\theta}_{s}]$, wherein $\theta$ serves as the neural network parameters. As elucidated in Figure\ref{fig2}, by training the auxiliary DNN, $\chi^{\theta}$ is expected to yield predictions that approximate the stage values of IRKs in the governing equations of the system, represented as $\chi^{f}$:
\begin{equation*}
\chi^{\theta}_{i}(t_{k}) \approx \chi^{f}_{i}(t_{k}), \quad i=1, \dots, s, \quad k=0,1, \dots, m-1.
\end{equation*}
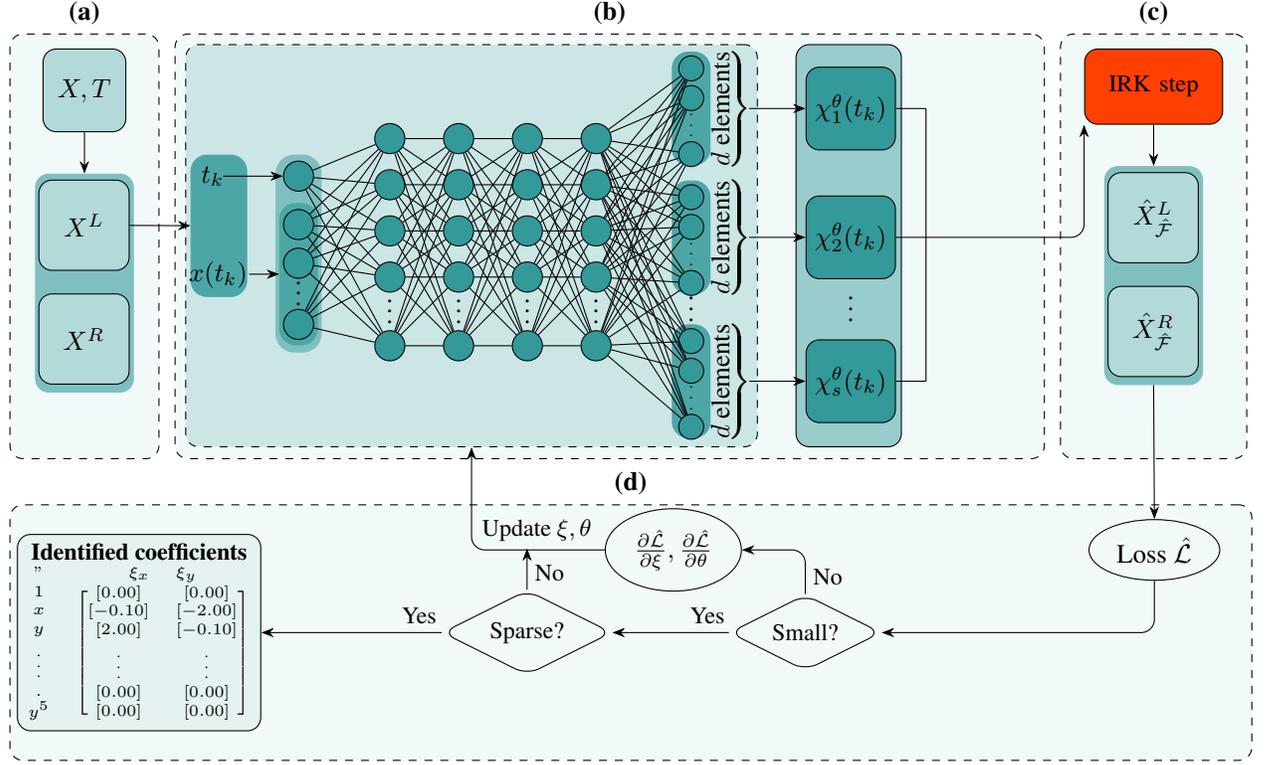
\begin{figure}
\flushleft
\begin{tikzpicture}[
node distance=0.1cm and 0.5cm,
neuron/.style={circle, fill=teal!80, draw, minimum size=0.4cm},
block/.style={rectangle, draw, rounded corners, fill=teal!30, minimum width=2cm, minimum height=1cm},
>=Stealth,
every node/.style={align=center} 
]

\begin{scope}[local bounding box=NN, rounded corners=5pt, minimum width=2cm, minimum height=0cm]

\node[block, fill=teal!5,, rounded corners=5pt, minimum width=0.12\textwidth, minimum height=5.65cm, dashed] (datablock) {};

\node[above=of datablock.north, yshift=-0.1cm]{\textbf{(a)}};

\node[block,fill=teal!30,rounded corners=5pt, minimum width=1.1cm, minimum height=1.1cm, below=0.2cm of datablock.north] (data) {$X, T$};

\node[block, fill=teal!50, draw=teal!50, rounded corners=5pt, minimum width=1.3cm, minimum height=2.9cm, below=0.55cm of data.south] (xpred) {};

\node[block, fill=teal!30, rounded corners=5pt, minimum width=1.2cm, minimum height=1.2cm , above left=0.15cm of xpred.south east, xshift=0.05cm] (xrf) {$X^{R}$};

\node[block, fill=teal!30, rounded corners=5pt, minimum width=1.2cm, minimum height=1.2cm, above=of xrf, yshift=0.2cm] (xlf) {$X^{L}$};

\end{scope}

\begin{scope}[local bounding box=NN, rounded corners=5pt, minimum width=5cm, minimum height=0cm]

\node[block, fill=teal!5, right=0.2cm of datablock.east, rounded corners=5pt, minimum width=0.7\textwidth, minimum height=5.65cm, dashed] (genblock) {};

\node[above=of genblock.north, yshift=-0.1cm]{\textbf{(b)}};

\node[block, fill=teal!20,, rounded corners=5pt, minimum width=7.6cm, minimum height=5.35cm, dashed, below right=0.2cm of genblock.north west] (nnblock) {};
    

\node[block, fill=teal!50, draw=teal!50, rounded corners=7pt, minimum width=0.6cm, minimum height=2.62cm, below right=of nnblock.north west, yshift=-1.37cm, xshift=0.7cm] (inputblock) {};

\node[block, fill=teal!70, draw=teal!70 , rounded corners=6pt, minimum width=0.74cm, minimum height=1.88cm, below right=0.1cm of nnblock.north west, yshift=-1.4cm] (datak) {};

\node[neuron, below =of inputblock.north, yshift=0.02cm] (t) {};

\node[left=of t, xshift=-0.16cm, minimum width=0.1cm, minimum height=0.1cm](tk){$t_{k}$};

\node[block, fill=teal!70, draw=teal!70 , rounded corners=6pt, minimum width=0.5cm, minimum height=1.88cm, below=of t, yshift=-0.05cm] (x0block) {};

\node[neuron, below=0.1cm of x0block.north, yshift=0.02cm] (x01) {};
\node[neuron, below= of x01] (x02) {};
\node[neuron, below= of x02, yshift=-0.3cm] (x0d) {};
\node[below= of x02, yshift=0.4cm] (vdots) {$\vdots$};
\node[left=of x0block, xshift=0.2cm, minimum width=0.2cm, minimum height=0.2cm](xtk){$x(t_{k})$};

\node[neuron, right=of t, yshift=0.5cm, xshift=0.3cm] (h11) {};
\node[neuron, below=of h11, yshift=-0.1cm] (h12){};
\node[neuron, below=of h12, yshift=-0.1cm] (h13){};
\node[neuron, below=of h13, yshift=-0.1cm] (h14){};
\node[neuron, below=of h14, yshift=-0.4cm] (h15){};
\node[below= of h14, yshift=0.35cm] (vdots) {$\vdots$};

\node[neuron, right=of h11] (h21) {};
\node[neuron, below=of h21, yshift=-0.1cm] (h22){};
\node[neuron, below=of h22, yshift=-0.1cm] (h23){};
\node[neuron, below=of h23, yshift=-0.1cm] (h24){};
\node[neuron, below=of h24, yshift=-0.4cm] (h25){};
\node[below= of h24, yshift=0.35cm] (vdots) {$\vdots$};

\node[neuron, right=of h21] (h31) {};
\node[neuron, below=of h31, yshift=-0.1cm] (h32){};
\node[neuron, below=of h32, yshift=-0.1cm] (h33){};
\node[neuron, below=of h33, yshift=-0.1cm] (h34){};
\node[neuron, below=of h34, yshift=-0.4cm] (h35){};
\node[below= of h34, yshift=0.35cm] (vdots) {$\vdots$};

\node[neuron, right=of h31] (h41) {};
\node[neuron, below=of h41, yshift=-0.1cm] (h42){};
\node[neuron, below=of h42, yshift=-0.1cm] (h43){};
\node[neuron, below=of h43, yshift=-0.1cm] (h44){};
\node[neuron, below=of h44, yshift=-0.4cm] (h45){};
\node[below= of h44, yshift=0.35cm] (vdots) {$\vdots$};

\node[block, fill=teal!70, draw=teal!70, rounded corners=6pt, minimum width=0.5cm, minimum height=1.5cm, right=of h41, yshift= 0.4cm, xshift=0.3cm] (xblock1) {};

\node[neuron, below=0.05cm of xblock1.north, minimum size=0.1cm] (chi11) {};
\node[neuron, below= of chi11, minimum size=0.1cm, yshift=0.05cm] (chi12) {};
\node[below= of chi12.south, yshift=0.4cm] (vdots2) {{\tiny $\vdots$}};
\node[neuron, below= of chi12, minimum size=0.2cm, yshift=-0.3cm] (chi1d) {};

\node[block, fill=teal!70, draw=teal!70, rounded corners=6pt, minimum width=0.5cm, minimum height=1.5cm, below=of xblock1, yshift= -0.1cm] (xblock2) {};
\node[neuron, below=0.05cm of xblock2.north, minimum size=0.1cm] (chi21) {};
\node[neuron, below= of chi21, minimum size=0.1cm, yshift=0.05cm] (chi22) {};
\node[below= of chi22, yshift=0.4cm] (vdots3) {{\tiny $\vdots$}};
\node[neuron, below= of chi22, minimum size=0.2cm, yshift=-0.3cm] (chi2d) {};

\node[below= of xblock2.south, yshift=0.4cm] (vdots3) {$\vdots$};

\node[block, fill=teal!70, draw=teal!70, rounded corners=6pt, minimum width=0.5cm, minimum height=1.5cm, below=of xblock2, yshift=-0.3cm] (xblocks) {};
\node[neuron, below=0.05cm of xblocks.north, minimum size=0.2cm] (chis1) {};
\node[neuron, below= of chis1, minimum size=0.2cm, yshift=0.05cm] (chis2) {};
\node[below= of chis2, yshift=0.4cm] (vdots2) {{\tiny $\vdots$}};
\node[neuron, below= of chis2, minimum size=0.2cm, yshift=-0.3cm] (chisd) {};

\node[right=of xblock1, xshift=-2.75cm](c1){\rotatebox[origin=c]{90}{$\underbrace{d \text{ elements}}$}};

\node[right=of xblock2, xshift=-2.75cm] (c2) {\rotatebox[origin=c]{90}{$\underbrace{d \text{ elements}}$}};

\node[right=of xblocks, xshift=-2.75cm] (cs) {\rotatebox[origin=c]{90}{$\underbrace{d \text{ elements}}$}};

\node[block,fill=teal!40,rounded corners=5pt, minimum width=1.4cm, minimum height=5.35cm, right=of nnblock.east, xshift=-0.0cm] (stage) {};

\node[block,fill=teal!80,rounded corners=5pt, minimum width=1.1cm, minimum height=1.1cm, right=of c1, xshift=-2.0cm] (stage1) {$\chi^{\theta}_{1}(t_{k})$};

\node[block,fill=teal!80,rounded corners=5pt, minimum width=1.1cm, minimum height=1.1cm, right=of c2, xshift=-2.0cm] (stage2) {$\chi^{\theta}_{2}(t_{k})$};

\node[below=of stage2, yshift=0.2cm] (vdots4) {$\vdots$};

\node[block,fill=teal!80,rounded corners=5pt, minimum width=1.1cm, minimum height=1.1cm, right=of cs, xshift=-2.0cm] (stages) {$\chi^{\theta}_{s}(t_{k})$};
    
\foreach \i in {t, x01, x02, x0d}
	\foreach \j in {1, 2, 3, 4, 5}
		\draw[-] (\i) -- (h1\j);

\foreach \i in {1, 2, 3, 4, 5}
	\foreach \j in {1, 2, 3, 4, 5}
		\draw[-] (h1\i) -- (h2\j);

\foreach \i in {1, 2, 3, 4, 5}
	\foreach \j in {1, 2, 3, 4, 5}
		\draw[-] (h2\i) -- (h3\j);
		
\foreach \i in {1, 2, 3, 4, 5}
	\foreach \j in {1, 2, 3, 4, 5}
		\draw[-] (h3\i) -- (h4\j);
            
\foreach \i in {1, 2, 3, 4, 5}
	\foreach \j in {1, 2, d}
		\draw[-] (h4\i) -- (chi1\j);

\foreach \i in {1, 2, 3, 4, 5}
	\foreach \j in {1, 2, d}
		\draw[-] (h4\i) -- (chi2\j);

\foreach \i in {1, 2, 3, 4, 5}
	\foreach \j in {1, 2, d}
		\draw[-] (h4\i) -- (chis\j);

\draw[->] (data) -- (xpred);
\draw[->] (xlf) -- (datak);
\draw[->] (tk)++(0.14cm, 0.0cm) -- (t);
\draw[->] (xtk.east) ++(-0.1cm, 0.0cm) -- (x0block);

\draw[->, solid] (c1.east)++(-2.25cm, 0.0cm) -- (stage1.west);
\draw[->, solid] (c2.east)++(-2.25cm, 0.0cm) -- (stage2.west);
\draw[->, solid] (cs.east)++(-2.25cm, 0.0cm) -- (stages.west);

\end{scope}

\begin{scope}[local bounding box=NN, rounded corners=5pt, minimum width=5cm, minimum height=0cm]

\node[block, fill=teal!5, right=0.2cm of genblock.east, rounded corners=5pt, minimum width=0.15\textwidth, minimum height=5.65cm, dashed] (irkblock) {};

\node[above=of irkblock.north, yshift=-0.1cm]{\textbf{(c)}};

\node[block, fill=orange!50!red, rounded corners=5pt, minimum width=1.6cm, below=0.2cm of irkblock.north, text width=1.6cm, ,align=center] (irkstep) {{\small IRK step}};

\node[block, fill=teal!50, draw=teal!50, rounded corners=5pt, minimum width=1.3cm, minimum height=2.9cm, below=0.55cm of irkstep.south] (xpred1) {};

\node[block, fill=teal!30, rounded corners=5pt, minimum width=1.2cm, minimum height=1.2cm , above left=0.15cm of xpred1.south east, xshift=0.05cm] (xr) {$\hat{X}^{R}_{\hat{\mathcal{F}}}$};

\node[block, fill=teal!30, rounded corners=5pt, minimum width=1.2cm, minimum height=1.2cm, above=of xr, yshift=0.2cm] (xl) {$\hat{X}^{L}_{\hat{\mathcal{F}}}$};

\draw[-] (stage1.east) -- ++(0.4cm,0);

\draw[->] (stage2.east) -| (irkstep.south west);

\draw[-] (stages.east) -- ++(0.4cm,0);
\draw[-] (stage1.east) ++(0.4cm,0) -- ++(0.0cm,-3.63cm);

\draw[->] (irkstep) -- (xpred1);

\end{scope}

\begin{scope}[local bounding box=NN, rounded corners=5pt, minimum width=5cm, minimum height=0cm]

\node[block, fill=teal!5, below right=0cm of datablock.south west, rounded corners=5pt, minimum width=\textwidth, minimum height=3.4cm, dashed, yshift=-0.6cm] (diagblock) {};

\node[above=of diagblock.north, yshift=-0.1cm]{\textbf{(d)}};

\node[block, fill=teal!10, rounded corners=5pt, minimum width=3.0cm, right=0.1cm of diagblock.west, text width=3.0cm, ,align=left] (coefs1) {{\tiny
\begin{equation*}
\begin{matrix}
" &
\begin{matrix}\xi_{x} &\xi_{y}
\end{matrix}\\
\begin{matrix}
1\\
x\\
y\\
\vdots\\
.\\
y^{5}
\end{matrix} &
\begin{bmatrix}
[0.00] & [0.00]\\
[-0.10] & [-2.00] \\
[2.00] & [-0.10] \\
\vdots & \vdots\\
[0.00] & [0.00]\\
[0.00] & [0.00]
\end{bmatrix}
\end{matrix}
\end{equation*}}};
\node[below=of coefs1.north, yshift=0.1cm]{{\footnotesize \textbf{Identified coefficients}}};

\node[diamond, minimum height=0.4cm, minimum width=1.5cm, draw, right=2.4cm of coefs1.east, shape aspect=2] (sparse) {{\small Sparse?}};

\node[diamond, minimum height=0.4cm, minimum width=1.5cm, draw, right=1.55cm of sparse, shape aspect=2] (small) {{\small Small?}};

\node[ellipse, draw, right=2.75cm of small, minimum width=1.1cm, minimum height=0.6cm, yshift=1.1cm] (loss) {Loss $\hat{\mathcal{L}}$};

\node[ellipse, draw, right=-0.1cm of sparse, yshift=1.1cm, minimum width=1.4cm, minimum height=0.5cm] (dloss) {${\partial\hat{\mathcal{L}} \over \partial\xi}, {\partial\hat{\mathcal{L}} \over \partial\theta}$};

\draw[->] (xpred1.south) -- (loss.north);
\draw[->] (loss.south) |- (small);
\draw[->] (small) node[above, xshift=-1.3cm]{{\small Yes}}  -- (sparse);
\draw[->] (small.north) node[above, xshift=0.3cm]{{\small No}}  |- (dloss.east);
\draw[->] (sparse.north) node[above, xshift=0.3cm]{{\small No}} -- ++(0.0cm, 0.55cm);
\draw[->] (dloss.west) node[above, xshift=-0.9cm]{{\small Update $\xi, \theta$}}  -| (nnblock.south);
\draw[->] (sparse.west) node[above, xshift=-0.3cm]{{\small Yes}} -- (coefs1.east);
\end{scope}
\end{tikzpicture}
\caption{Overview of the deep IRK-SINDy framework: a. The dataset is prepared for the purpose of training the neural network. b. The inputs to the neural network are assigned into two distinct variables: time and state variables. The neurons located in the output layer of the network are partitioned into $s$ segments, each containing $d$ neurons. The $i$'th segment predicts the $d$ stage values corresponding to the $\chi_{i}$. c. Through the process of forward propagation within the DNN, the stage values are predicted, and these predictions are subsequently employed in the IRK steps, i.e. eq.\eqref{Firk}, facilitating both forward and backward predictions. d. By comparing the predictions against the data, the loss is computed, followed by the optimization step. Upon reaching a specified number of epochs, at which point the loss is sufficiently minimized, the sparsity-promotion algorithm is exclusively applied to the coefficient matrix $\xi$. Finally, the non-zero coefficients of the $\xi$ denote the active terms in the nonlinear feature library.}\label{fig2}
\end{figure}
To facilitate the combination of the auxiliary DNN and IRKs within the framework of the optimization process, it becomes imperative to reformulate eq.\eqref{Firk}. Thus, by subtracting eq.\eqref{eqrka} from eq.\eqref{eqrkb}, we derive eqns. \eqref{eqrkap} and \eqref{eqrkbp}:
\begin{subequations}
\begin{equation}
\hat{\mathcal{F}}_{irk}^{L}(f, x(t_{k}), h_{k}, i) := \chi_{i}(t_{k}) - h_{k}\sum_{j=1}^{s}a_{ij} f(\chi_{j}(t_{k})), \quad i=1, \dots, s,\label{eqrkap}
\end{equation}
\begin{equation}
\hat{\mathcal{F}}_{irk}^{R}(f, x(t_{k}), h_{k}, i) :=  \chi_{i}(t_{k}) + h_{k} \sum_{j=1}^{s}(b_{j} - a_{ij})f(\chi_{j}(t_{k})), \quad i=1, \dots, s,\label{eqrkbp}
\end{equation}
\end{subequations}
It subsequently becomes evident that:
\begin{subequations} \label{eqequival}
\begin{equation}
x(t_{k}) \approx \hat{\mathcal{F}}_{irk}^{L}(f, x(t_{k}), h_{k}, i), \quad i=1,\dots, s,
\end{equation}
\begin{equation}
x(t_{k+1}) = x(t_{k}+h_{k}) \approx \hat{\mathcal{F}}_{irk}^{R}(f, x(t_{k}), h_{k}, i),\quad i=1,\dots, s.
\end{equation}
\end{subequations}
In a manner similar to eqns.\eqref{xlr} and \eqref{xlrf}, the irk network can be systematically defined as eq.\eqref{xrlfnn}:
\begin{equation}\label{xrlfnn}
\hat{X}_{\hat{\mathcal{F}}}^{i}(f) =
\begin{bmatrix}
\hat{\mathcal{F}}^{i}_{irk}(f, x(t_{0}), h_{0}, 0) & \dots & \hat{\mathcal{F}}^{i}_{irk}(f, x(t_{0}), h_{0}, s)\\
\hat{\mathcal{F}}^{i}_{irk}(f, x(t_{1}), h_{1}, 0) & \dots & \hat{\mathcal{F}}^{i}_{irk}(f, x(t_{1}), h_{1}, s)\\
\vdots &  & \vdots\\
\hat{\mathcal{F}}^{i}_{irk}(f, x(t_{m-1}), h_{m-1}, 0)  & \dots & \hat{\mathcal{F}}^{i}_{irk}(f, x(t_{m-1}), h_{m-1}, s)
\end{bmatrix},
\quad i=L,R.
\end{equation}

Now, to select the most active terms among the nonlinear features of the $\Phi$ library, the loss function is formulated in eq.\eqref{lossnn} through the integration of three ingredients, including sparse identification, IRKs, and the auxiliary DNN, with the aim of simultaneously determining the parameters associated with the neural network as well as the coefficient matrix:
\begin{equation}\label{lossnn}
\hat{\mathcal{L}}(\xi) = \alpha \| X^{L} - \hat{X}^{L}_{\hat{\mathcal{F}}}(\Phi(.)\xi) \|_{2}^{2} + (1 - \alpha) \| X^{R} - \hat{X^{R}}_{\hat{\mathcal{F}}}(\Phi(.)\xi) \|_{2}^{2},
\end{equation}
where $0 \leq \alpha \leq 1$. In accordance with eq.\eqref{lossnn} the corresponding regularized optimization problem can be formulated in the following manner \cite{Goyal2022, Chen2023, Messenger2021a, Xu2021}
\begin{equation}\label{optprobnn}
\xi^{*} = \argmin_{\xi, \theta} \{ \hat{\mathcal{L}}(\xi, \theta) +  \mathcal{R}(\xi, \lambda)\},
\end{equation}

This approach enables the simultaneous determination of both $\xi$ and $\theta$, while effectively integrating the auxiliary DNN, IRKs, and sparse identification, referred to as deep IRK-SINDy,  perform the optimization process of discovering the governing differential equations without the need for derivative information. This characteristic is particularly significant as our proposed approach demonstrates not only a remarkable resistance to noise, but also an impressive efficiency in scenarios characterized by data scarcity \cite{Fasel2022, Messenger2021a, Messenger2021b}.

\subsection{Sparsity-promoting procedure}
When the nonlinear optimization problems delineated in eqns.\eqref{optprob} and \eqref{optprobnn} are rigorously formulated, the goal is to seek an approximate sparse solution denoted as $\xi^{*}$ \cite{Goyal2022}. Although several sparse regression methodologies, including but not limited to LASSO \cite{Tibshirani1996}, are frequently employed to promote the sparsity in the resultant solution, it is crucial to note that a significant number of these algorithms are predominantly tailored for linear and convex optimization problems \cite{Golden2024, Messenger2024}. In order to select the active terms within the frameworks of problems \eqref{optprob} and \eqref{optprobnn}, analogous to the sparse regression employed in the sequential thresholding least squares (STLS) algorithm utilized in SINDy \cite{Brunton2016a} (for linear and convex optimization problems), we adopt a gradient-based sequential thresholding procedure \cite{Goyal2022}, as outlined in Figure\ref{fig1} and Figure\ref{fig2}. 

During each iteration of this procedure, the loss functions articulated in eqns.\eqref{loss} and \eqref{lossnn} is initially minimized through the application of a gradient-descent method during the training phase \cite{Luthen2021},which is conducted with respect to the coefficient matrix $\xi$ as well asthe parameter $\theta$ in eq.\eqref{lossnn}. Following a certain number of epochs, sparsity-promoting modifications are made to arrive at a sparse coefficient matrix $\xi$. This stands in contrast to prevalent sparse regression methods like LASSO \cite{Tibshirani1996, Friedman2010} and elastic net \cite{Sun2020}, which modify the loss function by sparsity-promoting penalties, represented as $loss + \lambda_{1} \Vert \xi \Vert_{1} + \lambda_{2} \Vert \xi \Vert_{2}$, and cause the coefficients continuously tend to zero. The $l_{p}-\text{regularization}$ \cite{McCulloch2024}, $p \in \mathbb{N}$, serves as an approximation for the more theoretically desirable $l_{0}$ penalty. In the context of $l_{0}-\text{regularization}$, which is recognized as an NP-hard problem \cite{Natarajan1995}, the magnitude of the coefficients should be matter, a concern that is addressed in a more straightforward manner through the sequential thresholding method \cite{Zhang2019}. Our proposed procedure initiates with the establishment of an initial guess for $\xi$ coupled with setting a threshold value $\lambda$, and at the end of each iteration, any coefficients with absolute value smaller than $\lambda$ are set to zero. This iterative process is sustained until convergence is attained. In practice, when reasonable values of $\lambda$ are employed, the sequential thresholding surprisingly requires few number of iterations to achieve convergence, ultimately leading to the derivation of the optimal coefficient matrix $\xi^{*}$.

It is imperative to underscore the point that when we set $\lambda$ to zero, every term within the nonlinear feature library is recognized as an active term; this scenario is particularly pronounced in instances characterized by measurement errors and numerical round-off effects, which can be deemed non-physical in nature. Furthermore, by specifying $\lambda=1$, the regularization term effectively overcomes the loss function, compelling $\xi$ to approach zero and, consequently, the model to $\dot{x}=0$.  Therefore, the reasonable determination of the sparsity-promoting parameter $\lambda$ through concepts such as cross-validation \cite{Quade2018} and the analysis of the Pareto front \cite{Mangan2017}, which aims to balance the trade-off between loss minimization and model complexity, play a pivotal role in correctly identifying the reference dynamics. In this context, our objective is to impose a penalty on the number of terms within the model, while simultaneously striving to minimize the loss function, thereby yielding the most parsimonious model. To achieve this end, we can utilize the information criterion for model selection delineated in \cite{Mangan2017, Kaptanoglu2023, Dong2023}, a process that has been successfully applied across various sparse identification problems, with each case resulting in the correct identification of the reference model.

\section{Results}
In this section, the efficacy of the proposed methodologies for the data-driven discovery of governing differential equations is demonstrated and examined through a series of numerical experiments on benchmark problems exhibiting varying classified complexity, from linear and nonlinear oscillators to noisy measurement of predator-prey dynamics. The robustness to noise and data scarcity is illustrated in comparison with conventional SINDy \cite{Brunton2016a} (referred to as Conv-SINDy) and RK4-SINDy \cite{Goyal2022} without access to derivative information. The two proposed methodologies have implemented in the PyTorch deep-learning module \cite{Paszke2019} utilizing a gradient descent optimization method alongside Gauss methods up to 500 stages \cite{Raissi2019}, as IRKs with the highest accuracy, to address the optimization problems delineated in eqns.\eqref{optprob} and \eqref{optprobnn}. In this implementation, the Adam optimizer \cite{Kingma2014} has employed to iteratively update the coefficient matrix $\xi$ and the parameters of the neural network $\theta$. All models have trained using an Core i5-7400 CPU with 8 GB memory. Synthetic data are generated by forward-solving the system of differential equations utilizing numerical methods, specifically the fourth-order Runge-Kutta method and Backward Differentiation Formulas (BDF) \cite{Butcher2016, Hairer1993, Wanner1996}, as implemented in the solve$\_$ivp function of SciPy module, followed by introducing various levels of noise into the dataset. Finally, after the model discovery process, the identified model is subjected to a comparative analysis against a reference model in each numerical experiment by evaluating the solutions corresponding to distinct initial conditions.

\subsection{Linear damped oscillator}
As a first illustrative example, we consider discovering the governing equations of a two-dimensional linear damped oscillatory system using data with different noise levels and scarcity, whose dynamics is given by \eqref{loscilator}:
\begin{subequations} \label{loscilator}
\begin{equation}
\dot{x}_{1}(t) = -0.1 x_{1}(t) + 2.0 x_{2}(t),
\end{equation}
\begin{equation}
\dot{x}_{2}(t) = -2.0 x_{1}(t) - 0.1 x_{2}(t).
\end{equation}
\end{subequations}

\begin{figure}[]
\centering
\begin{subfigure}{\textwidth}
\caption{Collected data}\label{fig3a}
\includegraphics[width=0.9\textwidth]{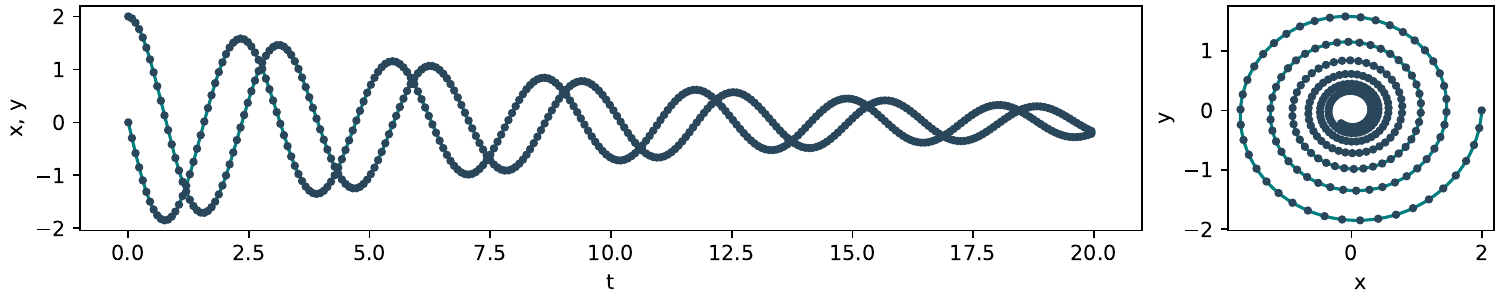}
\end{subfigure}
\begin{subfigure}{\textwidth}
\caption{$m = 801$}
\includegraphics[width=0.9\textwidth, height=3.4cm]{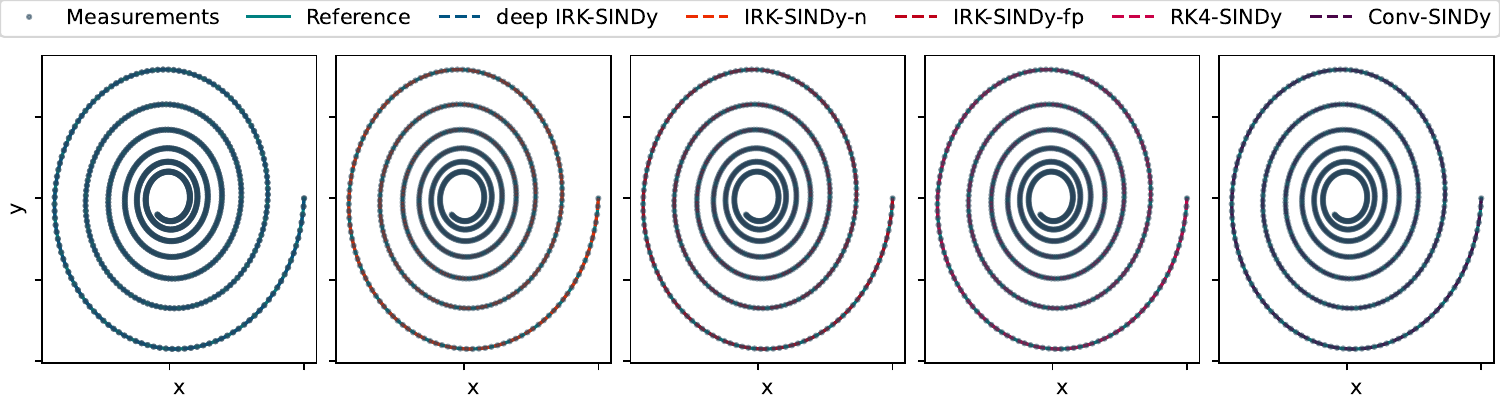}
\end{subfigure}
\begin{subfigure}{\textwidth}
\caption{$m = 201$}
\includegraphics[width=0.9\textwidth, height=3.0cm]{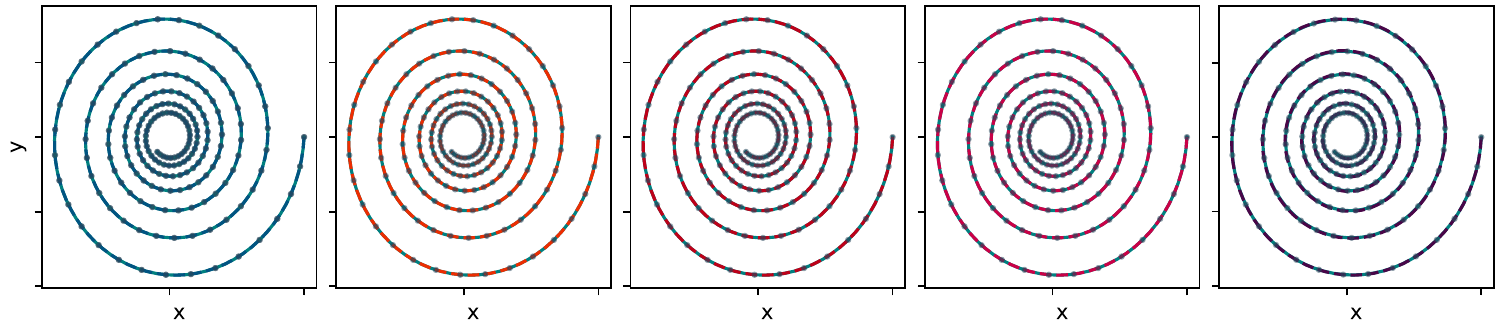}
\end{subfigure}
\begin{subfigure}{\textwidth}
\caption{$m = 51$}
\includegraphics[width=0.9\textwidth, height=3.0cm]{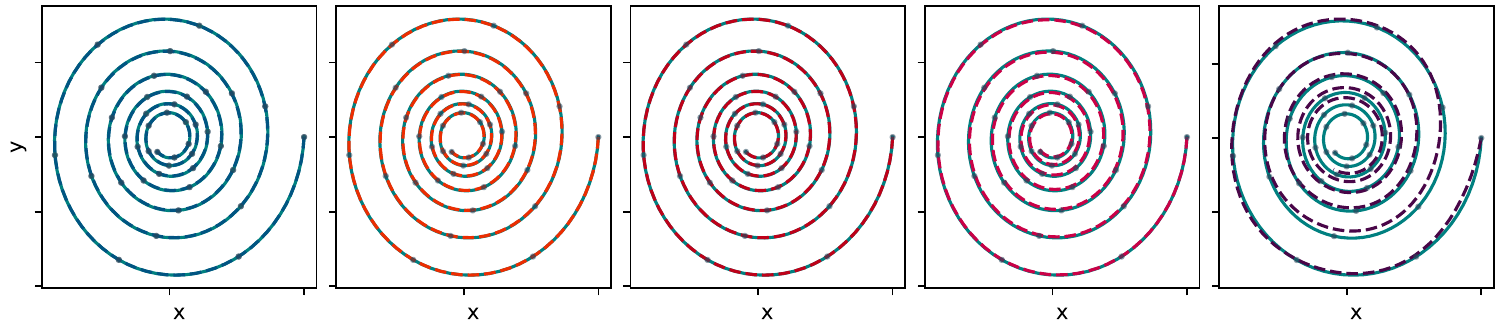}
\end{subfigure}
\begin{subfigure}{\textwidth}
\caption{$m = 31$}
\includegraphics[width=0.9\textwidth, height=3.0cm]{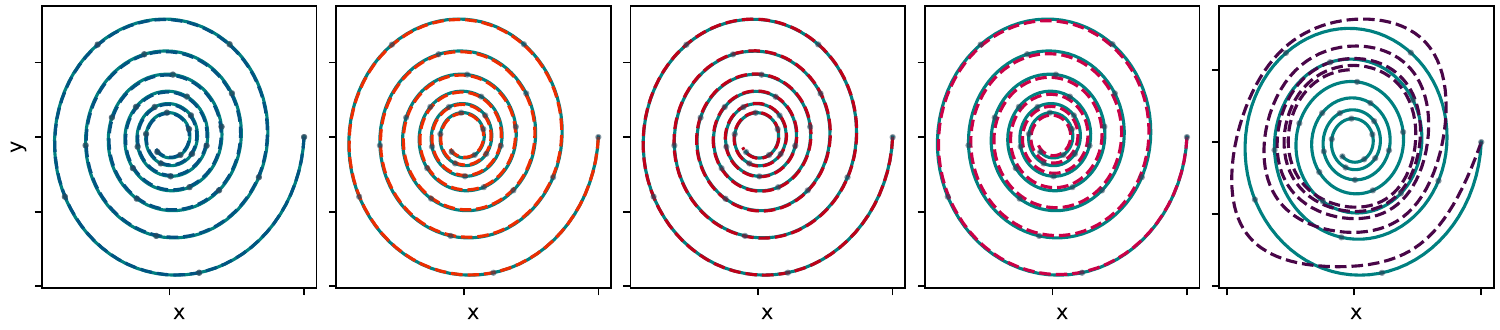}
\end{subfigure}
\caption{Linear damped oscillator: Comparing identified models under various levels of data scarcity with reference model. a. Data, b. sample size $m=801$, c. sample size $m=201$, d. sample size $m=51$, e. sample size $m=31$.}\label{fig3}
\end{figure}

To initiate the data collection process, we establish the initial condition as $\begin{bmatrix}
x_{1}(0) & x_{2}(0)
\end{bmatrix}^{T} = \begin{bmatrix}
2.0 & 0.0
\end{bmatrix}^{T}$, and subsequently sample a total of $m+1$ data points utilizing a fixed stepsize $20/m$, within the time interval $t\in [0, 20]$. In the first scenario, we use a variety of values for $m$, with the objective of rigorously assessing the robustness of the proposed methodologies to data scarcity, without adding noise to measurements. It is noteworthy to mention that the algorithm is also flexible for data with variable stepsize, however, for the sake of simplicity and clarity in our analytical framework, we have opted to utilize a constant stepsize in the current analysis. The data from the sampling process are plotted in Figure\ref{fig3a}.

\begin{table}[ht]
\centering
\begin{tabular}{llll}
\hline
\multicolumn{1}{|l|}{Sample-size} & \multicolumn{1}{l|}{Deep IRK-SINDy} & \multicolumn{1}{l|}{RK4-SINDy} & \multicolumn{1}{l|}{Conv-SINDy} \\ \hline
\multicolumn{1}{|l|}{$m = 801$} & \multicolumn{1}{l|}{
\begin{minipage}{4.75cm}{\begin{equation*}
\begin{split}
\dot{x}_{1}(t) = -0.100 x_{1}(t) + 2.000 x_{2}(t)\\
\dot{x}_{2}(t) = -2.001 x_{1}(t) - 0.100 x_{2}(t)
\end{split}
\end{equation*}}\end{minipage}} & \multicolumn{1}{l|}{
\begin{minipage}{4.75cm}{\begin{equation*}
\begin{split}
\dot{x}_{1}(t) = -0.100 x_{1}(t) + 2.000 x_{2}(t)\\
\dot{x}_{2}(t) = -2.001 x_{1}(t) - 0.100 x_{2}(t)
\end{split}
\end{equation*}}\end{minipage}} & \multicolumn{1}{l|}{
\begin{minipage}{4.75cm}{\begin{equation*}
\begin{split}
\dot{x}_{1}(t) = -0.100 x_{1}(t) + 2.000 x_{2}(t)\\
\dot{x}_{2}(t) = -2.000 x_{1}(t) - 0.100 x_{2}(t)
\end{split}
\end{equation*}}\end{minipage}} \\ \hline
\multicolumn{1}{|l|}{$m = 201$} & \multicolumn{1}{l|}{
\begin{minipage}{4.75cm}{\begin{equation*}
\begin{split}
\dot{x}_{1}(t) = -0.100 x_{1}(t) + 2.000 x_{2}(t)\\
\dot{x}_{2}(t) = -2.001 x_{1}(t) - 0.100 x_{2}(t)
\end{split}
\end{equation*}}\end{minipage}} & \multicolumn{1}{l|}{
\begin{minipage}{4.75cm}{\begin{equation*}
\begin{split}
\dot{x}_{1}(t) = -0.100 x_{1}(t) + 2.001 x_{2}(t)\\
\dot{x}_{2}(t) = -2.001 x_{1}(t) - 0.100 x_{2}(t)
\end{split}
\end{equation*}}\end{minipage}} & \multicolumn{1}{l|}{
\begin{minipage}{4.75cm}{\begin{equation*}
\begin{split}
\dot{x}_{1}(t) = -0.099 x_{1}(t) + 1.987 x_{2}(t)\\
\dot{x}_{2}(t) = -1.989 x_{1}(t) - 0.098 x_{2}(t)
\end{split}
\end{equation*}}\end{minipage}} \\ \hline
\multicolumn{1}{|l|}{$m = 41$} & \multicolumn{1}{l|}{
\begin{minipage}{4.75cm}{\begin{equation*}
\begin{split}
\dot{x}_{1}(t) = -0.101 x_{1}(t) + 2.003 x_{2}(t)\\
\dot{x}_{2}(t) = -2.003 x_{1}(t) - 0.101 x_{2}(t)
\end{split}
\end{equation*}}\end{minipage}} & \multicolumn{1}{l|}{
\begin{minipage}{4.75cm}{\begin{equation*}
\begin{split}
\dot{x}_{1}(t) = -0.103 x_{1}(t) + 2.011 x_{2}(t)\\
\dot{x}_{2}(t) = -2.011 x_{1}(t) - 0.103 x_{2}(t)
\end{split}
\end{equation*}}\end{minipage}} & \multicolumn{1}{l|}{
\begin{minipage}{4.75cm}{\begin{equation*}
\begin{split}
&\dot{x}_{1}(t) = -0.066 x_{1}^2(t) - 0.083 x_{1}^{3}(t)\\
& + 1.680 x_{2}(t)\\
&\dot{x}_{2}(t) =   -1.545 x_{1}(t) - 0.080 x_{1}^{2}(t)\\
&- 0.140 x_{1}^{3}(t) + 0.063 x_{1}^{2}(t)x_{2}(t)\\
&- 0.096 x_{2}(t)
\end{split}
\end{equation*}}\end{minipage}} \\ \hline
\multicolumn{1}{|l|}{$m = 31$} & \multicolumn{1}{l|}{
\begin{minipage}{4.75cm}{\begin{equation*}
\begin{split}
\dot{x}_{1}(t) = -0.102 x_{1}(t) + 2.009 x_{2}(t)\\
\dot{x}_{2}(t) = -2.008 x_{1}(t) - 0.103 x_{2}(t)
\end{split}
\end{equation*}}\end{minipage}} & \multicolumn{1}{l|}{
\begin{minipage}{4.75cm}{\begin{equation*}
\begin{split}
\dot{x}_{1}(t) = -0.107 x_{1}(t) + 2.017 x_{2}(t)\\
\dot{x}_{2}(t) = -2.016 x_{1}(t) - 0.106 x_{2}(t)
\end{split}
\end{equation*}}\end{minipage}} & \multicolumn{1}{l|}{
\begin{minipage}{4.75cm}{\begin{equation*}
\begin{split}
&\dot{x}_{1}(t) = -0.158 x_{1}(t) -0.099 x_{1}^{2}(t)\\
& -0.225 x_{1}^{3}(t) -0.125 x_{1}(t)x_{2}(t)\\
&+ 0.103 x_{1}(t)x_{2}^{2}(t) + 1.424 x_{2}(t)\\
&+ 0.050 x_{2}^{3}(t)\\
&\dot{x}_{2}(t) = -1.309 x_{1}(t) -0.098 x_{1}^{2}(t)\\
&- 0.213 x_{1}^{3}(t) -0.126 x_{1}(t)x_{2}(t)\\
&+ 0.121 x_{1}(t)x_{2}^{2}(t) - 0.084 x_{2}(t)\\
&+ 0.062 x_{2}^{3}(t)
\end{split}
\end{equation*}}\end{minipage}} \\ \hline
\end{tabular}
\caption{Linear damped oscillator: the discovered governing equations using deep IRK-SINDy, RK4-SINDy, and Conv-SINDy for various sample-size $m$.}\label{tbl1}
\end{table}

We systematically explore the desired model in the model space of possible descriptions of the dynamical system under investigation, which, in the specific case of this particular example, is restricted to the space of polynomials up to degree $3$. Upon the careful selection of the nonlinear feature library, we proceed to set the threshold value across all approaches to constant value $\lambda=0.05$. Within the frameworks of the IRK-SINDy and RK4-SINDy approaches, we use a learning rate of $lr=0.01$ for updating coefficient matrix $\xi$ and regularly conduct sequential thresholding every $1000$ epochs. It is crucial to point out that in the context of the deep IRK-SINDy approach, the first approximation of stage values may require potentially more epochs initially to facilitate effective sequential thresholding. This requirement arises due to the need of the DNN to undergo adequate training in order to accurately predict the stage values corresponding to IRKs, which subsequently allows for the prediction of both next and preceding step values via eq.\eqref{eqequival}. Thus, it follows that the more rapidly the network acquires proficiency in learning the $\chi^{\theta}$ values, the fewer epochs will be necessitated by the algorithm for the procedure of sequential thresholding. It is imperative to emphasize again that within the framework of the deep IRK-SINDy approach, the Adam optimizer updates $\xi$ and $\theta$ simultaneously, while, the endeavor to learn the stage values imposes a limit on the permissible learning rate. For this reason, in the deep IRK-SINDy framework, we set two different learning rates of $10^{-3}$ and $0.01$ for the DNN parameters and the coefficient matrix $\xi$, respectively. In the first iteration, $15,000$ epochs and in subsequent iterations, $1,000$ epochs are employed to train the network. In this configuration, we use $4$ hidden layers, each comprising $32$ neurons, utilizing the $\tanh$ activation function within the architecture of the fully connected DNN. At the end of each iteration, the learning rate values are updated by multiplying them by a number between $0$ and $1$. We use $4$ and $3$ iterations in fixed point and Newton's iterations, respectly.

\begin{table}[ht]
\centering
\begin{tabular}{lll}
\hline
\multicolumn{1}{|l|}{Sample-size} & \multicolumn{1}{l|}{Newton's iterations} & \multicolumn{1}{l|}{Fixed point iterations} \\ \hline
\multicolumn{1}{|l|}{$m = 801$} & \multicolumn{1}{l|}{
\begin{minipage}{5.5cm}{\begin{equation*}
\begin{split}
\dot{x}_{1}(t) = -0.100 x_{1}(t) + 2.000 x_{2}(t)\\
\dot{x}_{2}(t) = -2.001 x_{1}(t) - 0.100 x_{2}(t)
\end{split}
\end{equation*}}\end{minipage}} & \multicolumn{1}{l|}{
\begin{minipage}{5.5cm}{\begin{equation*}
\begin{split}
\dot{x}_{1}(t) = -0.100 x_{1}(t) + 2.000 x_{2}(t)\\
\dot{x}_{2}(t) = -2.001 x_{1}(t) - 0.100 x_{2}(t)
\end{split}
\end{equation*}}\end{minipage}} \\ \hline
\multicolumn{1}{|l|}{$m = 201$} & \multicolumn{1}{l|}{
\begin{minipage}{5.5cm}{\begin{equation*}
\begin{split}
\dot{x}_{1}(t) = -0.100 x_{1}(t) + 2.000 x_{2}(t)\\
\dot{x}_{2}(t) = -2.001 x_{1}(t) - 0.100 x_{2}(t)
\end{split}
\end{equation*}}\end{minipage}} & \multicolumn{1}{l|}{
\begin{minipage}{5.5cm}{\begin{equation*}
\begin{split}
\dot{x}_{1}(t) = -0.100 x_{1}(t) + 2.000 x_{2}(t)\\
\dot{x}_{2}(t) = -2.001 x_{1}(t) - 0.100 x_{2}(t)
\end{split}
\end{equation*}}\end{minipage}} \\ \hline
\multicolumn{1}{|l|}{$m = 41$} & \multicolumn{1}{l|}{
\begin{minipage}{5.5cm}{\begin{equation*}
\begin{split}
\dot{x}_{1}(t) = -0.101 x_{1}(t) + 2.003 x_{2}(t)\\
\dot{x}_{2}(t) = -2.003 x_{1}(t) - 0.101 x_{2}(t)
\end{split}
\end{equation*}}\end{minipage}} & \multicolumn{1}{l|}{
\begin{minipage}{5.5cm}{\begin{equation*}
\begin{split}
&\dot{x}_{1}(t) = -0.101 x_{1}(t) + 2.004 x_{2}(t)\\
&\dot{x}_{2}(t) = -2.004 x_{1}(t) - 0.101 x_{2}(t)
\end{split}
\end{equation*}}\end{minipage}} \\ \hline
\multicolumn{1}{|l|}{$m = 31$} & \multicolumn{1}{l|}{
\begin{minipage}{5.5cm}{\begin{equation*}
\begin{split}
\dot{x}_{1}(t) = -0.102 x_{1}(t) + 2.009 x_{2}(t)\\
\dot{x}_{2}(t) = -2.008 x_{1}(t) - 0.103 x_{2}(t)
\end{split}
\end{equation*}}\end{minipage}} & \multicolumn{1}{l|}{
\begin{minipage}{5.5cm}{\begin{equation*}
\begin{split}
\dot{x}_{1}(t) = -0.100 x_{1}(t) + 2.015 x_{2}(t)\\
\dot{x}_{2}(t) = -2.014 x_{1}(t) - 0.100 x_{2}(t)
\end{split}
\end{equation*}}\end{minipage}} \\ \hline                  
\end{tabular}
\caption{Linear damped oscillator: the discovered governing equations using IRK-SINDy in the approach of Newton's iterations and fixed point iterations for various sample-size $m$.}\label{tbl2}
\end{table}

\begin{figure}[]
\centering
\begin{subfigure}{\textwidth}
\caption{Collected noisy data}
\includegraphics[width=0.9\textwidth]{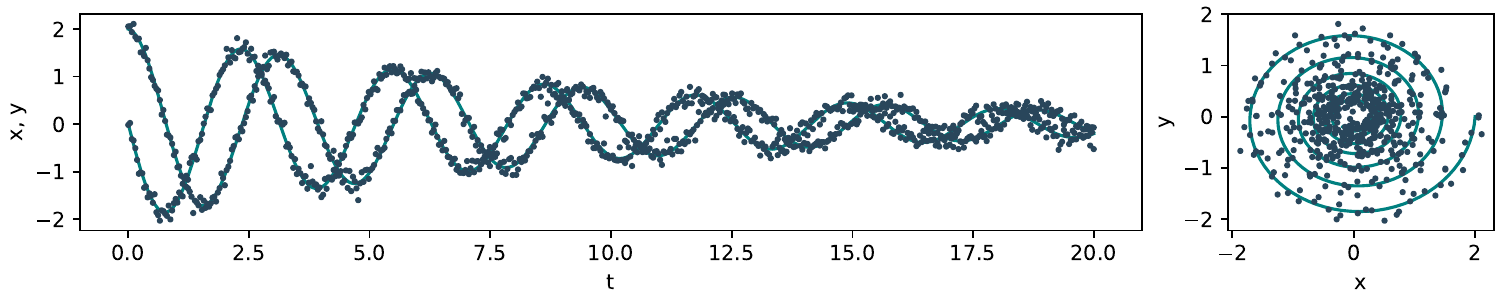}
\end{subfigure}
\begin{subfigure}{\textwidth}
\caption{$\sigma=0.01$}
\includegraphics[width=0.9\textwidth, height=3.4cm]{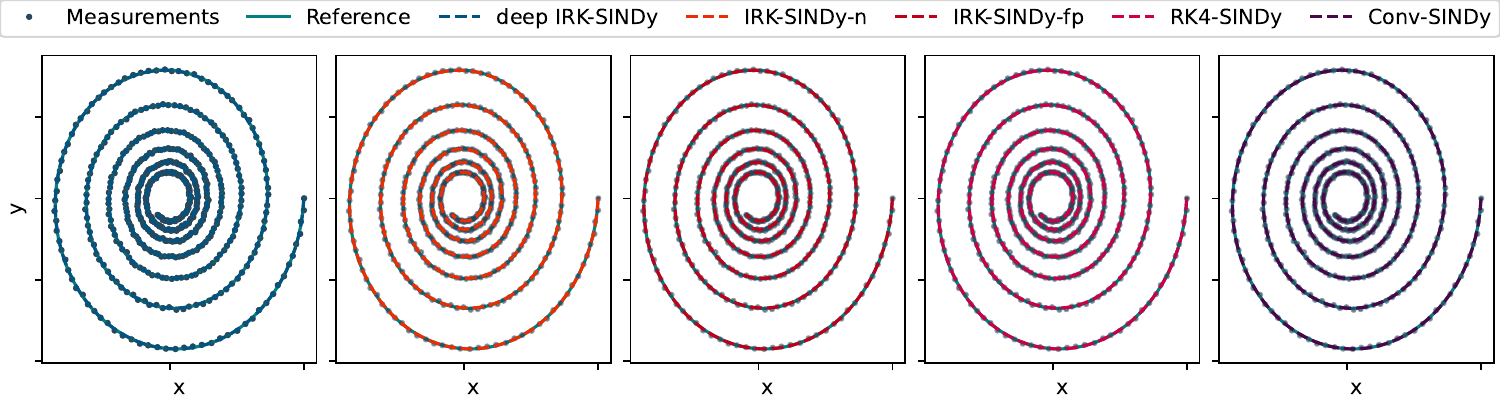}
\end{subfigure}
\begin{subfigure}{\textwidth}
\caption{$\sigma=0.04$}
\includegraphics[width=0.9\textwidth, height=3.0cm]{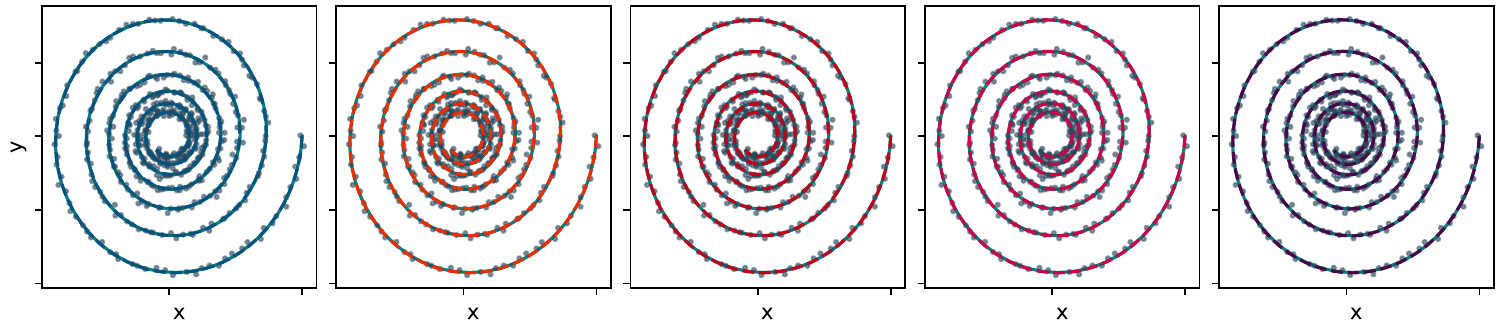}
\end{subfigure}
\begin{subfigure}{\textwidth}
\caption{$\sigma=0.08$}
\includegraphics[width=0.9\textwidth, height=3.0cm]{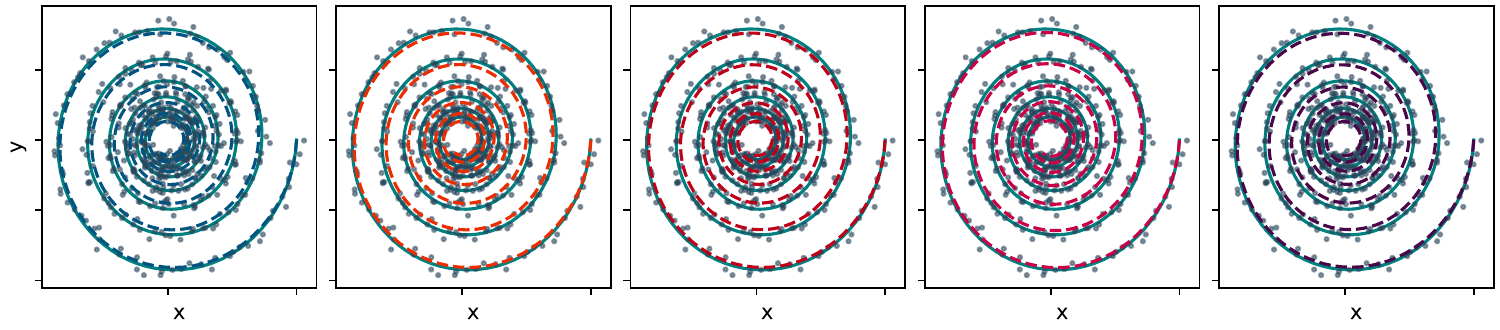}
\end{subfigure}
\begin{subfigure}{\textwidth}
\caption{$\sigma=0.16$}
\includegraphics[width=0.9\textwidth, height=3.0cm]{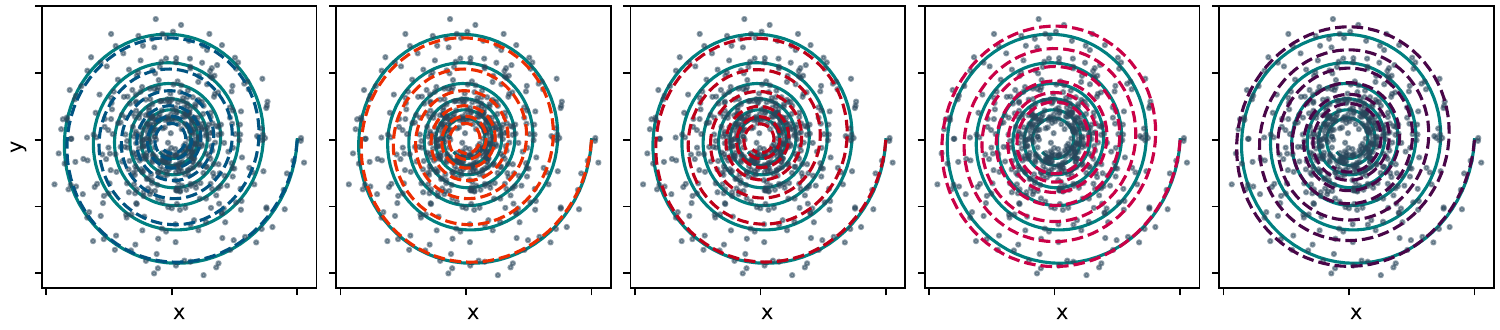}
\end{subfigure}
\caption{Linear damped oscillator: Comparing the response of identified models under various noise levels in measurements with reference model. a. Noisy data, b. noise level $\sigma = 0.01$, c. nise level $\sigma = 0.04$, d. noise level $\sigma=0.08$, e. noise level $\sigma=0.16$.}\label{fig4}
\end{figure}

Figure\ref{fig3} depicts a qualitative evaluation of the accuracy associated with the discovered dynamical systems by providing a comparative analysis between the reference trajectories and the predicted trajectories, alongside the resultant phase portrait. It reveals that all approaches, including the two versions of IRK-SINDy, are able to capture the dynamical evolution of the system given sufficient and clean data. As delineated in Table\ref{tbl1} and Table\ref{tbl2}, the equations derived from these approaches demonstrate a remarkable consistency with the reference equation. Furthermore, Figure\ref{fig3} confirms that for a limited amount of data, the two innovative proposed approaches exhibit superior performance relative to the previously mentioned approaches. For instance, in the scenario where $m = 31$, both the IRK-SINDy and deep IRK-SINDy approaches successfully capture the governing equations, outperforming the other approaches. Additionally, as anticipated, Figure\ref{fig3} and Table\ref{tbl2} show that applying IRK-SINDy with Newton’s itiration yields superior results when compared to its utilization with fixed-point iteration. Table\ref{tbl1} indicates that the two approaches Conv-SINDy and RK4-SINDy have encountered challenges in identifying correct and sparse models, particularly under conditions of significant data scarcity.

\begin{table}[ht]
\centering
\begin{tabular}{llll}
\hline
\multicolumn{1}{|l|}{Noise} & \multicolumn{1}{l|}{Deep IRK-SINDy} & \multicolumn{1}{l|}{RK4-SINDy} & \multicolumn{1}{l|}{Conv-SINDy} \\ \hline
\multicolumn{1}{|l|}{$\sigma = 0.01$} & \multicolumn{1}{l|}{
\begin{minipage}{4.6cm}{\begin{equation*}
\begin{split}
&\dot{x}_{1}(t) = -0.103 x_{1}(t) + 2.000 x_{2}(t)\\
&\dot{x}_{2}(t) = -2.001 x_{1}(t) - 0.099 x_{2}(t)
\end{split}
\end{equation*}}\end{minipage}} & \multicolumn{1}{l|}{
\begin{minipage}{4.6cm}{\begin{equation*}
\begin{split}
&\dot{x}_{1}(t) = -0.102 x_{1}(t) + 2.000 x_{2}(t)\\
&\dot{x}_{2}(t) = -2.001 x_{1}(t) - 0.099 x_{2}(t)
\end{split}
\end{equation*}}\end{minipage}} & \multicolumn{1}{l|}{
\begin{minipage}{4.6cm}{\begin{equation*}
\begin{split}
&\dot{x}_{1}(t) = -0.102 x_{1}(t) + 1.999 x_{2}(t)\\
&\dot{x}_{2}(t) = -1.999 x_{1}(t) - 0.099 x_{2}(t)
\end{split}
\end{equation*}}\end{minipage}} \\ \hline
\multicolumn{1}{|l|}{$\sigma = 0.04$} & \multicolumn{1}{l|}{
\begin{minipage}{4.6cm}{\begin{equation*}
\begin{split}
&\dot{x}_{1}(t) = -0.105 x_{1}(t) + 1.995 x_{2}(t)\\
&\dot{x}_{2}(t) = -2.009 x_{1}(t) - 0.097 x_{2}(t)
\end{split}
\end{equation*}}\end{minipage}} & \multicolumn{1}{l|}{
\begin{minipage}{4.6cm}{\begin{equation*}
\begin{split}
&\dot{x}_{1}(t) = -0.104 x_{1}(t) + 1.995 x_{2}(t)\\
&\dot{x}_{2}(t) = -2.008 x_{1}(t) - 0.097 x_{2}(t)
\end{split}
\end{equation*}}\end{minipage}} & \multicolumn{1}{l|}{
\begin{minipage}{4.6cm}{\begin{equation*}
\begin{split}
&\dot{x}_{1}(t) = -0.104 x_{1}(t) + 1.993 x_{2}(t)\\
&\dot{x}_{2}(t) = -2.008 x_{1}(t) - 0.097 x_{2}(t)
\end{split}
\end{equation*}}\end{minipage}} \\ \hline
\multicolumn{1}{|l|}{$\sigma = 0.08$} & \multicolumn{1}{l|}{
\begin{minipage}{4.6cm}{\begin{equation*}
\begin{split}
&\dot{x}_{1}(t) = -0.117 x_{1}(t) + 2.013 x_{2}(t)\\
&\dot{x}_{2}(t) = -1.972 x_{1}(t) - 0.105 x_{2}(t)
\end{split}
\end{equation*}}\end{minipage}} & \multicolumn{1}{l|}{
\begin{minipage}{4.6cm}{\begin{equation*}
\begin{split}
&\dot{x}_{1}(t) = -0.115 x_{1}(t) + 2.013 x_{2}(t)\\
&\dot{x}_{2}(t) = -1.972 x_{1}(t) - 0.103 x_{2}(t)
\end{split}
\end{equation*}}\end{minipage}} & \multicolumn{1}{l|}{
\begin{minipage}{4.6cm}{\begin{equation*}
\begin{split}
&\dot{x}_{1}(t) = -0.120 x_{1}(t) + 2.011 x_{2}(t)\\
&\dot{x}_{2}(t) = -1.964 x_{1}(t) - 0.102 x_{2}(t)
\end{split}
\end{equation*}}\end{minipage}} \\ \hline
\multicolumn{1}{|l|}{$\sigma = 0.16$} & \multicolumn{1}{l|}{
\begin{minipage}{4.6cm}{\begin{equation*}
\begin{split}
&\dot{x}_{1}(t) = -0.147 x_{1}(t) + 1.991 x_{2}(t)\\
&\dot{x}_{2}(t) = -2.006 x_{1}(t) - 0.086 x_{2}(t)
\end{split}
\end{equation*}}\end{minipage}} & \multicolumn{1}{l|}{
\begin{minipage}{4.6cm}{\begin{equation*}
\begin{split}
&\dot{x}_{1}(t) = -0.139 x_{1}(t) + 1.991 x_{2}(t)\\
&\dot{x}_{2}(t) = -2.002 x_{1}(t)
\end{split}
\end{equation*}}\end{minipage}} & \multicolumn{1}{l|}{
\begin{minipage}{4.6cm}{\begin{equation*}
\begin{split}
&\dot{x}_{1}(t) = -0.144 x_{1}(t) + 1.988 x_{2}(t)\\
&\dot{x}_{2}(t) = -1.999 x_{1}(t)
\end{split}
\end{equation*}}\end{minipage}} \\ \hline
\end{tabular}
\caption{Linear damped oscillator: the discovered governing equations using deep IRK-SINDy, RK4-SINDy, and Conv-SINDy for various noise level $\sigma$.}\label{tbl3}
\end{table}

In the second scenario, $m$ is kept constant and equal to $551$, while various levels of noise are introduced to the measurements in order to examine the robustness of the proposed methodologies to noise. Herein, we utilize $\sigma \in \{0.01, 0.04, 0.08, 0.16\}$ to produce Gaussian noise $\mathcal{N}(\mu, \sigma^{2})$ with a zero mean $\mu = 0$ and variance $\sigma^{2}$, whereby the noise level is controlled by the standard deviation $sigma$. We employ $3$ thresholding iterations with $2,000$ epochs per iteration, using a thresholding parameter of $\lambda = 0.06$. For deep IRK-SINDy, $20,000$ epochs are allocated in the first iteration and, $2,000$ epochs in the subsequent iterations, employing the identical DNN architecture. Before the training process, we employ the Savitzky–Golay \cite{Savitzky1964} filter for data preprocessing to obtain denoised data \cite{Goyal2022, Naozuka2022}. This preprocessing phase is employed by default in the PySINDy module \cite{Kaptanoglu2021}. For unbiased comparative analysis, we apply the finite difference derivative approximation in the Conv-SINDy simulation step.

\begin{table}[ht]
\centering
\begin{tabular}{lll}
\hline
\multicolumn{1}{|l|}{Noise} & \multicolumn{1}{l|}{Newton's iterations} & \multicolumn{1}{l|}{Fixed point iterations} \\ \hline
\multicolumn{1}{|l|}{$\sigma = 0.01$} & \multicolumn{1}{l|}{
\begin{minipage}{4.6cm}{\begin{equation*}
\begin{split}
&\dot{x}_{1}(t) = -0.102 x_{1}(t) + 2.000 x_{2}(t)\\
&\dot{x}_{2}(t) = -2.001 x_{1}(t) - 0.099 x_{2}(t)
\end{split}
\end{equation*}}\end{minipage}} & \multicolumn{1}{l|}{
\begin{minipage}{4.6cm}{\begin{equation*}
\begin{split}
&\dot{x}_{1}(t) = -0.102 x_{1}(t) + 2.000 x_{2}(t)\\
&\dot{x}_{2}(t) = -2.001 x_{1}(t) - 0.099 x_{2}(t)
\end{split}
\end{equation*}}\end{minipage}} \\ \hline
\multicolumn{1}{|l|}{$\sigma = 0.04$} & \multicolumn{1}{l|}{
\begin{minipage}{4.6cm}{\begin{equation*}
\begin{split}
&\dot{x}_{1}(t) = -0.105 x_{1}(t) + 1.995 x_{2}(t)\\
&\dot{x}_{2}(t) = -2.009 x_{1}(t) - 0.097 x_{2}(t)
\end{split}
\end{equation*}}\end{minipage}} & \multicolumn{1}{l|}{
\begin{minipage}{4.6cm}{\begin{equation*}
\begin{split}
&\dot{x}_{1}(t) = -0.105 x_{1}(t) + 1.995 x_{2}(t)\\
&\dot{x}_{2}(t) = -2.008 x_{1}(t) - 0.097 x_{2}(t)
\end{split}
\end{equation*}}\end{minipage}} \\ \hline
\multicolumn{1}{|l|}{$\sigma = 0.08$} & \multicolumn{1}{l|}{
\begin{minipage}{4.6cm}{\begin{equation*}
\begin{split}
&\dot{x}_{1}(t) = -0.117 x_{1}(t) + 2.013 x_{2}(t)\\
&\dot{x}_{2}(t) = -1.972 x_{1}(t) - 0.105 x_{2}(t)
\end{split}
\end{equation*}}\end{minipage}} & \multicolumn{1}{l|}{
\begin{minipage}{4.6cm}{\begin{equation*}
\begin{split}
&\dot{x}_{1}(t) = -0.117 x_{1}(t) + 2.013 x_{2}(t)\\
&\dot{x}_{2}(t) = -1.972 x_{1}(t) - 0.105 x_{2}(t)
\end{split}
\end{equation*}}\end{minipage}} \\ \hline
\multicolumn{1}{|l|}{$\sigma = 0.16$} & \multicolumn{1}{l|}{
\begin{minipage}{4.6cm}{\begin{equation*}
\begin{split}
&\dot{x}_{1}(t) = -0.147 x_{1}(t) + 1.991 x_{2}(t)\\
&\dot{x}_{2}(t) = -2.006 x_{1}(t) - 0.086 x_{2}(t)
\end{split}
\end{equation*}}\end{minipage}} & \multicolumn{1}{l|}{
\begin{minipage}{4.6cm}{\begin{equation*}
\begin{split}
&\dot{x}_{1}(t) = -0.148 x_{1}(t) + 1.991 x_{2}(t)\\
&\dot{x}_{2}(t) = -2.006 x_{1}(t) - 0.088 x_{2}(t)
\end{split}
\end{equation*}}\end{minipage}} \\ \hline                  
\end{tabular}
\caption{Linear damped oscillator: the discovered governing equations using IRK-SINDy in the approach of Newton's iterations and fixed point iterations for various noise levels $\sigma$.}\label{tbl4}
\end{table}

Our simulations demonstrate that our methodology is robust to noise. As illustrated in Figure\ref{fig4} and Table\ref{tbl3}, our approach is more robust to noise than previous approaches. Table\ref{tbl4} distinctly indicates that the application of Newton’s iterative method in the implementation of IRK-SINDy surpasses fixed point methods concerning noise. Despite the efficiency of the IRK-SINDy framework against noise, taking into account the capabilities of DNNs, appropriate architectures may be developed to enhance the network’s noise resistance, which can be the subject of future studies.

\subsection{Cubic damped oscillator}
Now, let us consider the two-dimensional damped harmonic oscillator characterized by cubic dynamics given in eq.\eqref{nloscilator}:
\begin{subequations} \label{nloscilator}
\begin{equation}
\dot{x}_{1}(t) = -0.1 x_{1}^{3}(t) + 2.0 x_{2}^{3}(t),
\end{equation}
\begin{equation}
\dot{x}_{2}(t) = -2.0 x_{1}^{3}(t) - 0.1 x_{2}^{3}(t).
\end{equation}
\end{subequations}
In this experiment, we employ the initial condition $
\begin{bmatrix}
x_{1}(0) & x_{2}(0)
\end{bmatrix}^{T} =
\begin{bmatrix}
2.0 & 0.0
\end{bmatrix}]^{T}$ to collect data across the temporal interval $t \in [0, 20]$ for the cases where $m$ is assumed the values of $801$, $401$, $101$, and $51$. The objective is to successfully recover the governing equations that dictate nonlinear behavior of the system, utilizing both a sufficient and scarce data. The architecture of DNN comprises a total of four hidden layers, each incorporating $32$ neurons. Furthermore, we establish a thresholding value of $\lambda=0.05$, alongside a learning rate $lr=10^{-3}$ for the parameters associated with the DNN, while concurrently applying a learning rate of $10^{-2}$ for the coefficient matrix $\xi$. Throughout this simulation, we incorporate a total of three thresholding iterations, each consisting of $1,000$ epochs conducted within the polynomial space up to degree $3$. $15,000$ epochs are used for the first iteration in the deep IRK-SINDy training process.

\begin{figure}[]
\flushleft
\begin{subfigure}{0.46\textwidth}
\caption{}
\begin{tikzpicture}[
node distance=0.0cm and 0.0cm,
block/.style={rectangle, draw, rounded corners, fill=teal!0, minimum width=2cm, minimum height=1cm},
>=Stealth,
every node/.style={align=center} 
]

\begin{scope}[local bounding box=NN, rounded corners=5pt, minimum width=2cm, minimum height=0cm]

\node[block, rounded corners=5pt, minimum width=0.99\textwidth, minimum height=7.5cm, dashed] (b1) {};

\node[below right=0.2cm of b1.north west] (img1) {\includegraphics[width=0.45\textwidth]{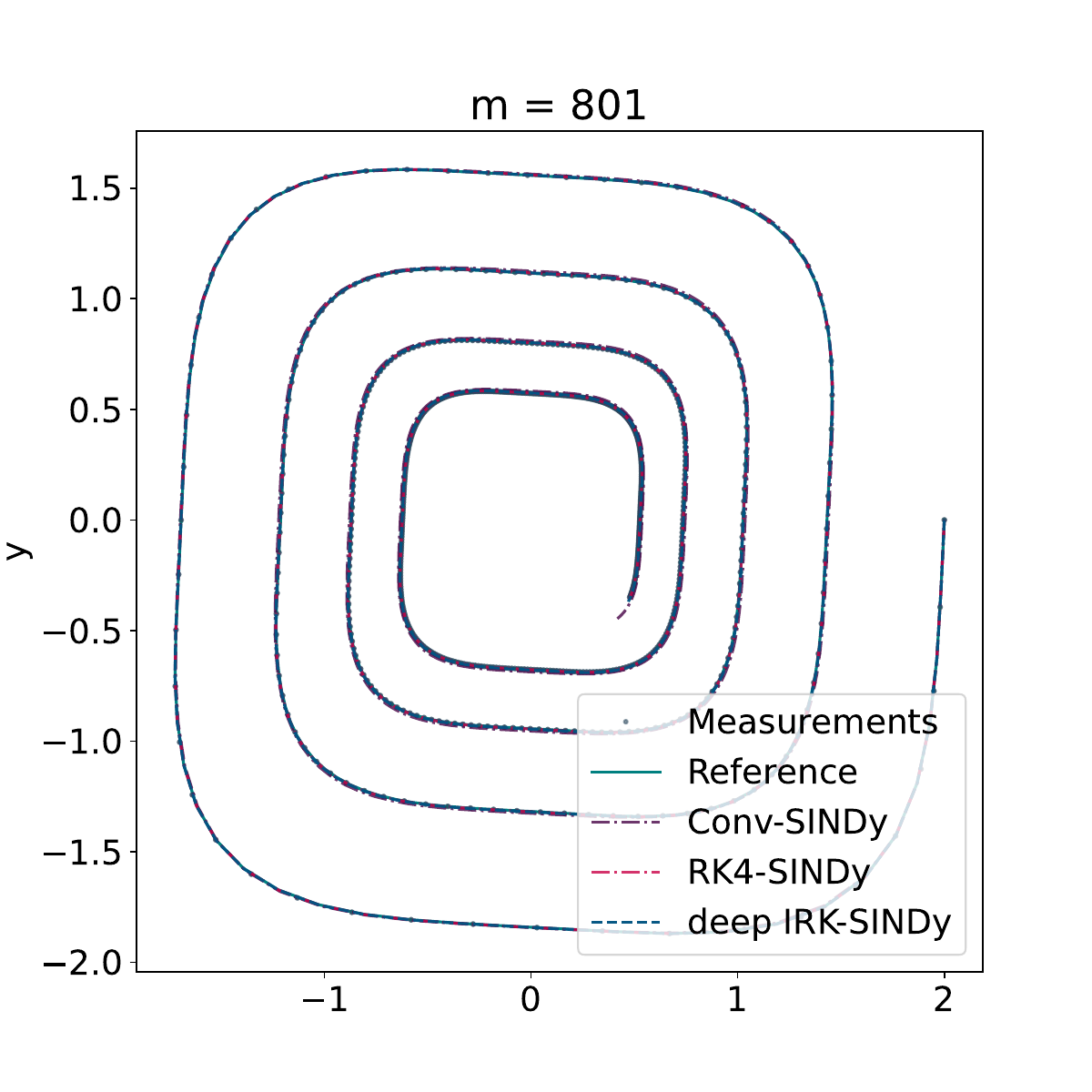}};
\node[below=of img1] (img2) {\includegraphics[width=0.45\textwidth]{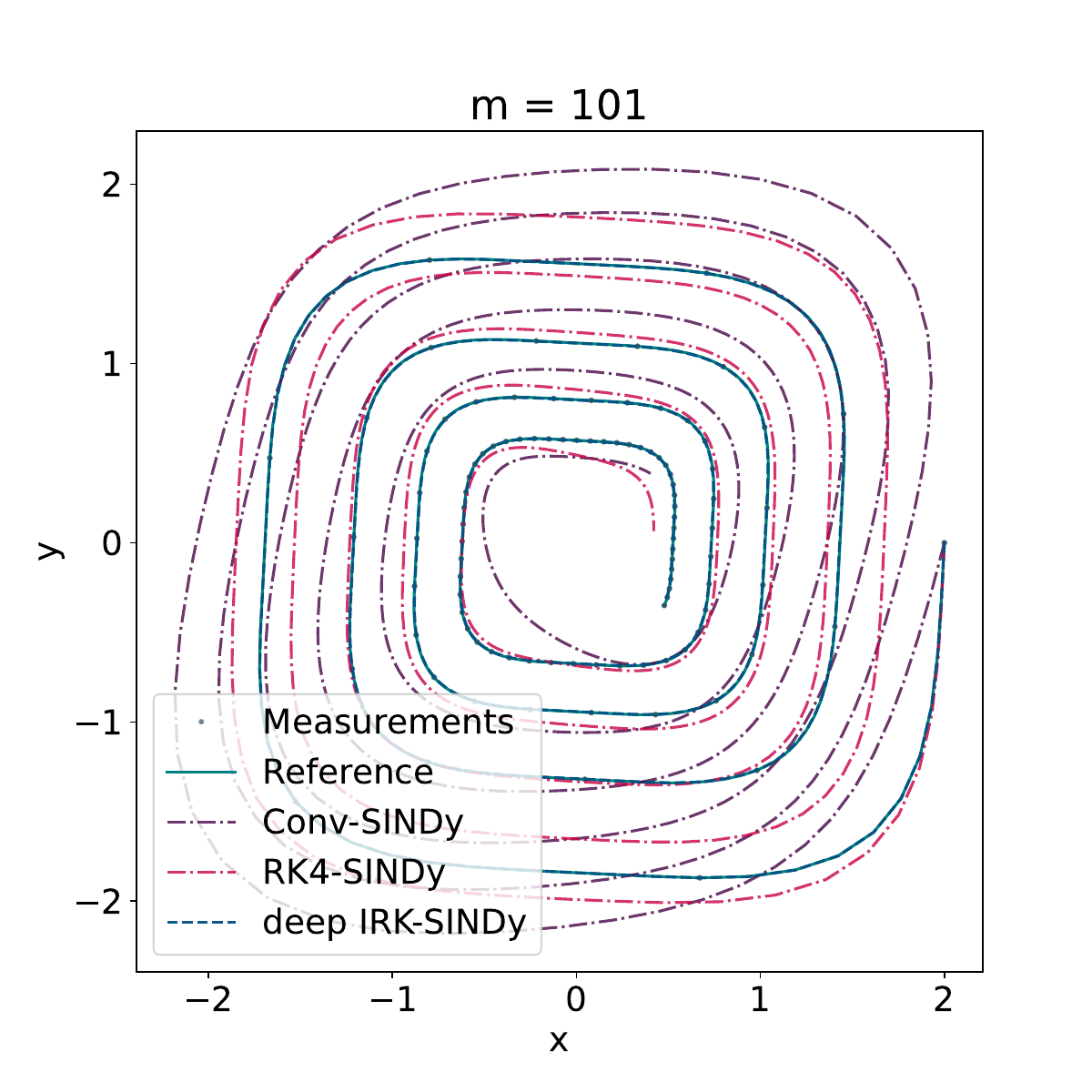}};
\node[right=of img1] (img3) {\includegraphics[width=0.45\textwidth]{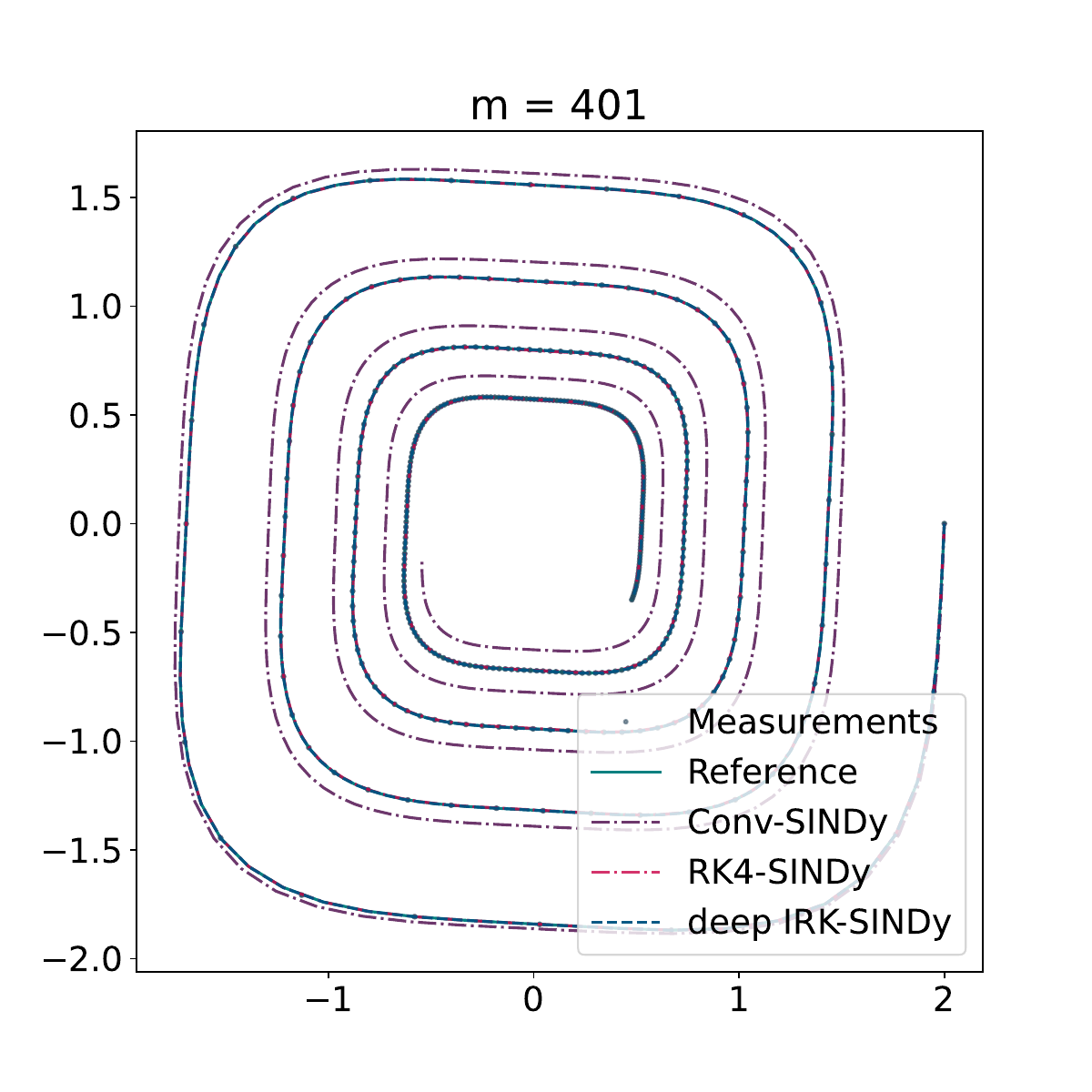}};
\node[below=of img3] (img4) {\includegraphics[width=0.45\textwidth]{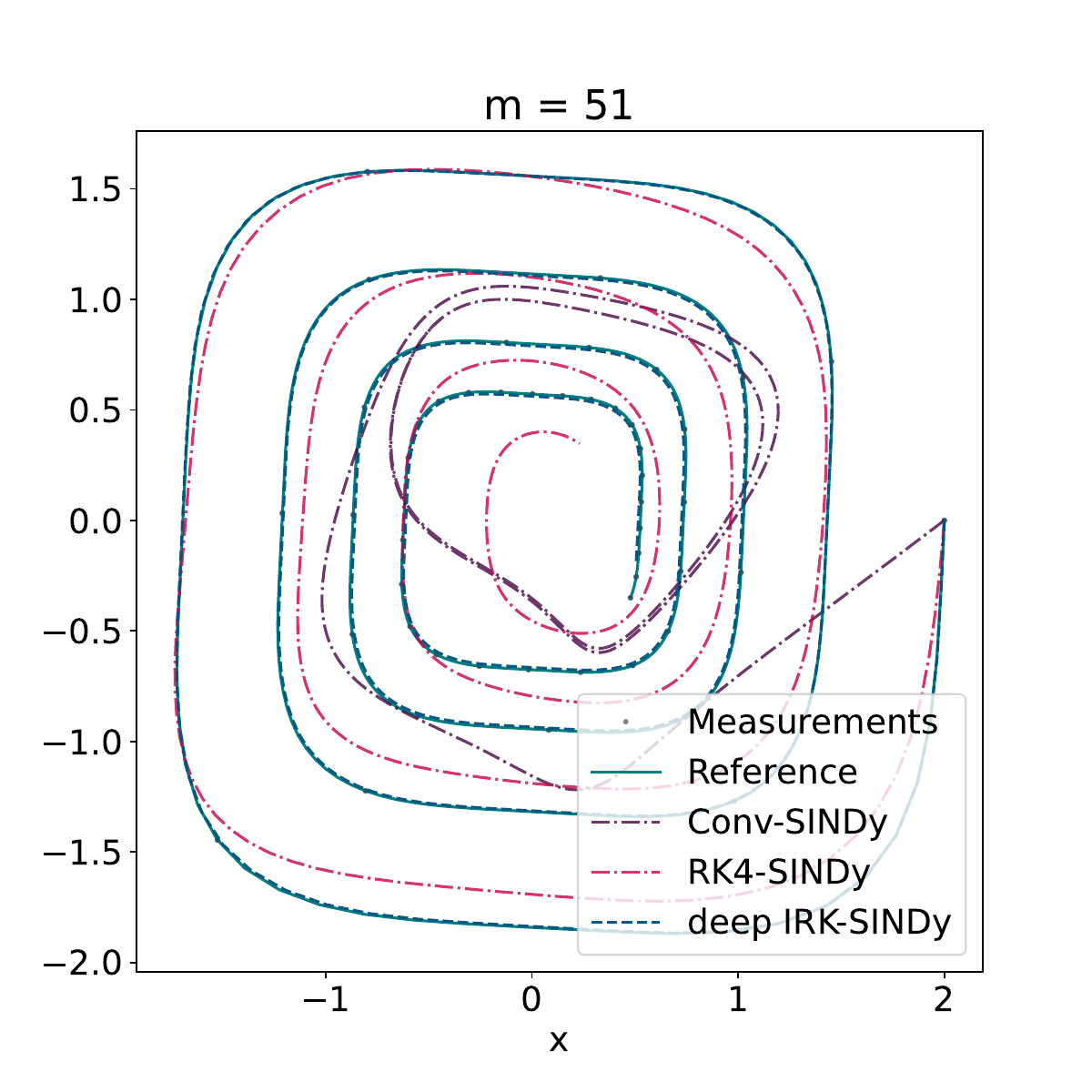}};
\end{scope}
\end{tikzpicture}
\end{subfigure}
\begin{subfigure}{0.47\textwidth}
\caption{}
\begin{tikzpicture}[
node distance=0.0cm and 0.0cm,
neuron/.style={circle, fill=teal!0, draw, minimum size=0.4cm},
block/.style={rectangle, draw, rounded corners, fill=teal!0, minimum width=2cm, minimum height=1cm},
>=Stealth,
every node/.style={align=center} 
]

\begin{scope}[local bounding box=NN, rounded corners=5pt, minimum width=2cm, minimum height=0cm]

\node[block, rounded corners=5pt, minimum width=0.99\textwidth, minimum height=7.5cm, dashed] (b1) {};

\node[below right=0.2cm of b1.north west] (img1) {\includegraphics[width=0.45\textwidth]{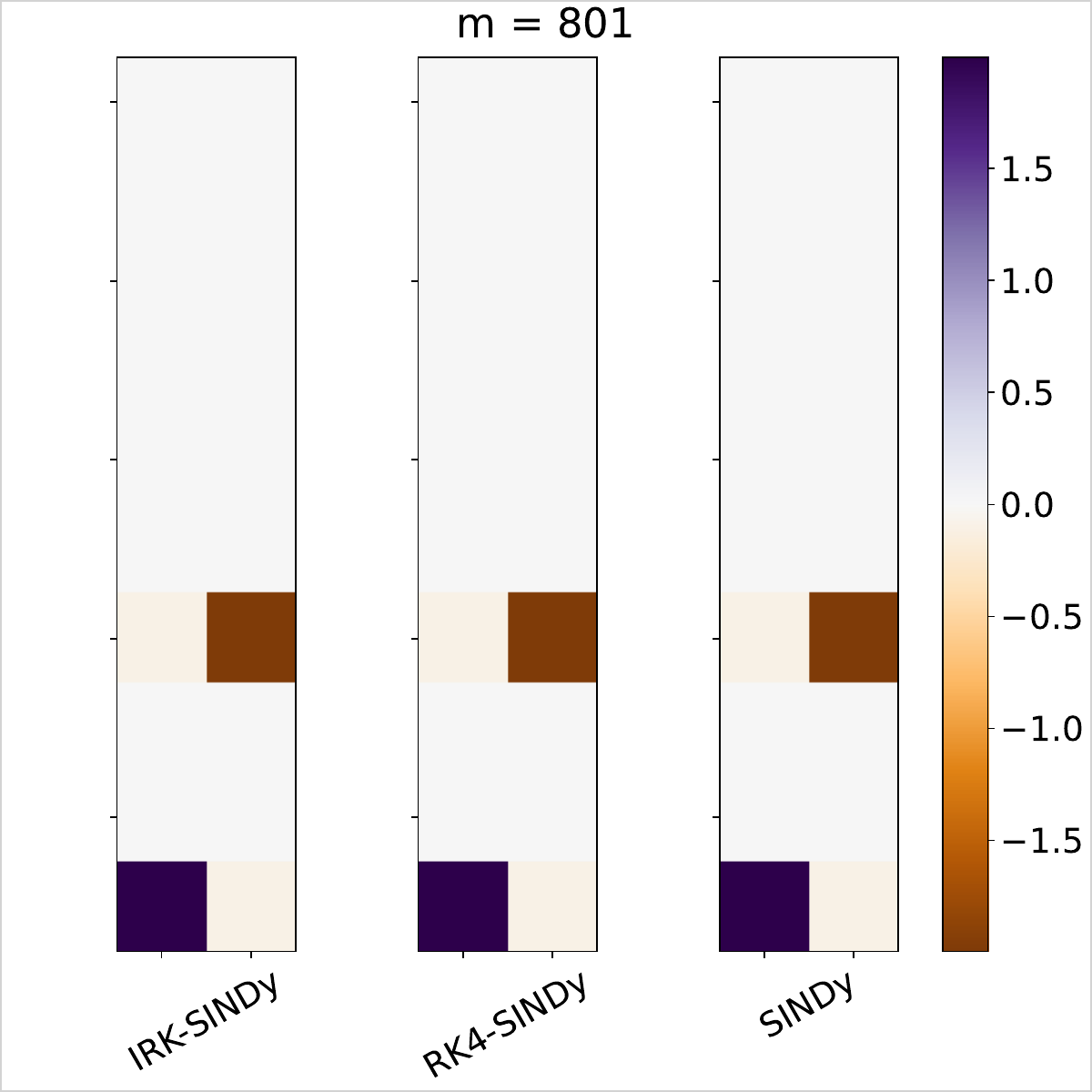}};
\node[below=of img1] (img2) {\includegraphics[width=0.45\textwidth]{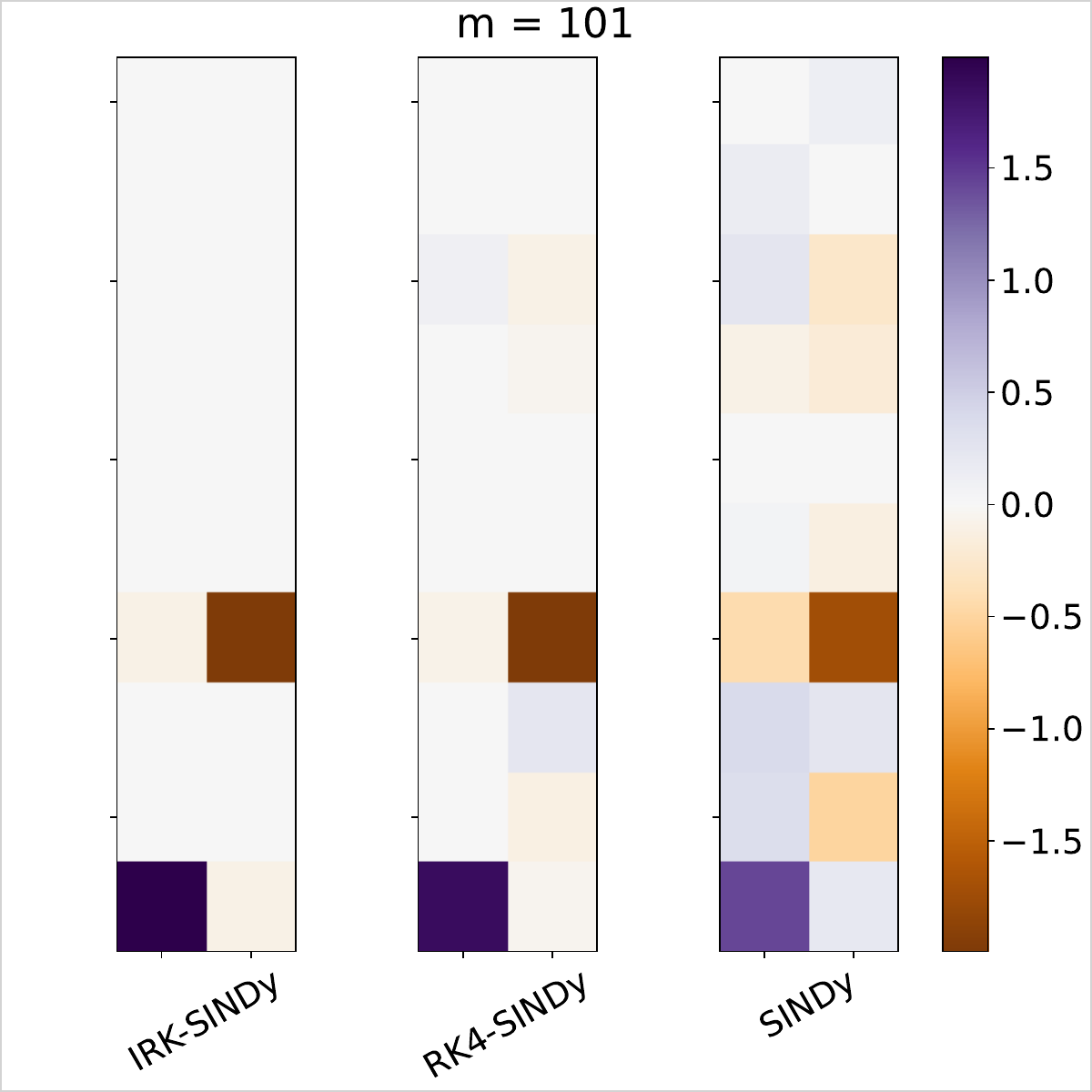}};
\node[right=of img1] (img3) {\includegraphics[width=0.45\textwidth]{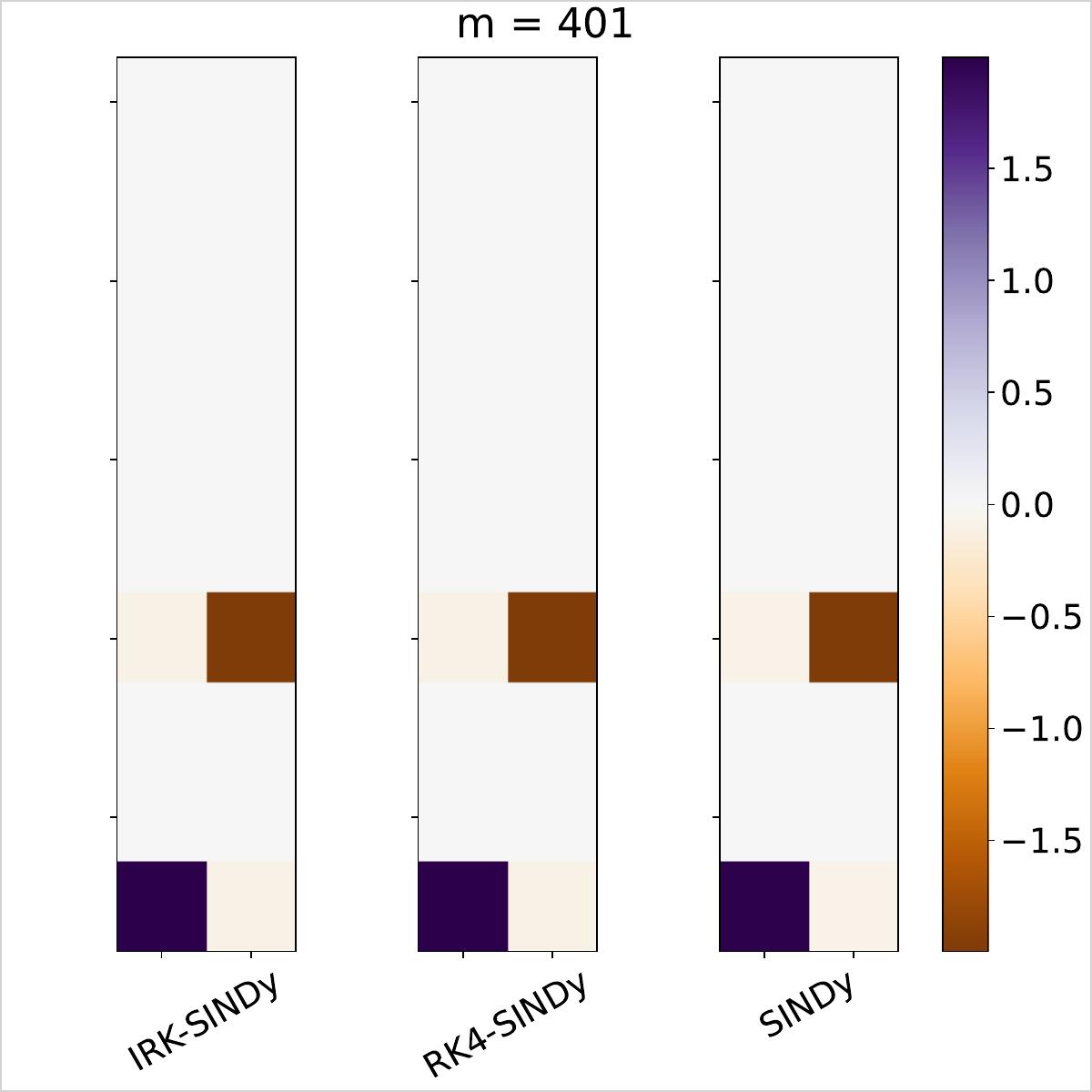}};
\node[below=of img3] (img4) {\includegraphics[width=0.45\textwidth]{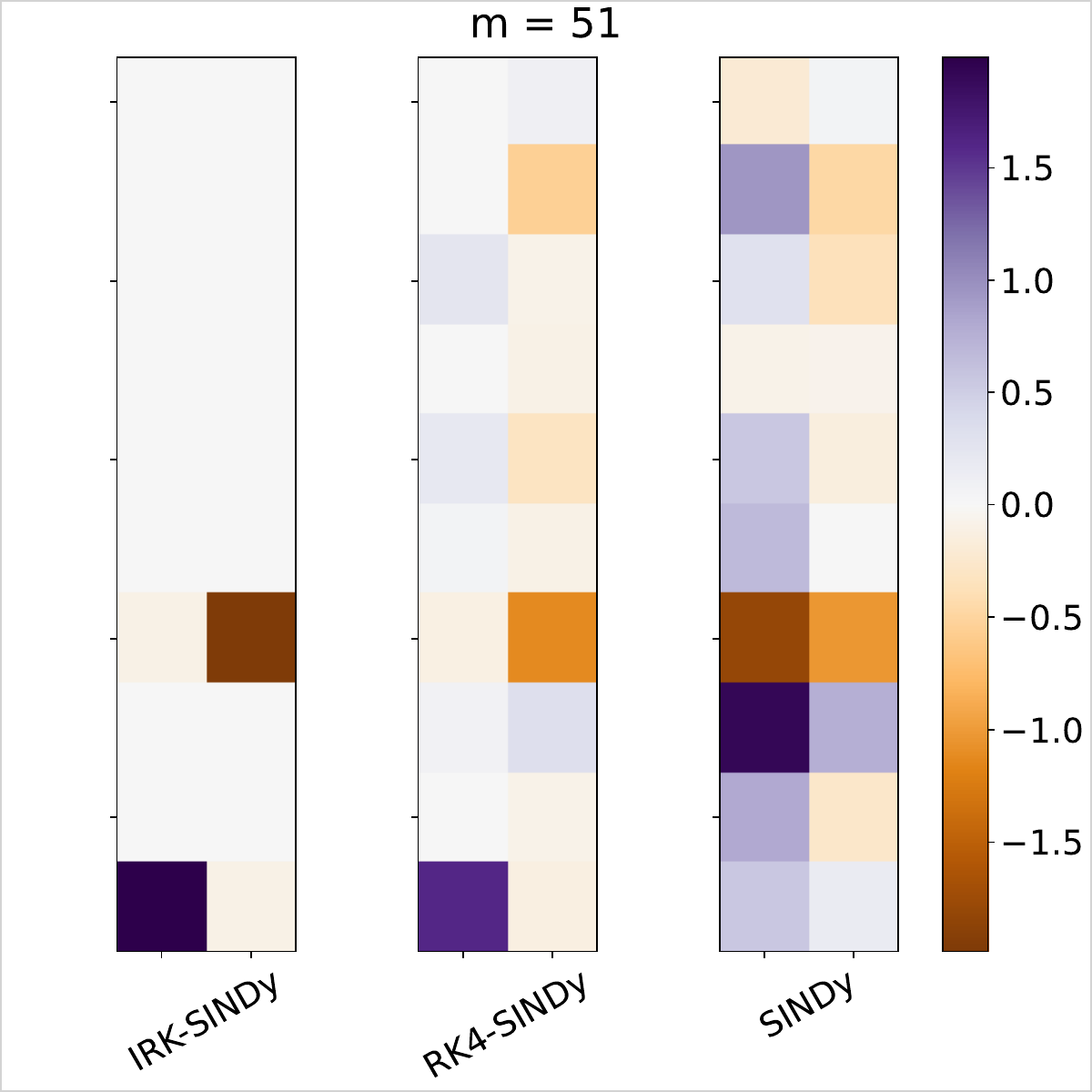}};
\end{scope}
\end{tikzpicture}
\end{subfigure}
\caption{Cubic damped oscillator: Comparing identified models under various levels of data scarcity with reference model. a. Phase portraits, b. coefficient matrices for $m\in\{ 801, 401, 101, 51 \}$. IRK-SINDy provides a more parsimonious and generalizable model compared to RK4-SINDy and Conv-SINDy.}\label{fig5}
\end{figure}

In the subsequent analysis presented in Figure\ref{fig5}, we provide a comprehensive comparison between IRK-SINDy, RK4-SINDy, and Conv-SINDy. It becomes readily apparent that the IRK-SINDy aaproach accurately identify the dynamics of system \eqref{nloscilator}, whereas the Conv-SINDy approach exhibits considerable difficulties when confronted with smaller values of $m$. Moreover, it is noteworthy that IRK-SINDy demonstrates a superior capability in discovering interpretable equations in comparison to RK4-SINDy, particularly in scenarios characterized by scarce data availability. This capability can be attributed to the A-stability properties in Gauss methods \cite{Butcher2016, Wanner1996}, which ensures that the stability region of these methods contains the entire left half-plane of the coordinate system, in contrast to the finite stability region associated with explicit methods such as RK4. It is pertinent to mention that, for the purposes of this comparative analysis, we employ an IRK method of the same order as the RK4 algorithm.

\subsection{FitzHugh-Nagumo}
The subsequent benchmark problem pertains to a relatively simple (in its formulation) yet important model in the field of mathematical neuroscience that characterizes the oscillatory and nonlinear dynamics in the electrical activity of neurons, known as the FitzHugh-Nagumo model \cite{FitzHugh1961}, commonly abbreviated as FHN model. In spite of its simplicity, this mathematical model is utilized in various neuroscience applications, particularly in illustrating how neurons are capable of generating action potentials in response to stimuli. FHN can simulate the periodic oscillatory behavior observed in neuronal activity, such as brain rhythms, and the propagation of electrical waves in a network of neurons. FHN serves as a foundational framework, establishing a basis for the subsequent development of more sophisticated models that seek to capture the complexities associated with neuronal activity. The system of differential equations that encapsulates the underlying dynamics of this model is expressed in eq.\eqref{FHN}:.\eqref{FHN}:
\begin{subequations} \label{FHN}
\begin{equation}
\dot{v}(t) = v(t) - w(t) - {1 \over 3}v^{3}(t) + 0.5,
\end{equation}
\begin{equation}
\dot{w}(t) = 0.04 v(t) -0.028w(t) +0.032.
\end{equation}
\end{subequations}
By employing the initial condition $
\begin{bmatrix}
v(0) & w(0)
\end{bmatrix}^{T}
\begin{bmatrix}
0.0 & 0.0
\end{bmatrix}^{T}$, we generate time-series data on the interval $t \in [0, 200]$. We use deep IRK-SINDy with a fully connected neural network architecture characterized by a periodic activation function, SIREN \cite{Sitzmann2020}, consisting of $3$ hidden layers each containing $32$ neurons, thereby providing a robust framework for capturing the periodic dynamics of the system on the long time intervals. We establish the thresholding value to be $\lambda=0.01$ and proceed to learn the governing equations in the space of polynomials up to degree $3$, employing a total of $8$ sequential thresholding iterations with $20,000$ epochs allocated for the first iteration, followed by $5,000$ epochs for each subsequent iteration, with a learning rate of $10^{-4}$  and $10^{-3}$ (Except for the case $m = 101$ that we utilize $5\times 10^{-4}$) for the DNN and the coefficient matrix $\xi$ in the first iteration, respectively.
\begin{figure}[]
\flushleft
\begin{subfigure}{0.46\textwidth}
\caption{}
\begin{tikzpicture}[
node distance=0.0cm and 0.0cm,
block/.style={rectangle, draw, rounded corners, fill=teal!0, minimum width=2cm, minimum height=1cm},
>=Stealth,
every node/.style={align=center} 
]

\begin{scope}[local bounding box=NN, rounded corners=5pt, minimum width=2cm, minimum height=0cm]

\node[block, rounded corners=5pt, minimum width=0.99\textwidth, minimum height=7.5cm, dashed] (b1) {};

\node[below right=0.2cm of b1.north west] (img1) {\includegraphics[width=0.45\textwidth]{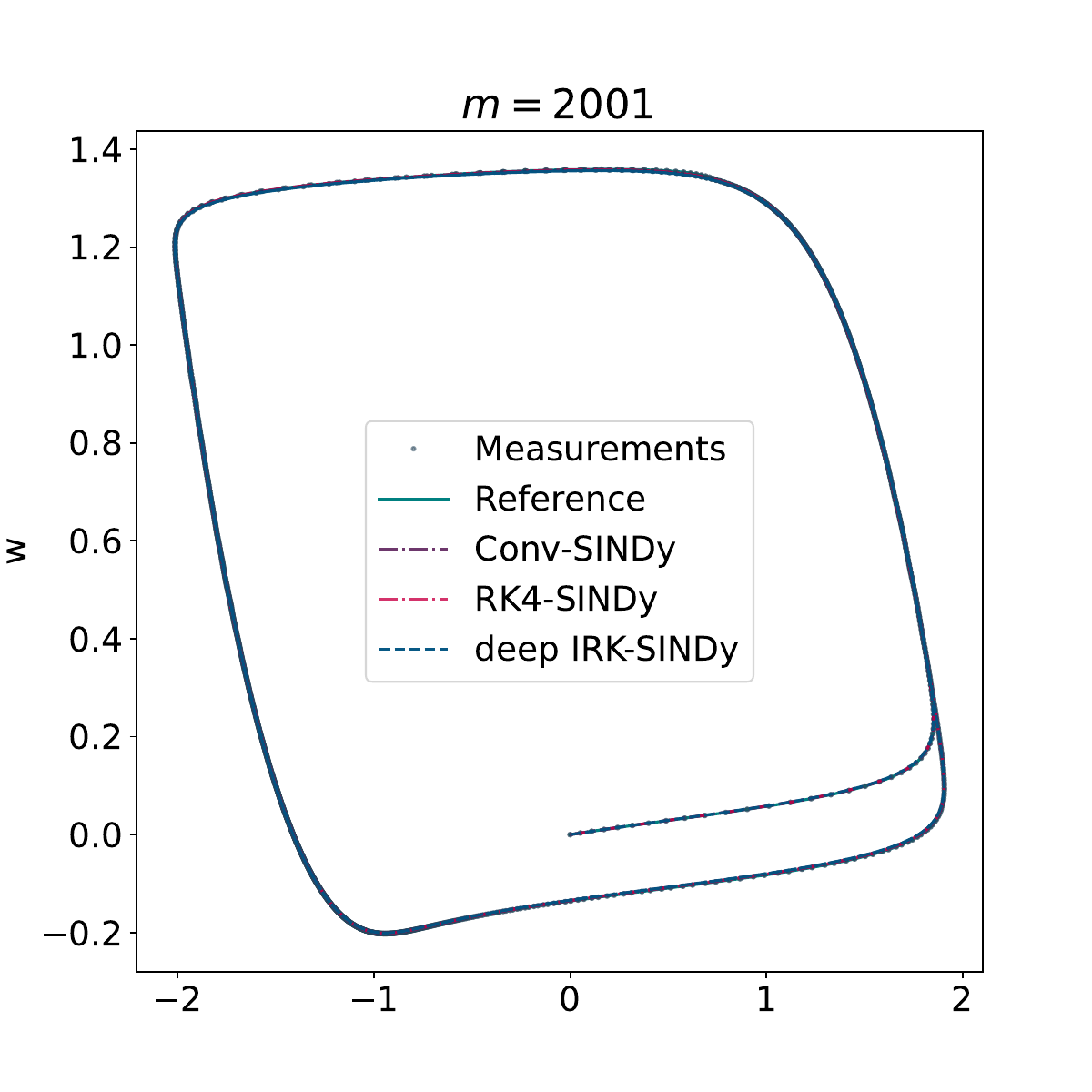}};
\node[below=of img1] (img2) {\includegraphics[width=0.45\textwidth]{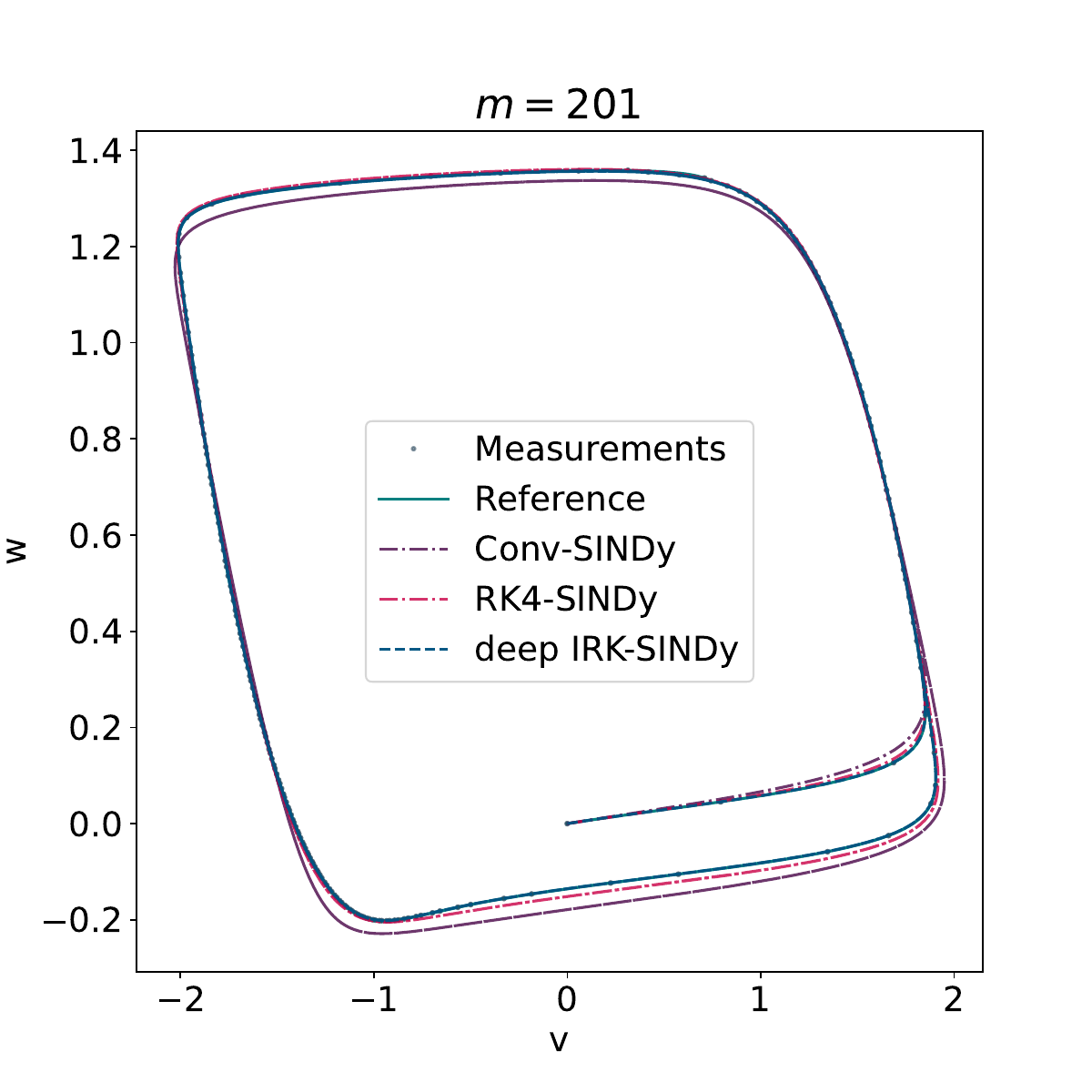}};
\node[right=of img1] (img3) {\includegraphics[width=0.45\textwidth]{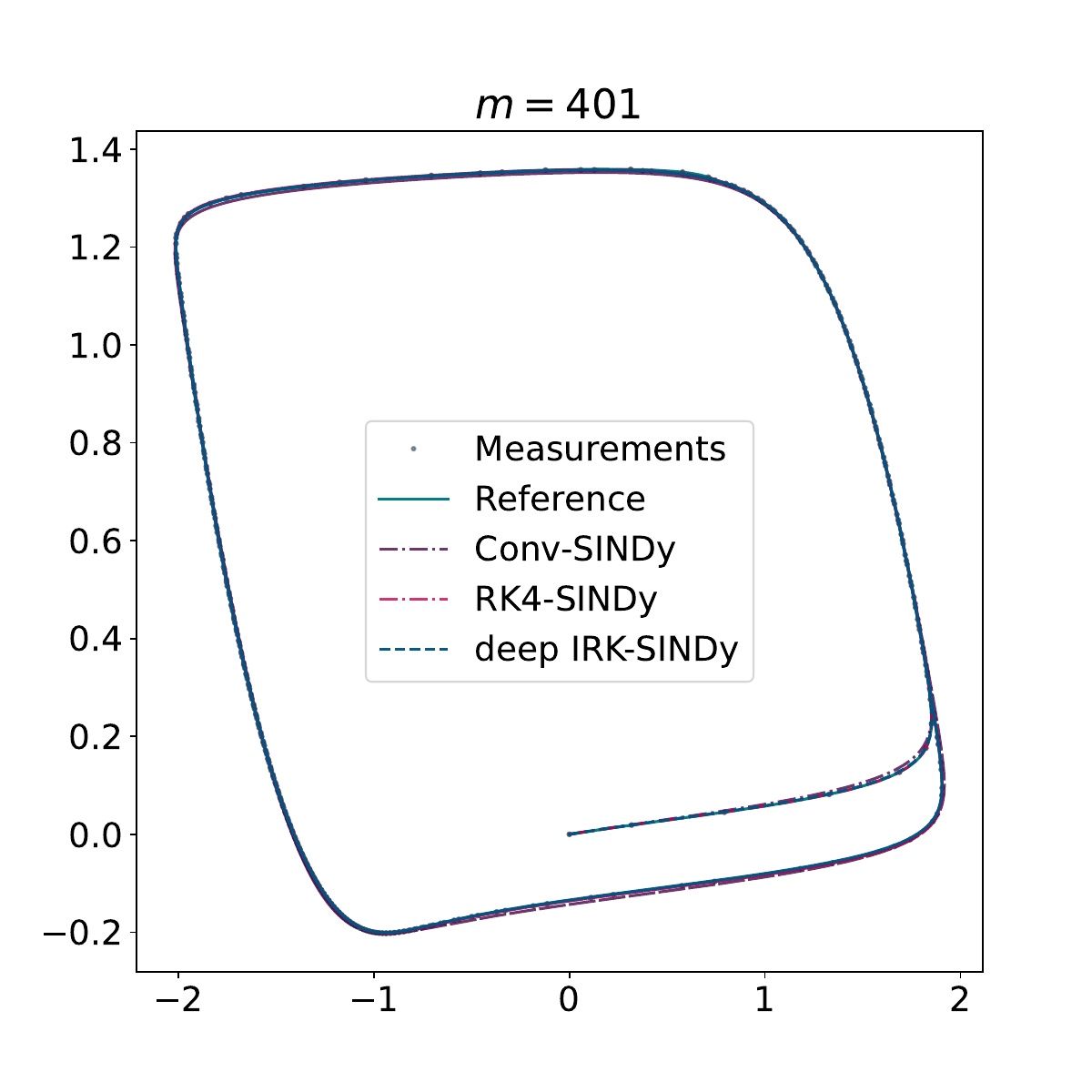}};
\node[below=of img3] (img4) {\includegraphics[width=0.45\textwidth]{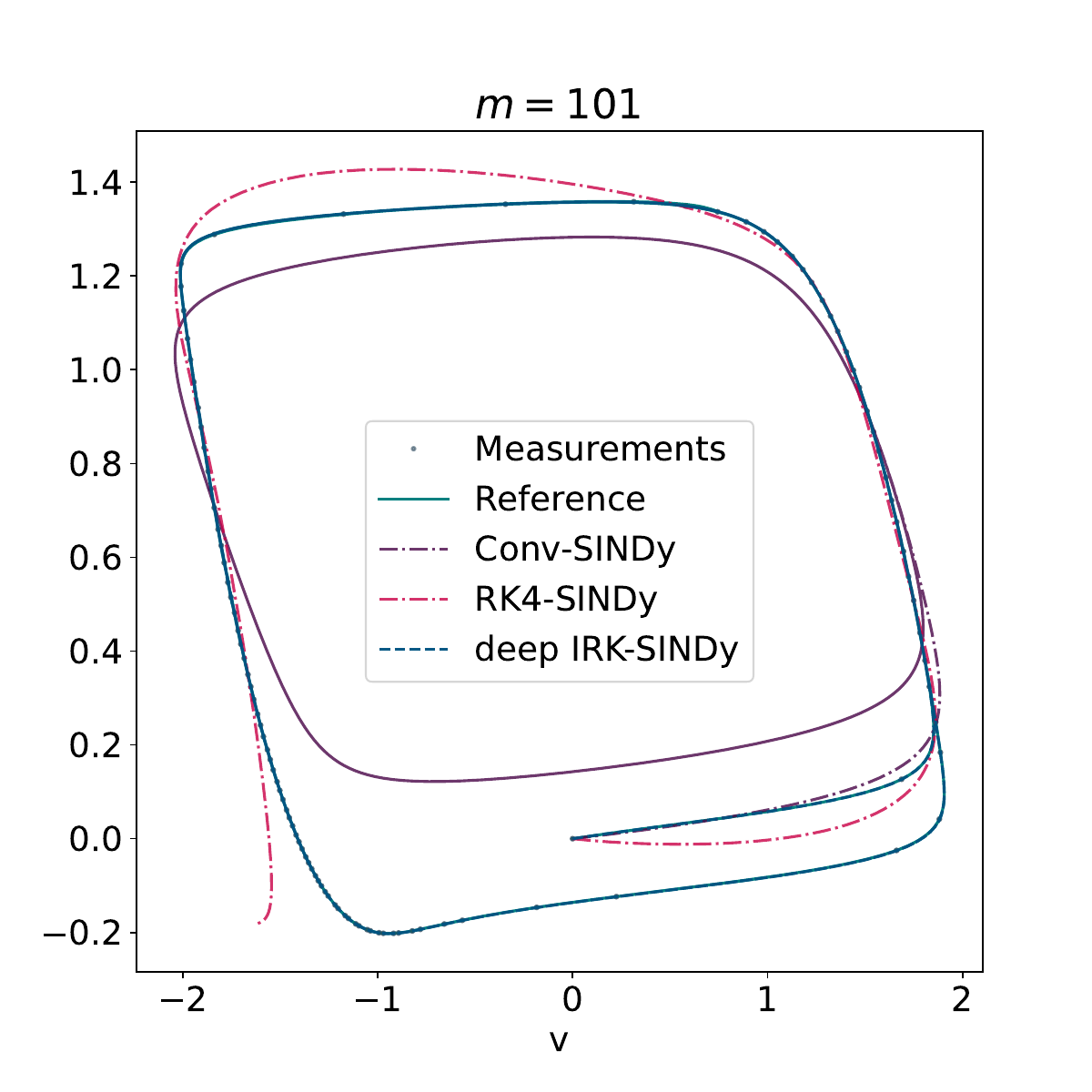}};
\end{scope}
\end{tikzpicture}
\end{subfigure}
\begin{subfigure}{0.47\textwidth}
\caption{}
\begin{tikzpicture}[
node distance=0.0cm and 0.0cm,
neuron/.style={circle, fill=teal!0, draw, minimum size=0.4cm},
block/.style={rectangle, draw, rounded corners, fill=teal!0, minimum width=2cm, minimum height=1cm},
>=Stealth,
every node/.style={align=center} 
]

\begin{scope}[local bounding box=NN, rounded corners=5pt, minimum width=2cm, minimum height=0cm]

\node[block, rounded corners=5pt, minimum width=0.99\textwidth, minimum height=7.5cm, dashed] (b1) {};

\node[below right=0.2cm of b1.north west] (img1) {\includegraphics[width=0.45\textwidth]{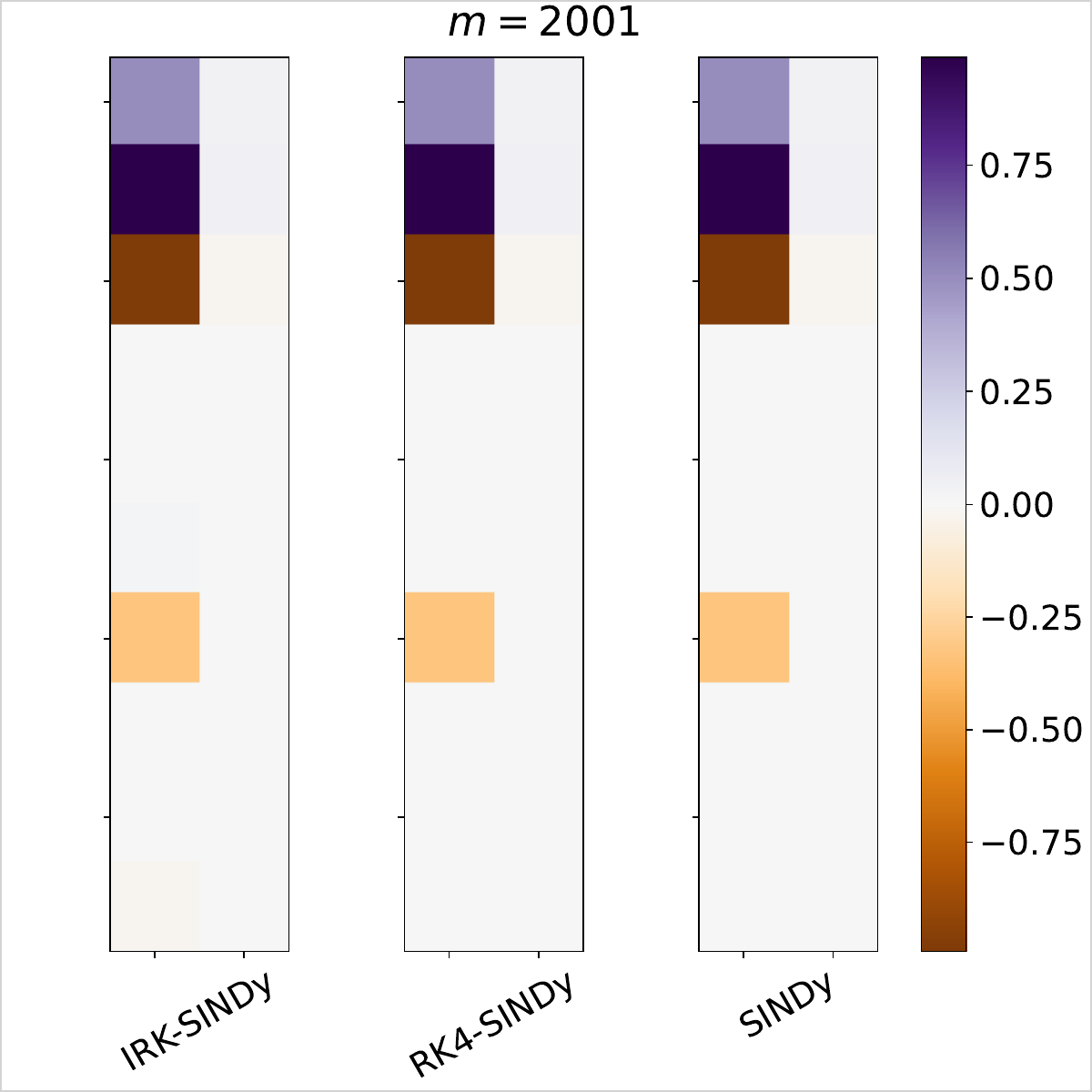}};
\node[below=of img1] (img2) {\includegraphics[width=0.45\textwidth]{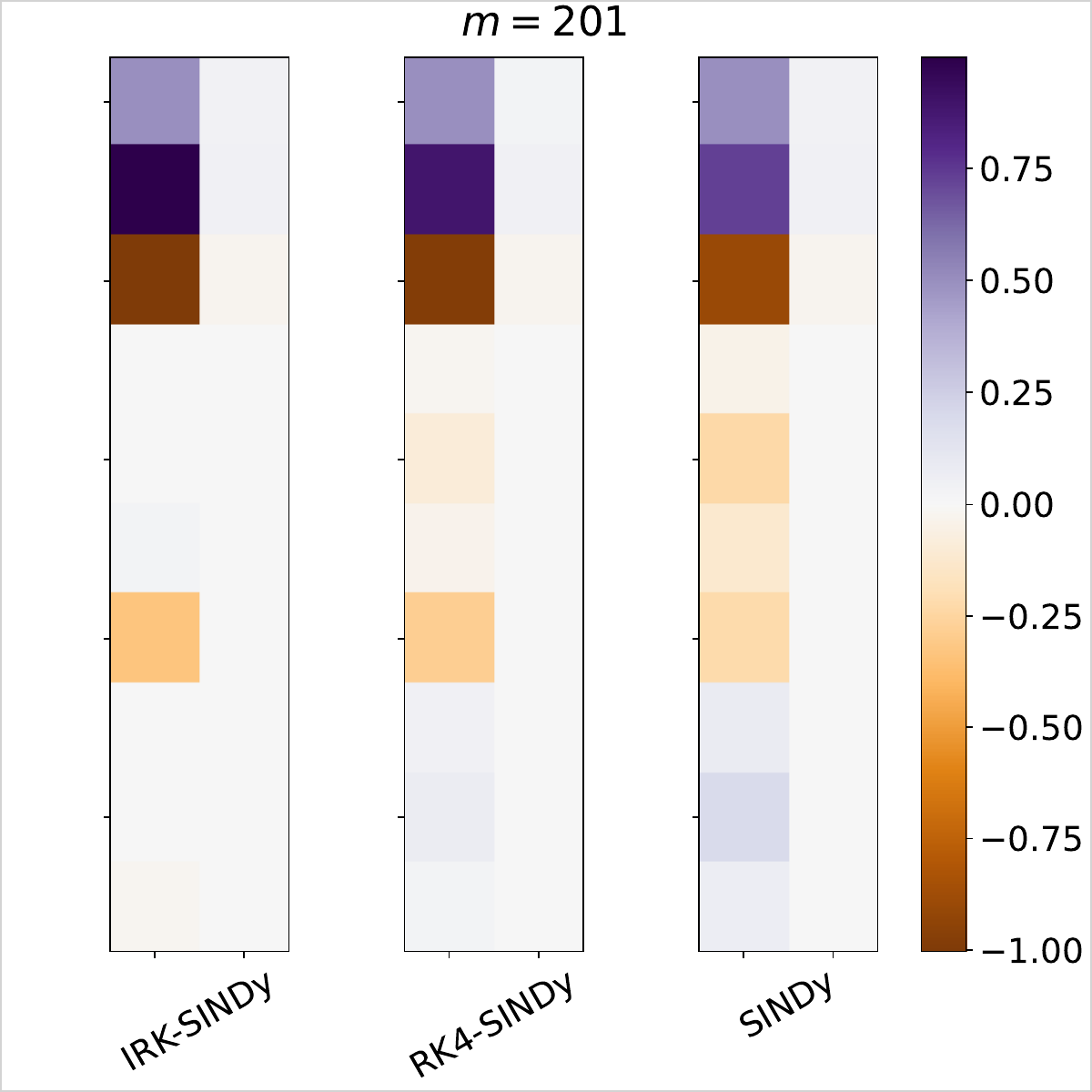}};
\node[right=of img1] (img3) {\includegraphics[width=0.45\textwidth]{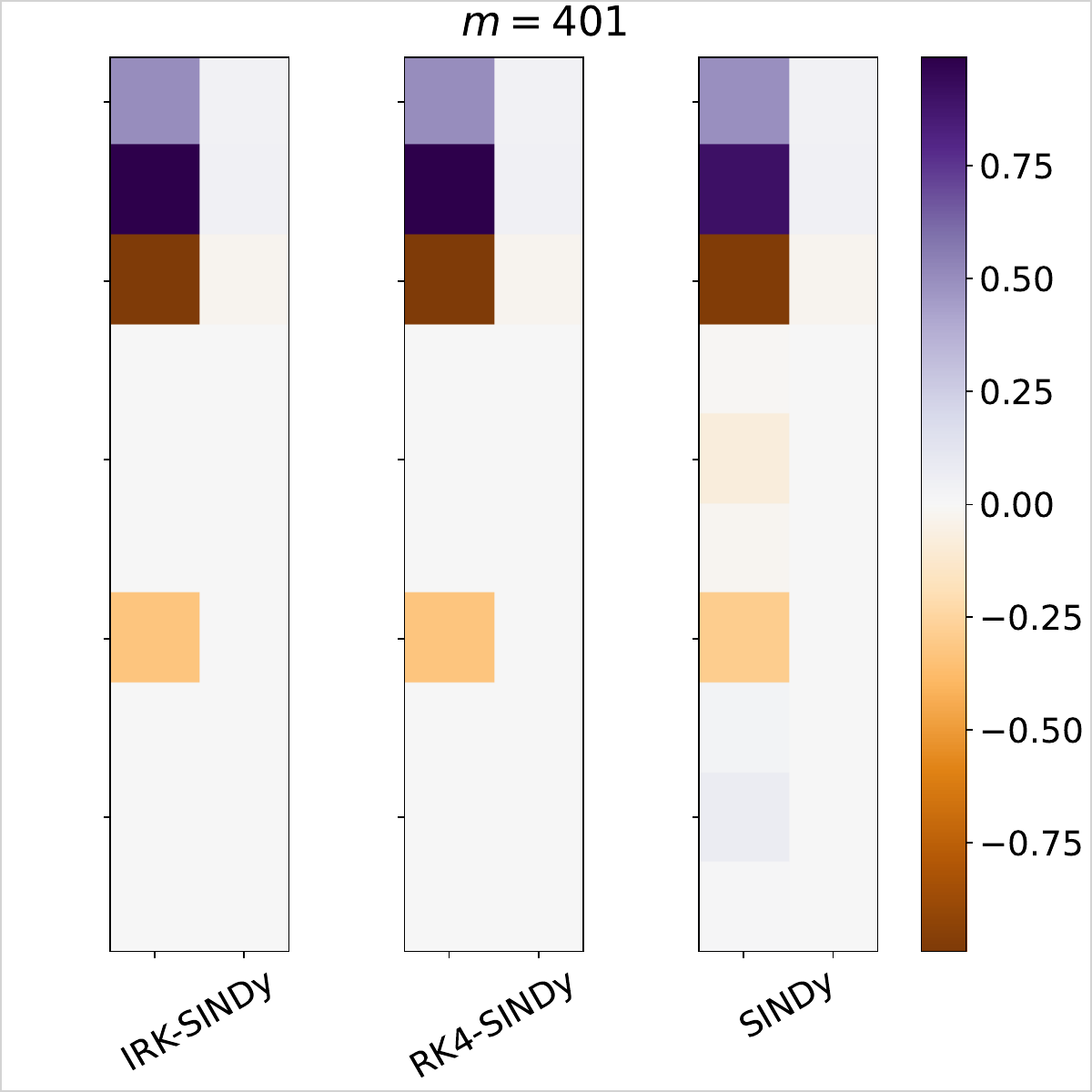}};
\node[below=of img3] (img4) {\includegraphics[width=0.45\textwidth]{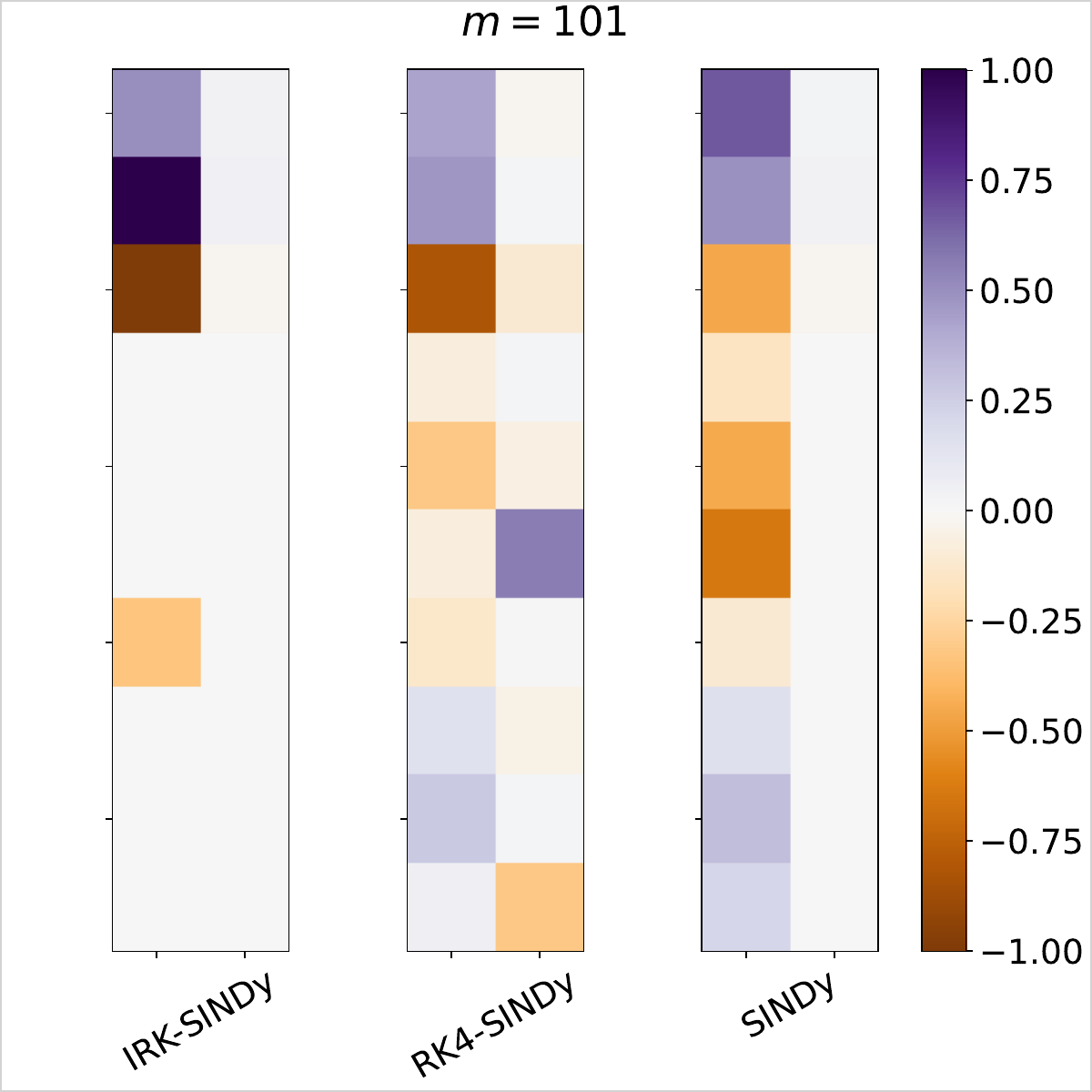}};
\end{scope}
\end{tikzpicture}
\end{subfigure}
\caption{Fitz–Hugh Nagumo model: a comparison of the reference model and recovered models using data collected at constant time stepsize. a. Phase portraits, b. coefficient matrices for $m\in\{ 2001, 401, 201, 101 \}$. IRK-SINDy results a parsimonious and generalizable model for biologically motivated models such as FHN even in data scarcity.}\label{fig6}
\end{figure}
In Figure\ref{fig6}, we present a comparative analysis of our proposed methodology against the performance of Conv-SINDy and RK4-SINDy, for various values of m. The incorporation of a periodic activation function within the DNN architecture significantly enhances the model's capability to effectively learn from periodic data\cite{Essakine2024}, which is of paramount importance for accurately predicting the stage values associated with the IRKs. Consequently, as illustrated in Figure\ref{fig6}, it becomes evident that deep IRK-SINDy demonstrates a markedly reduced dependence on the quantity of data points compared to the two alternative methodologies, thus revealing its superior efficacy in reconstructing the dynamics of the FHN model.

\subsection{Lorenz attractor}
Here, to explore the efficacy of deep IRK-SINDy in identifying chaotic dynamics, we examine the nonlinear $3D$ Lorenz system \cite{Lorenz1963} as the next illustrative example. A distinctive characteristic of this system is its sensitivity to initial conditions, which makes it a prominent candidate within the field of data-driven discovery of dynamical systems. The governing equations of this system are given as eq.\eqref{Lorenz}:
\begin{subequations} \label{Lorenz}
\begin{equation}
\dot{x}(t) = 10 (y(t) - x(t)),
\end{equation}
\begin{equation}
\dot{y}(t) = x(t) (28 - z(t)) - y(t),
\end{equation}
\begin{equation}
\dot{z}(t) = x(t) y(t) - {8 \over 3} z(t).
\end{equation}
\end{subequations}
Utilizing the initial condition$\begin{bmatrix}
x(0) & y(0) & z(0)
\end{bmatrix}^{T} =
\begin{bmatrix}
-8 & 7 & 27
\end{bmatrix}^{T}$, we conduct measurements on the time interval $t \in [0, 10]$ with different sample sizes $m$. We employ a DNN comprising one hidden layer with $256$ neurons to represent the nonlinear dynamics at the stage values of IRKs with the periodic activation function SIREN. For the sequential thresholding procedure, we use $10$ thresholding iterations, commencing with $20,000$ epochs in the first iteration and $2,000$ epochs in the subsequent iterations, with a thresholding parameter set at $\lambda=0.5$. During the training phase, we adopt a learning rate of $10^{-3}$ for the DNN parameters and a learning rate of $10^{-2}$ for the coefficient matrix while searching for active terms in the space of polynomials up to degree $2$. Moreover, at the end of each iteration, both learning rates are diminished by specific scales. It is important to note that due to the very large standard deviation of the state variables (significantly exceeding $1$), the library of polynomials may become ill-conditioned, thereby disrupting the optimization process. Therefore, prior to initiating the learning process, it is necessary to preprocess the data through scaling or normalization. This procedure is conducted such that the transformed data exhibits a mean of $0$ and a variance of $1$ \cite{Goyal2022}. It is crucial that the scaling of the data does not influence the interaction among the state variables, thus, the sparsity of the identified dynamics remains consistent with that of the reference dynamics. With a similar configuration, we repeat this numerical experiment for the scenario where $1\%$ noise is added to the data.

\begin{figure}[]
\centering
\begin{subfigure}{\textwidth}
\caption{Noise free}
\begin{tikzpicture}[
node distance=0.0cm and 0.0cm,
block/.style={rectangle, draw, rounded corners, fill=teal!0, minimum width=2cm, minimum height=1cm},
>=Stealth,
every node/.style={align=center} 
]

\begin{scope}[local bounding box=NN, rounded corners=5pt, minimum width=2cm, minimum height=0cm]

\node[block, rounded corners=5pt, minimum width=\textwidth, minimum height=4.5cm, dashed] (b1) {};

\node[below right=0.2cm of b1.north west] (img1) {\includegraphics[width=0.24\textwidth]{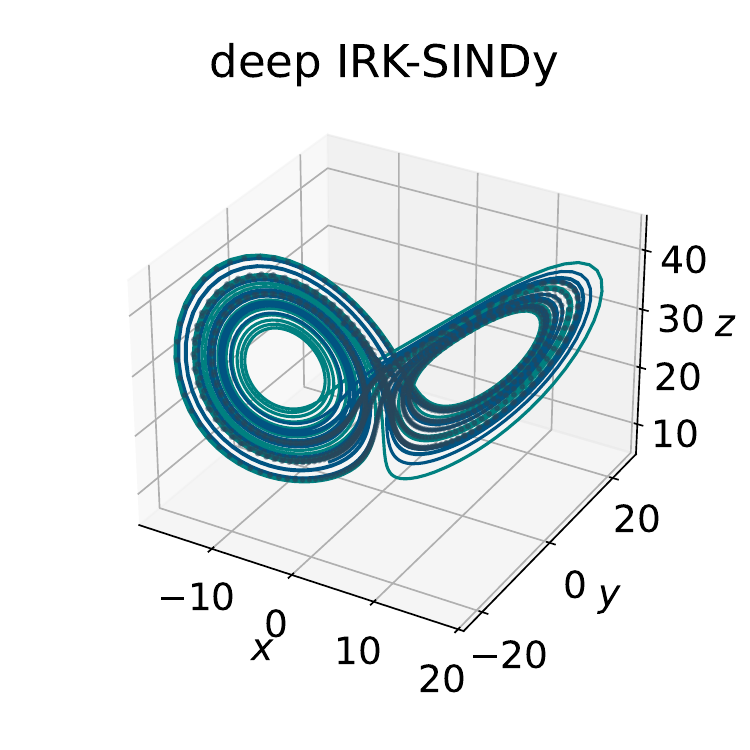}};
\node[right=of img1] (img2) {\includegraphics[width=0.24\textwidth]{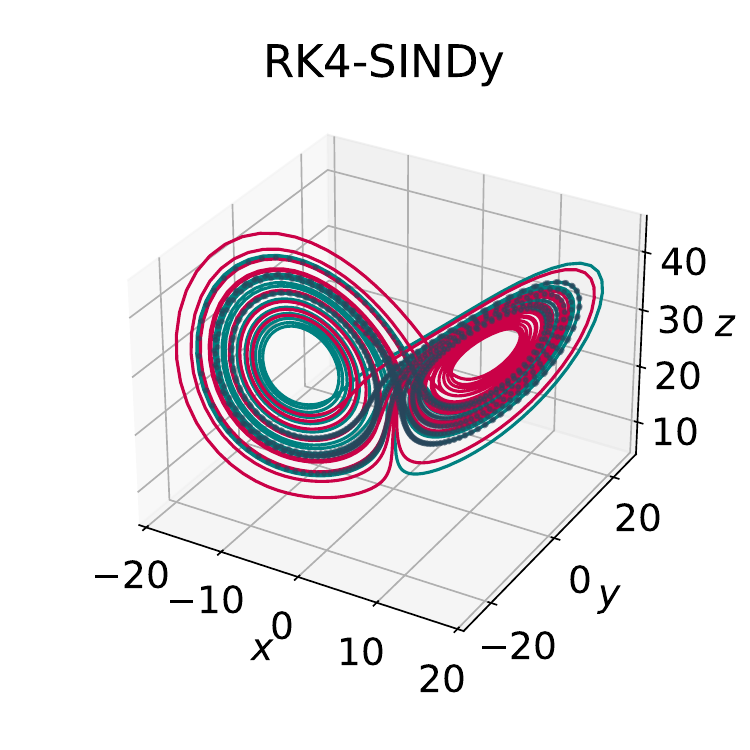}};
\node[right=of img2] (img3) {\includegraphics[width=0.24\textwidth]{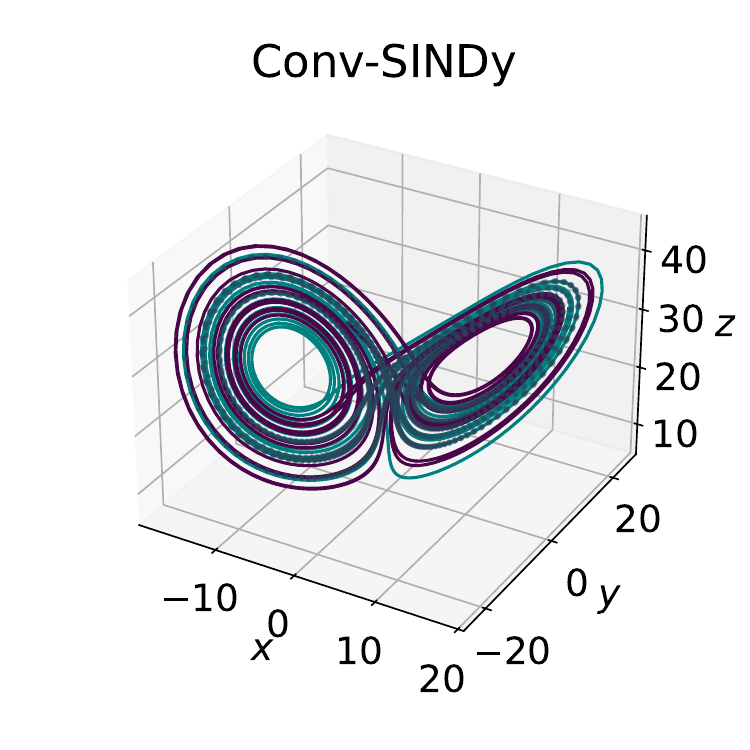}};
\node[right=of img3] (img4) {\includegraphics[width=0.21\textwidth]{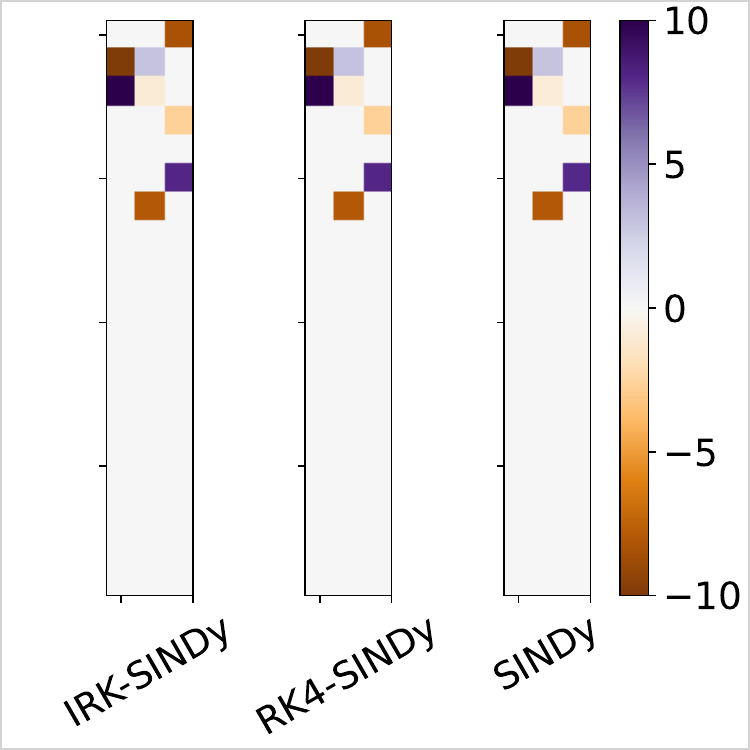}};
\end{scope}
\end{tikzpicture}
\end{subfigure}
\begin{subfigure}{\textwidth}
\caption{$\sigma=0.01$}
\begin{tikzpicture}[
node distance=0.0cm and 0.0cm,
neuron/.style={circle, fill=teal!0, draw, minimum size=0.4cm},
block/.style={rectangle, draw, rounded corners, fill=teal!0, minimum width=2cm, minimum height=1cm},
>=Stealth,
every node/.style={align=center} 
]

\begin{scope}[local bounding box=NN, rounded corners=5pt, minimum width=2cm, minimum height=0cm]

\node[block, rounded corners=5pt, minimum width=\textwidth, minimum height=4.5cm, dashed] (b1) {};

\node[below right=0.2cm of b1.north west] (img1) {\includegraphics[width=0.24\textwidth]{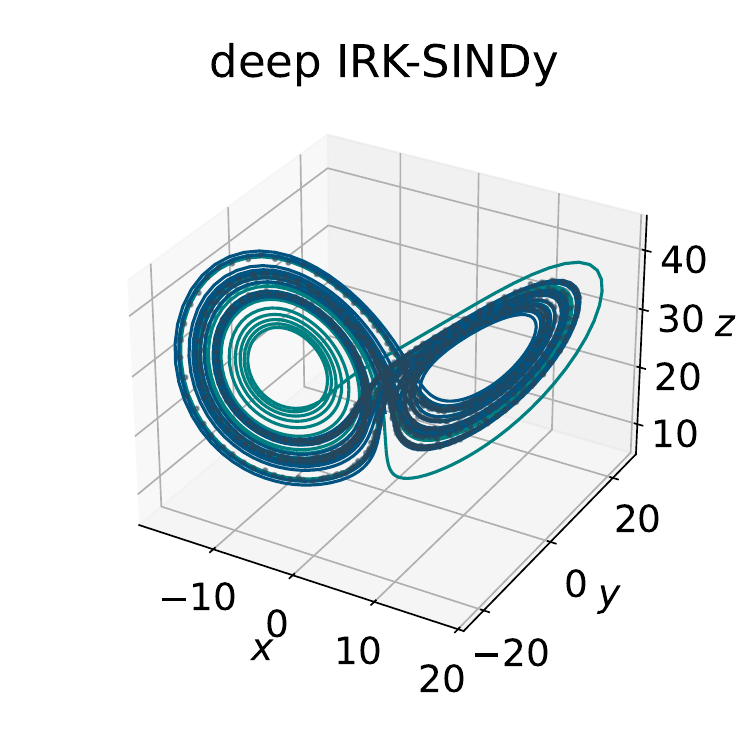}};
\node[right=of img1] (img2) {\includegraphics[width=0.24\textwidth]{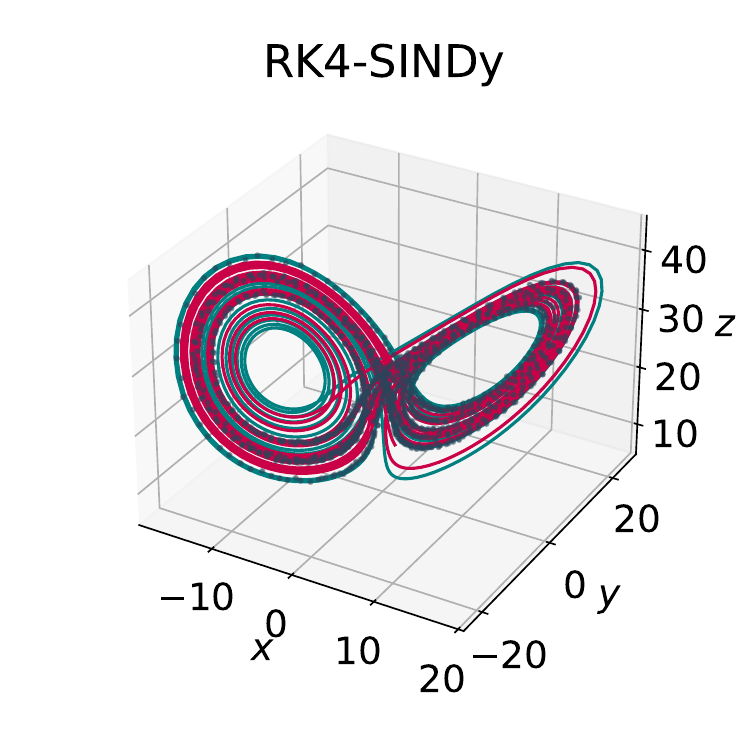}};
\node[right=of img2] (img3) {\includegraphics[width=0.24\textwidth]{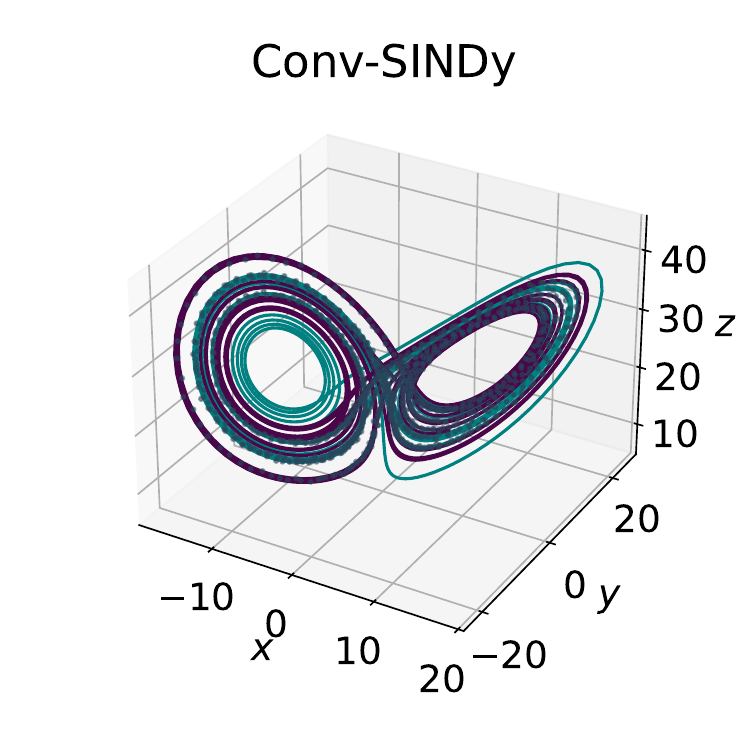}};
\node[right=of img3] (img4) {\includegraphics[width=0.21\textwidth]{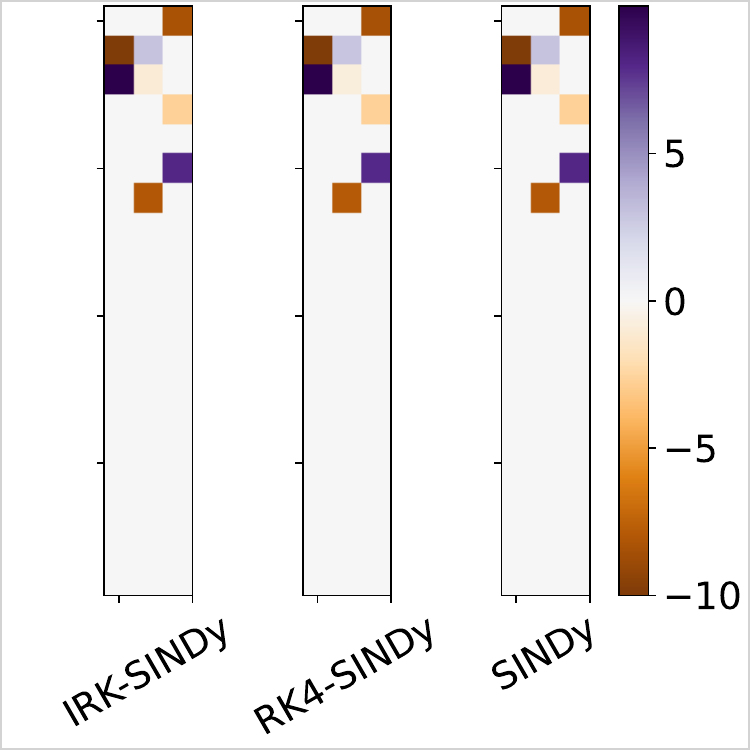}};
\end{scope}
\end{tikzpicture}
\end{subfigure}
\caption{Lorenz attractor model: a comparison of the reference model and recovered models using data collected at constant time stepsize in the cases a. noise free, b. $1$ percent noise $\sigma=0.01$.}\label{figlorenz}
\end{figure}

All three approaches IRK-SINDy, RK4-SINDy and Conv-SINDy are able to correctly identify the active nonlinear features, while the coefficient matrix obtained from IRK-SINDy is closer to the coefficients of the reference model. However, the Lorenz system has a positive Lyapunov exponent and small differences between the reference and discovered models cause exponential growth in the forecasted differences. As evidenced in Figure\ref{figlorenz}, although small deviations in the dynamic coefficients significantly affect the dynamics due to the highly chaotic behavior of the system, the bi-stable structure of the attractor is well captured even in presence of noise. While Figure \ref{figlorenz} shows that IRK-SINDy discovers more robust and parsimonious models against noise.

\subsection{Lotka–Volterra predator–prey model}
\begin{figure}[]
\flushleft
\begin{subfigure}{0.46\textwidth}
\caption{}
\begin{tikzpicture}[
node distance=0.0cm and 0.0cm,
block/.style={rectangle, draw, rounded corners, fill=teal!0, minimum width=2cm, minimum height=1cm},
>=Stealth,
every node/.style={align=center} 
]

\begin{scope}[local bounding box=NN, rounded corners=5pt, minimum width=2cm, minimum height=0cm]

\node[block, rounded corners=5pt, minimum width=0.97\textwidth, minimum height=7.5cm, dashed] (b1) {};

\node[below right=0.2cm of b1.north west] (img1) {\includegraphics[width=0.45\textwidth]{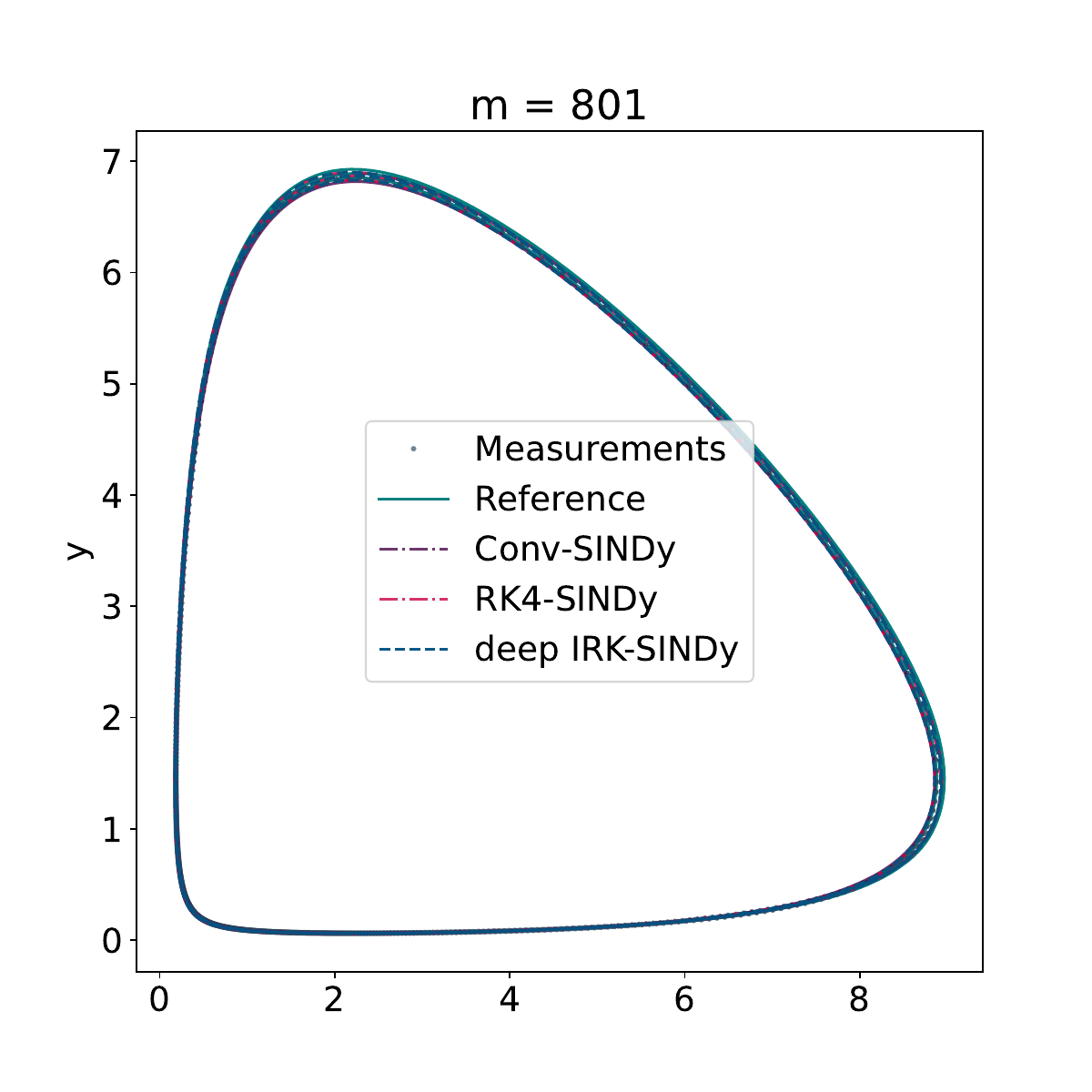}};
\node[below=of img1] (img2) {\includegraphics[width=0.45\textwidth]{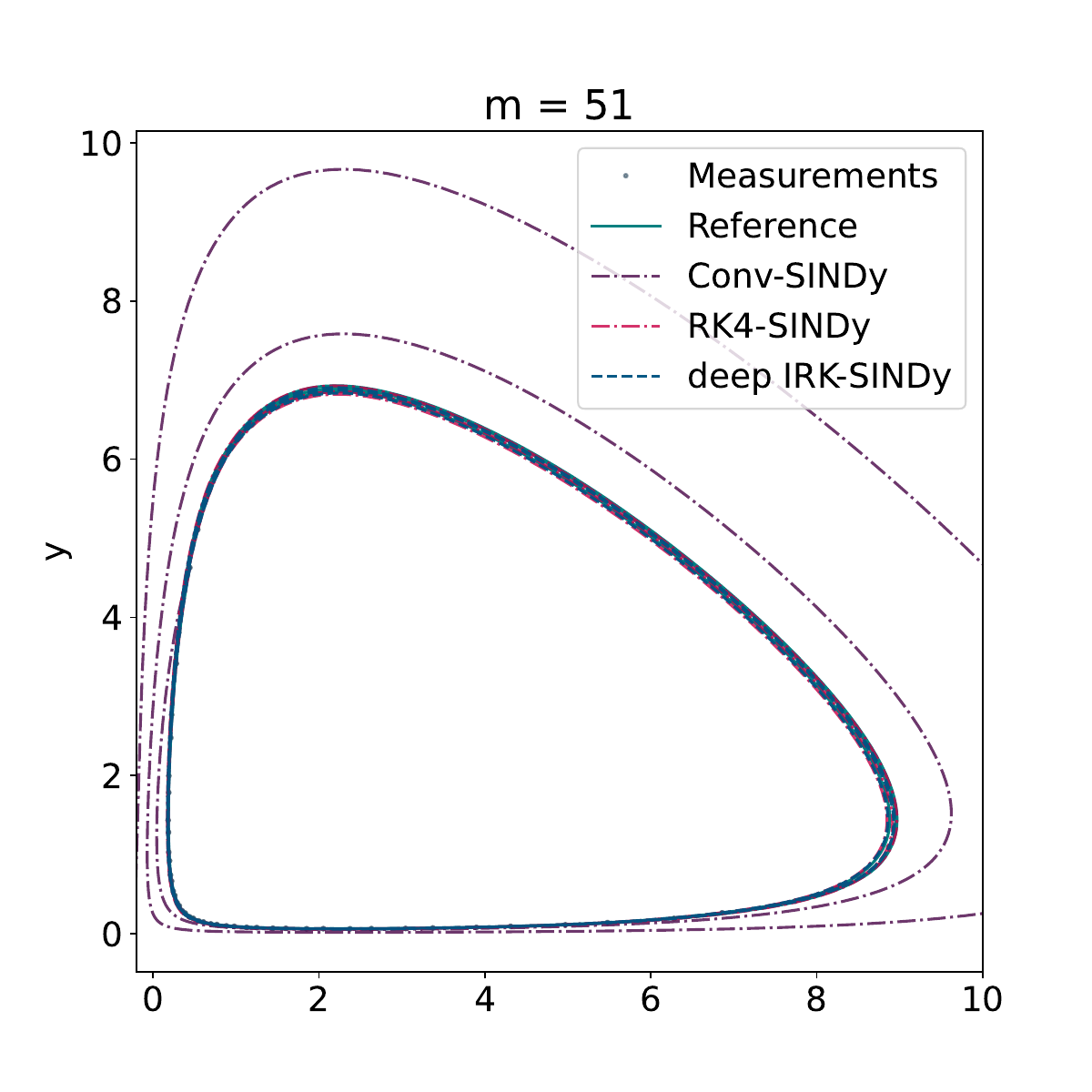}};
\node[right=of img1] (img3) {\includegraphics[width=0.45\textwidth]{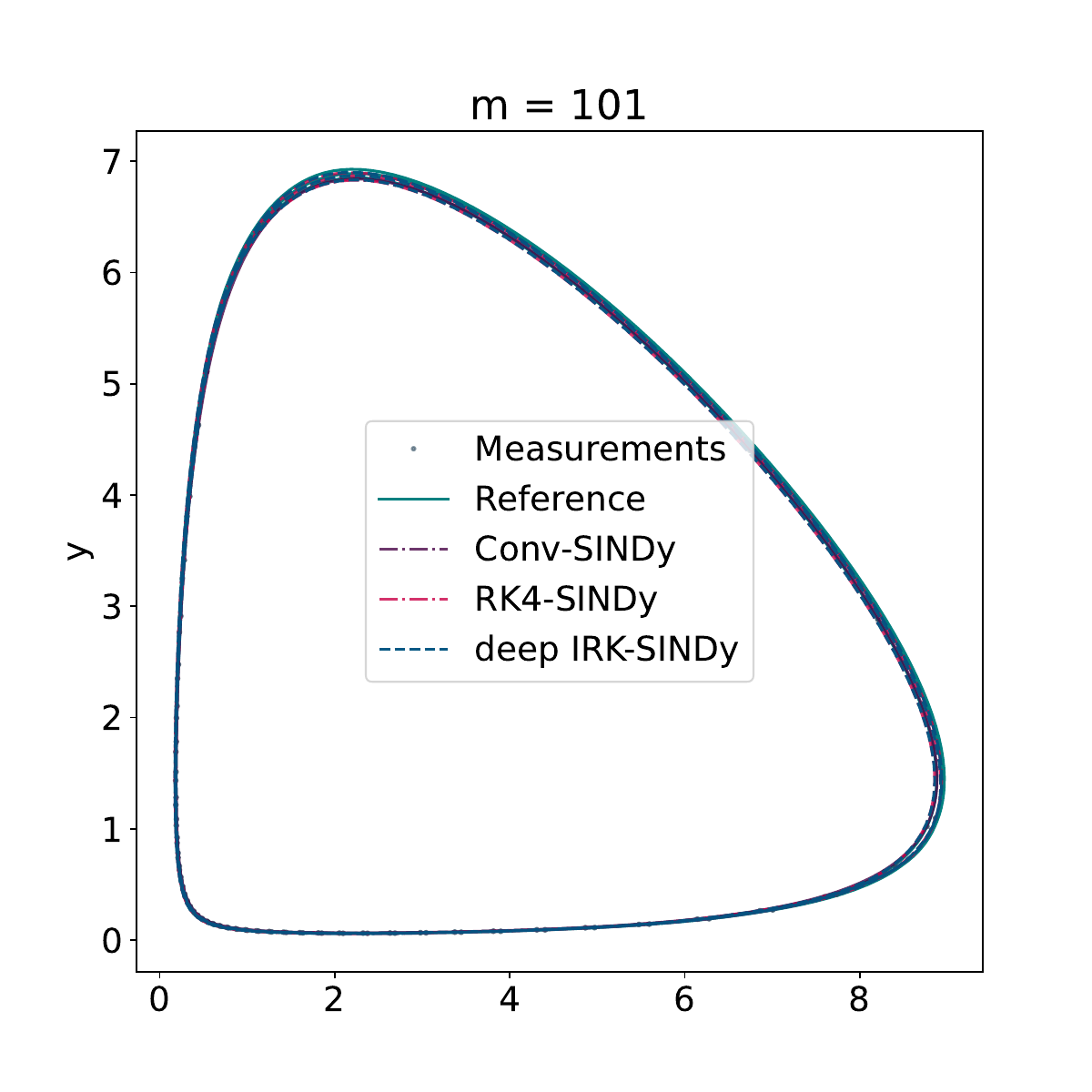}};
\node[below=of img3] (img4) {\includegraphics[width=0.45\textwidth]{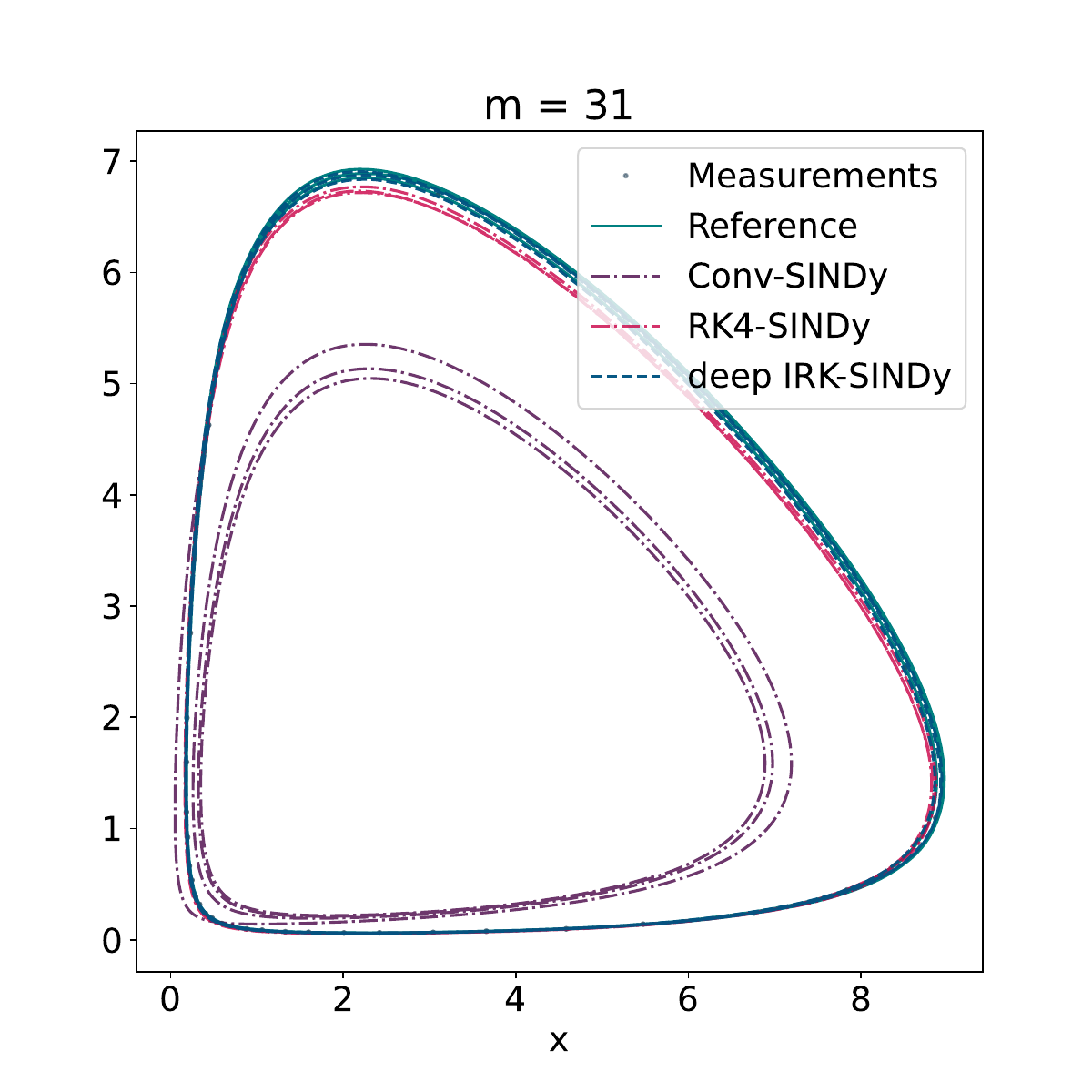}};
\end{scope}
\end{tikzpicture}
\end{subfigure}
\begin{subfigure}{0.47\textwidth}
\caption{}
\begin{tikzpicture}[
node distance=0.0cm and 0.0cm,
neuron/.style={circle, fill=teal!0, draw, minimum size=0.4cm},
block/.style={rectangle, draw, rounded corners, fill=teal!0, minimum width=2cm, minimum height=1cm},
>=Stealth,
every node/.style={align=center} 
]

\begin{scope}[local bounding box=NN, rounded corners=5pt, minimum width=2cm, minimum height=0cm]

\node[block, rounded corners=5pt, minimum width=0.99\textwidth, minimum height=7.5cm, dashed] (b1) {};

\node[below right=0.2cm of b1.north west] (img1) {\includegraphics[width=0.45\textwidth]{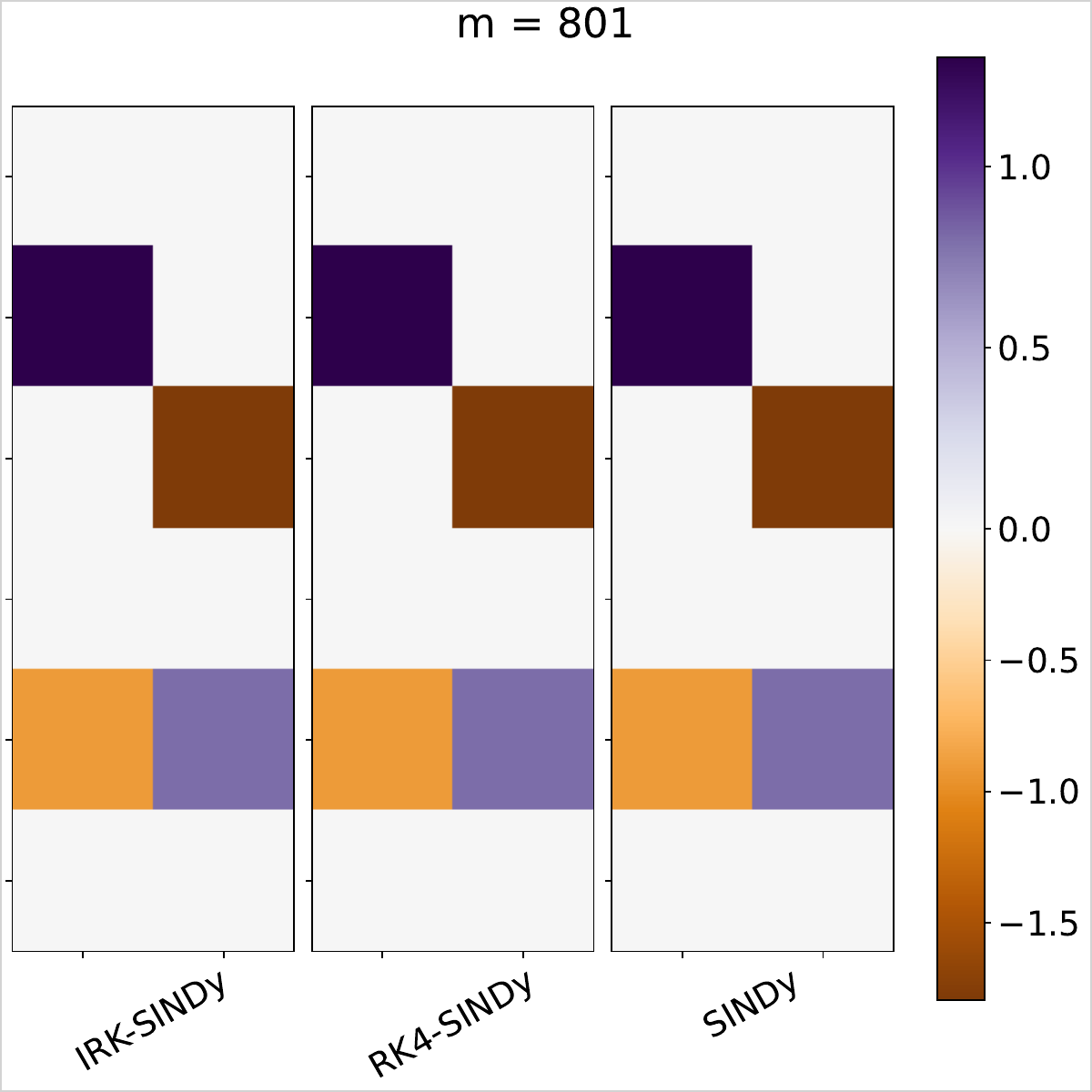}};
\node[below=of img1] (img2) {\includegraphics[width=0.45\textwidth]{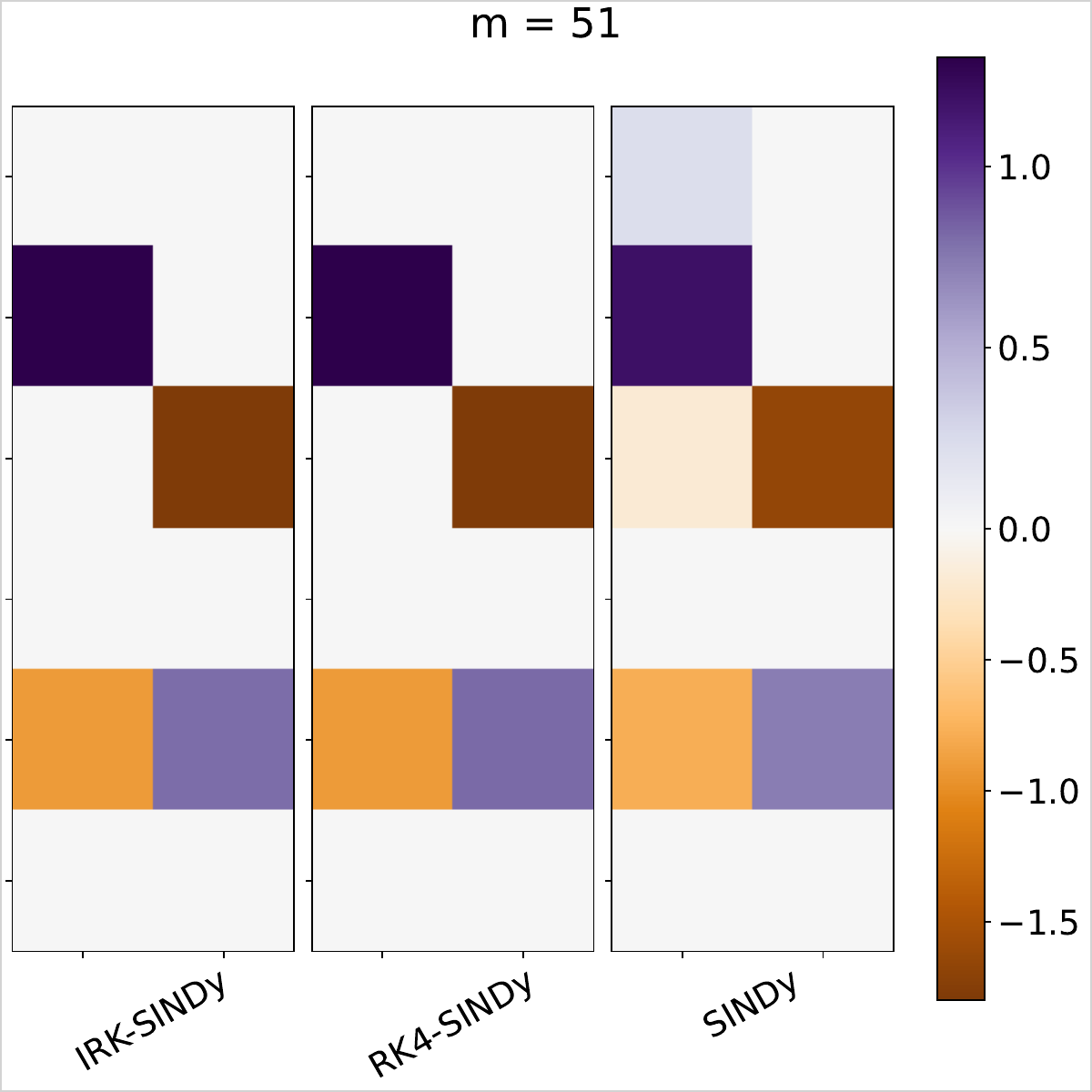}};
\node[right=of img1] (img3) {\includegraphics[width=0.45\textwidth]{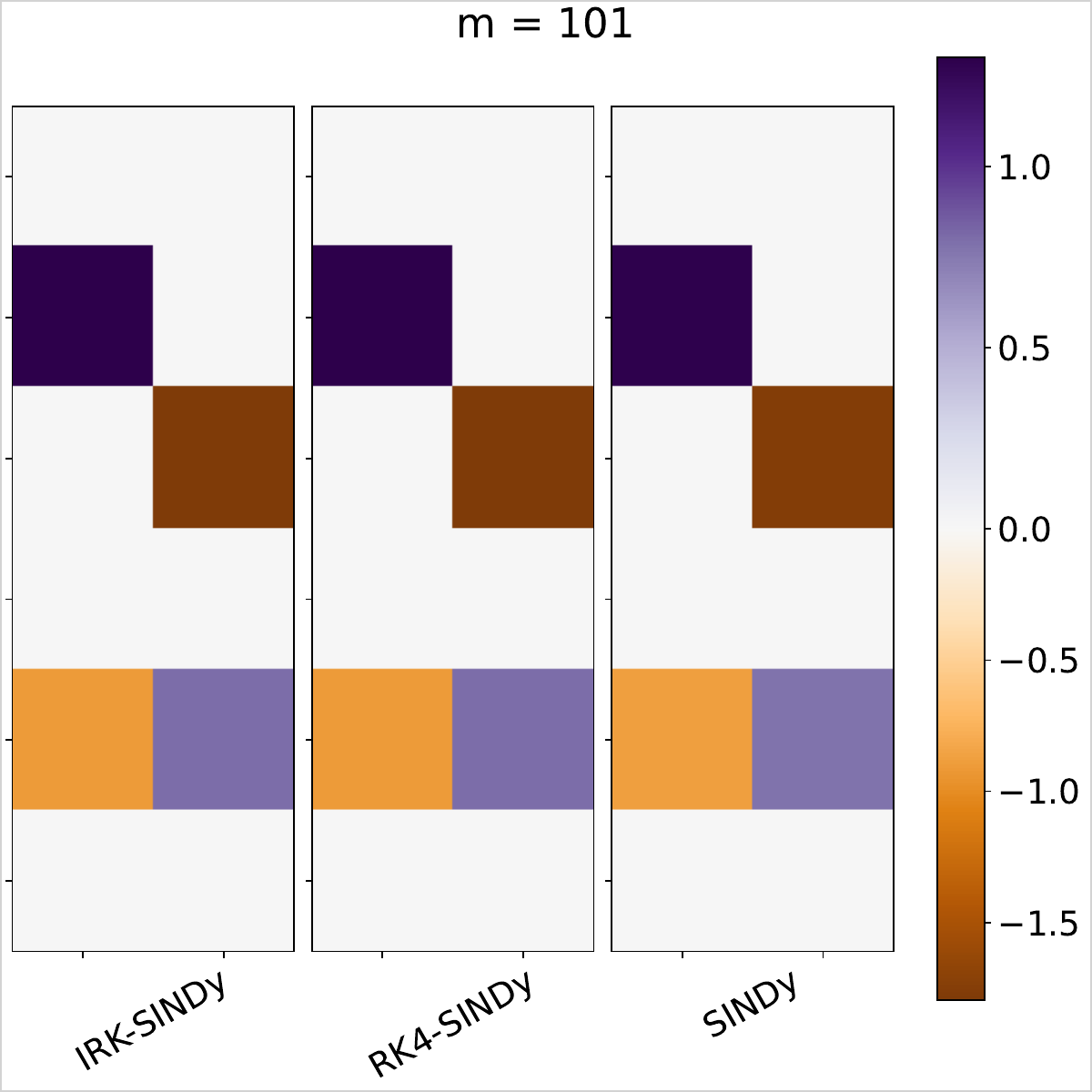}};
\node[below=of img3] (img4) {\includegraphics[width=0.45\textwidth]{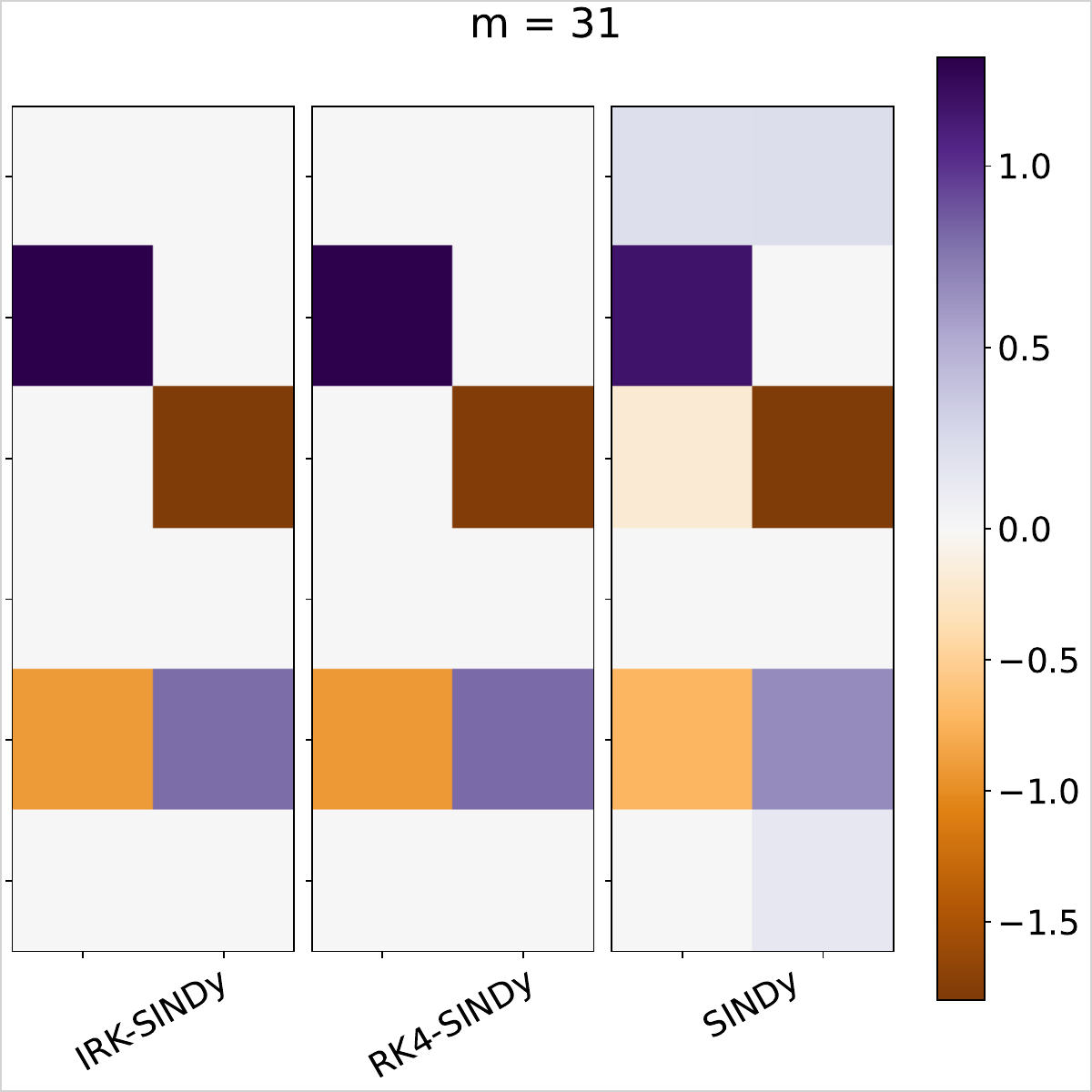}};
\end{scope}
\end{tikzpicture}
\end{subfigure}
\caption{Lotka–Volterra model: a comparison of the reference model and recovered models using data collected at constant time stepsize. a. Phase portraits, b. coefficient matrices for $m\in\{ 801, 101, 51, 31 \}$. Compared to Conv-SINDy and RK4-SINDy, IRK-SINDy produces a sufficiently sparse, interpretable and generalizable model.}\label{fig7}
\end{figure}

Next, we examine the Lotka-Volterra equations, which describe the predator-prey dynamics, and have recently been employed extensively as a significant biologically motivated benchmark problem in the data-driven discovery of the governing equations in biological systems\cite{Lejarza2022, Gutierrez2023, Wei2022}. This is crucial due to the fact that the predator-prey dynamics serve as the cornerstone for numerous mathematical models within the invwstigation of biological systems, particularly in the field of systems biology of cancer (where cancer cells are conceptualized as prey and the immune system as the predator)\cite{Hamilton2022, Kareva2015}. In this context, we demonstrate the usefulness of proposed approach in identifying the Lotka-Volterra equations given by eqns.\eqref{LV} that describe the interaction between two species, denoted as $u$ for prey and $v$ for predator:
\begin{subequations} \label{LV}
\begin{equation}
\dot{u}(t) = \alpha u(t) - \beta u(t) v(t),
\end{equation}
\begin{equation}
\dot{v}(t) = - \gamma v(t) + \delta u(t) v(t),
\end{equation}
\end{subequations}
The time-series data are generated through solving the reference differential equation by parameters $\{\alpha = {2 \over 3}, \beta = {4 \over 3}, \gamma = \delta = 1\}$ on the time interval $t \in [0, 10]$ with the initial condition $[u(0), v(0)]^{T} = [1.8, 1.8]^{T}$. By setting the thresholding value $\lambda=0.1$, we employ deep IRK-SINDy using the SIREN periodic activation function alongside an architecture comprising $2$ hidden layers, each containing $64$ neurons, to discover the governing differential equations across different data quantities m. This approach, incorporating learning rates of $10^{-3}$ and $10^{-4}$ during the first iteration for $\xi$ and $\theta$, respectively, through sequential thresholding with $3$ iterations that is employed $6,000$ epochs in each iteration (except $25,000$ epochs for the first iteration), successfully captures the correct active terms within the nonlinear feature library, which in this illustrative example is selected as a polynomial space of up to degree $2$. In Figure \ref{fig7}, the obtained models are tested on the time interval $t \in [0, 20]$ and the efficacy of deep IRK-SINDy under data scarcity is depicted.

\subsection{Logistic growth model}
\begin{figure}[]
\flushleft
\begin{subfigure}{0.46\textwidth}
\caption{}
\begin{tikzpicture}[
node distance=0.0cm and 0.0cm,
block/.style={rectangle, draw, rounded corners, fill=teal!0, minimum width=2cm, minimum height=1cm},
>=Stealth,
every node/.style={align=center} 
]

\begin{scope}[local bounding box=NN, rounded corners=5pt, minimum width=2cm, minimum height=0cm]

\node[block, rounded corners=5pt, minimum width=0.97\textwidth, minimum height=7.5cm, dashed] (b1) {};

\node[below right=0.2cm of b1.north west] (img1) {\includegraphics[width=0.45\textwidth]{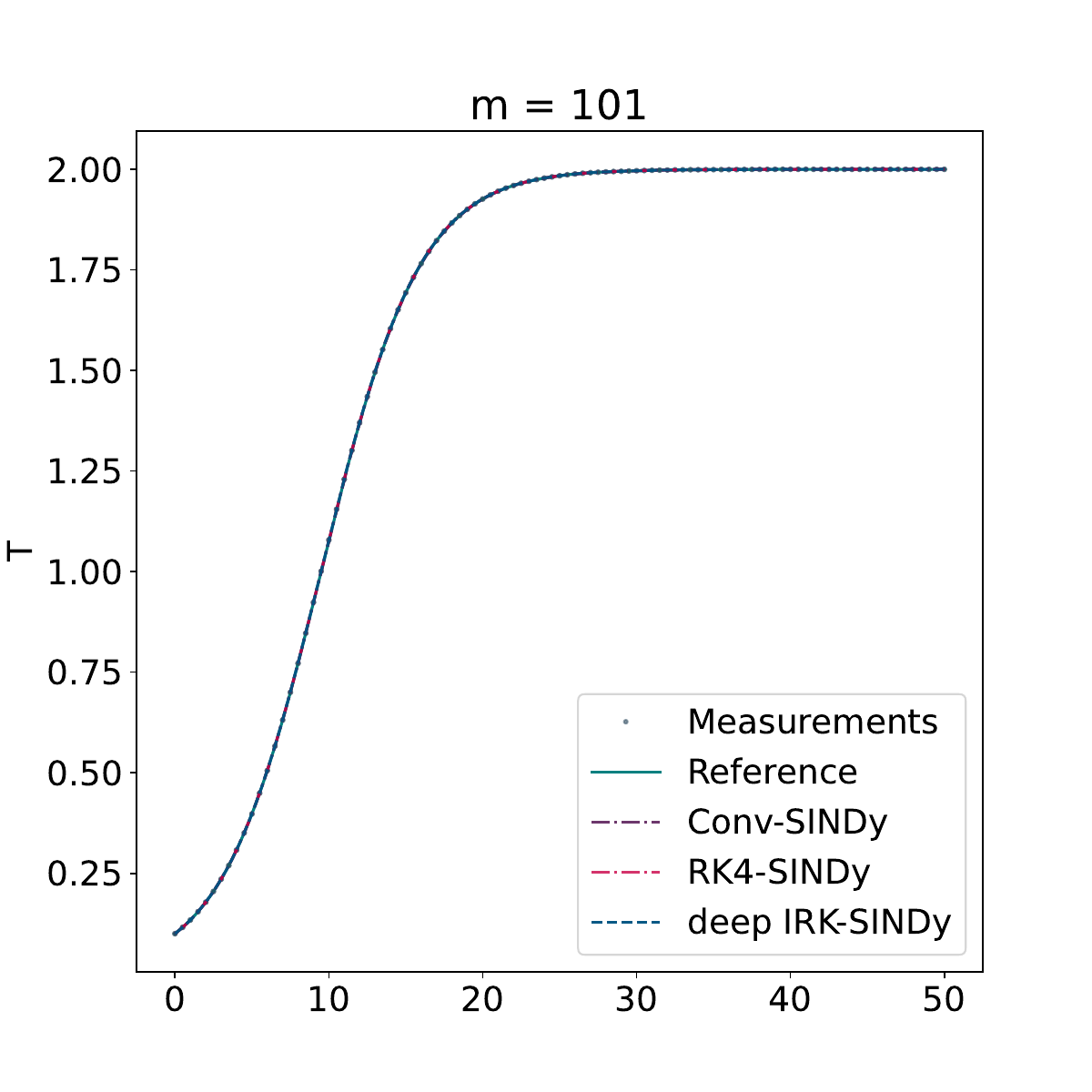}};
\node[below=of img1] (img2) {\includegraphics[width=0.45\textwidth]{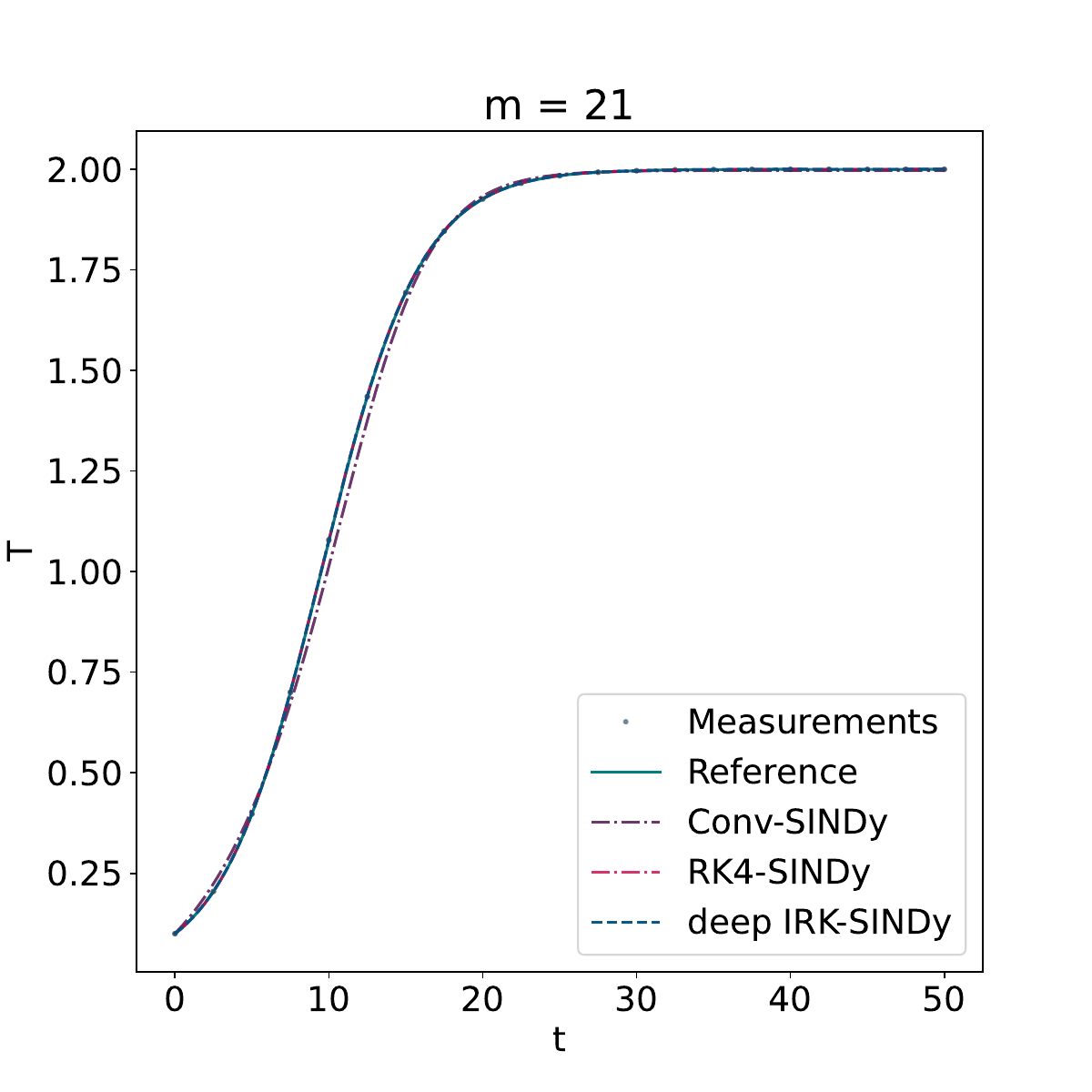}};
\node[right=of img1] (img3) {\includegraphics[width=0.45\textwidth]{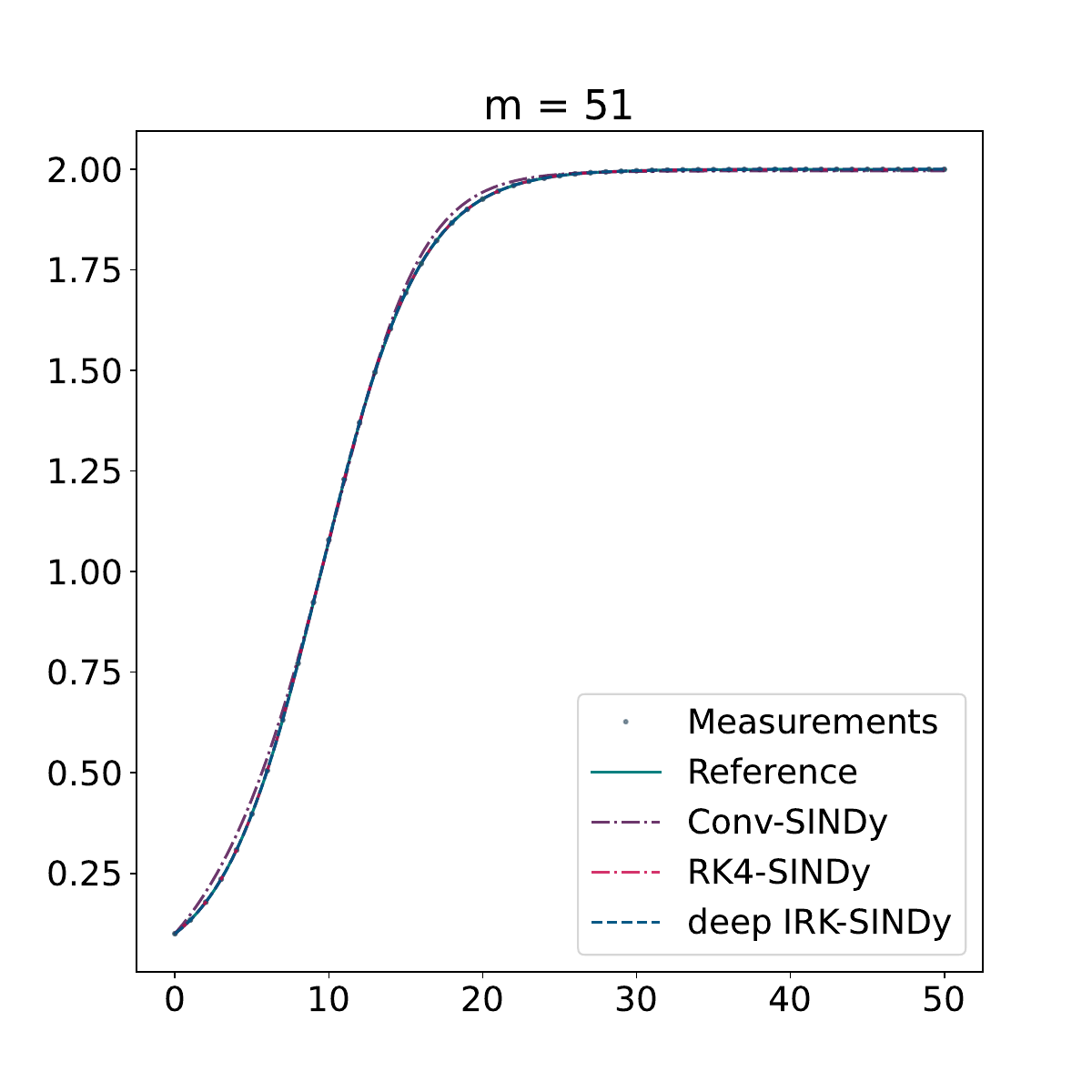}};
\node[below=of img3] (img4) {\includegraphics[width=0.45\textwidth]{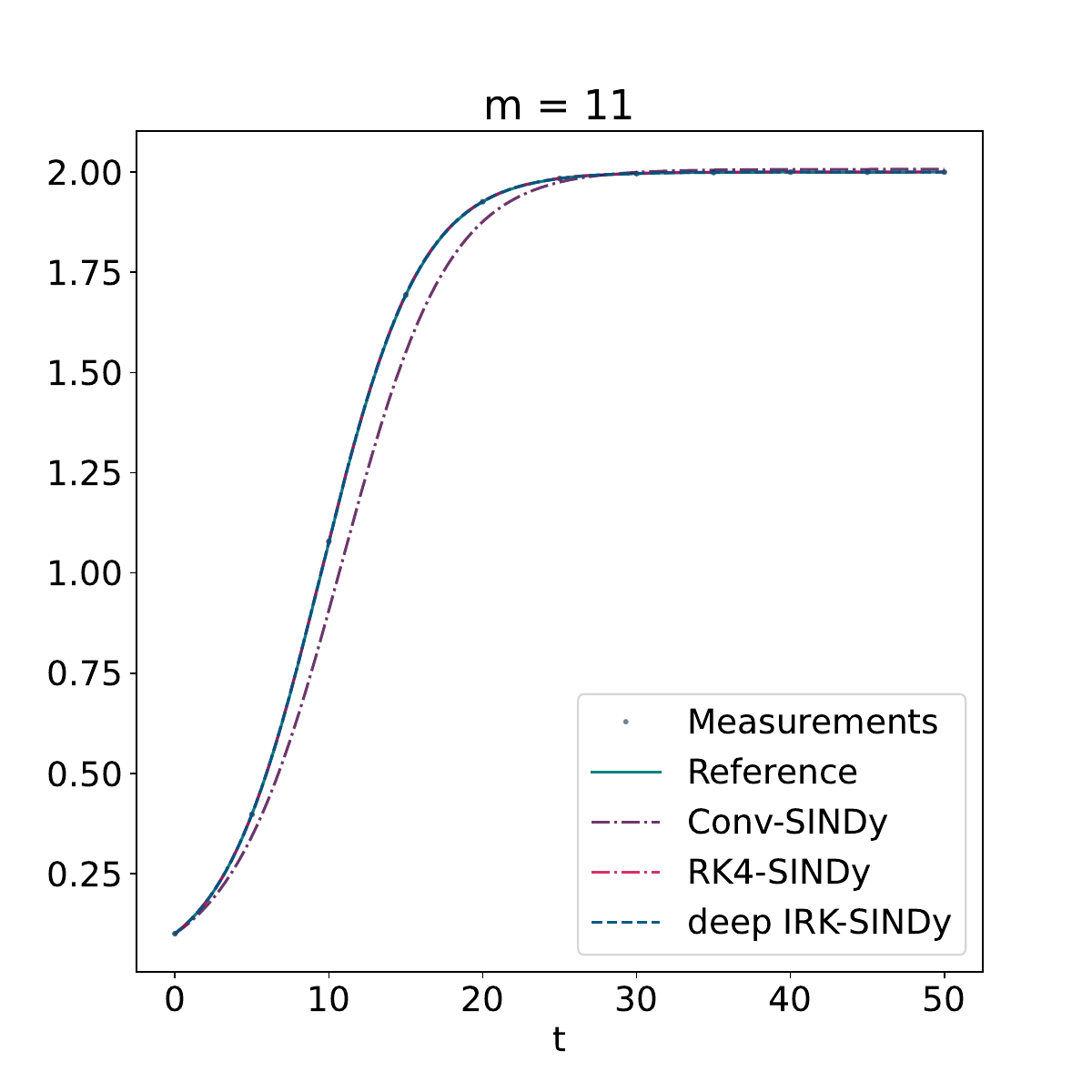}};
\end{scope}
\end{tikzpicture}
\end{subfigure}
\begin{subfigure}{0.47\textwidth}
\caption{}
\begin{tikzpicture}[
node distance=0.0cm and 0.0cm,
neuron/.style={circle, fill=teal!0, draw, minimum size=0.4cm},
block/.style={rectangle, draw, rounded corners, fill=teal!0, minimum width=2cm, minimum height=1cm},
>=Stealth,
every node/.style={align=center} 
]

\begin{scope}[local bounding box=NN, rounded corners=5pt, minimum width=2cm, minimum height=0cm]

\node[block, rounded corners=5pt, minimum width=0.99\textwidth, minimum height=7.5cm, dashed] (b1) {};

\node[below right=0.2cm of b1.north west] (img1) {\includegraphics[width=0.45\textwidth]{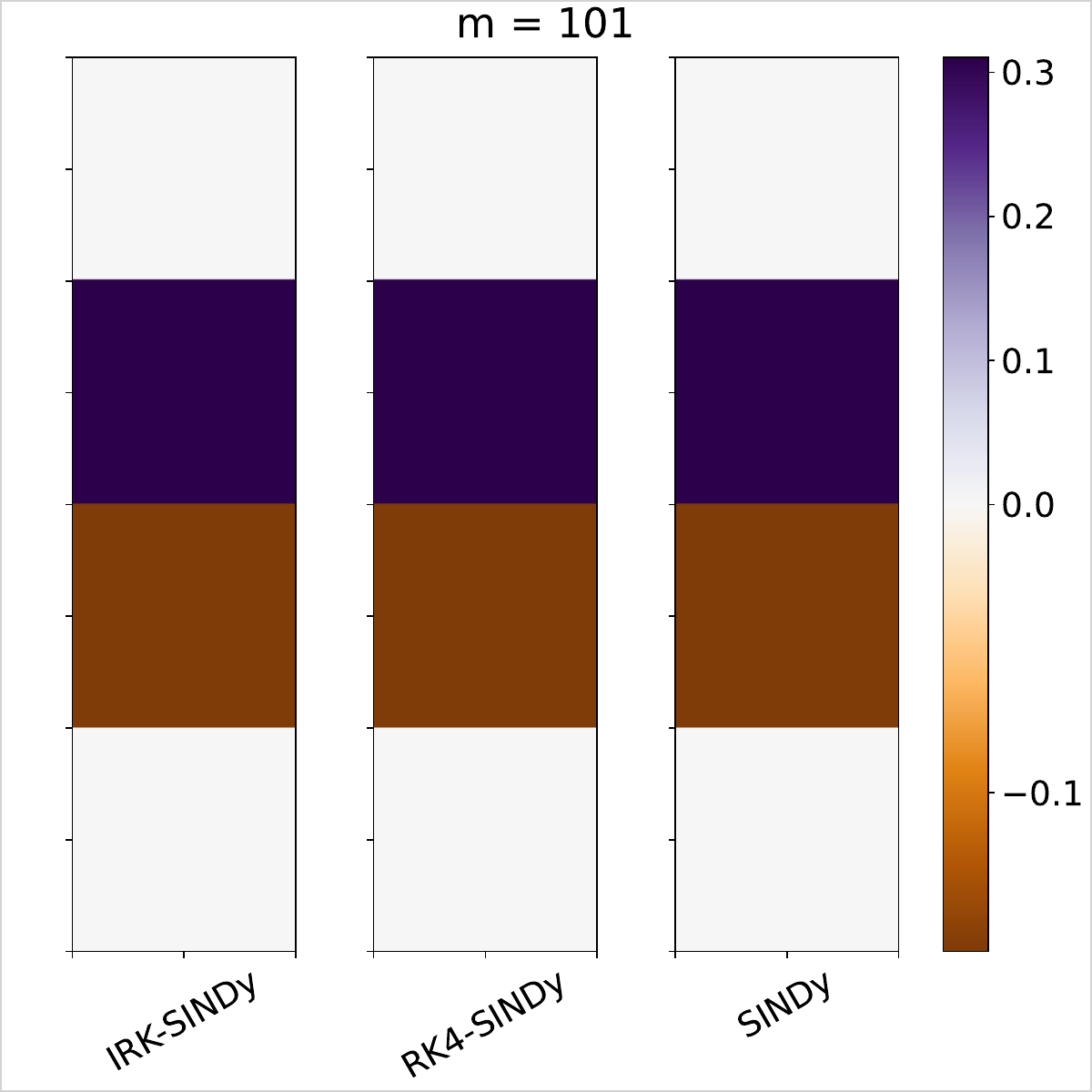}};
\node[below=of img1] (img2) {\includegraphics[width=0.45\textwidth]{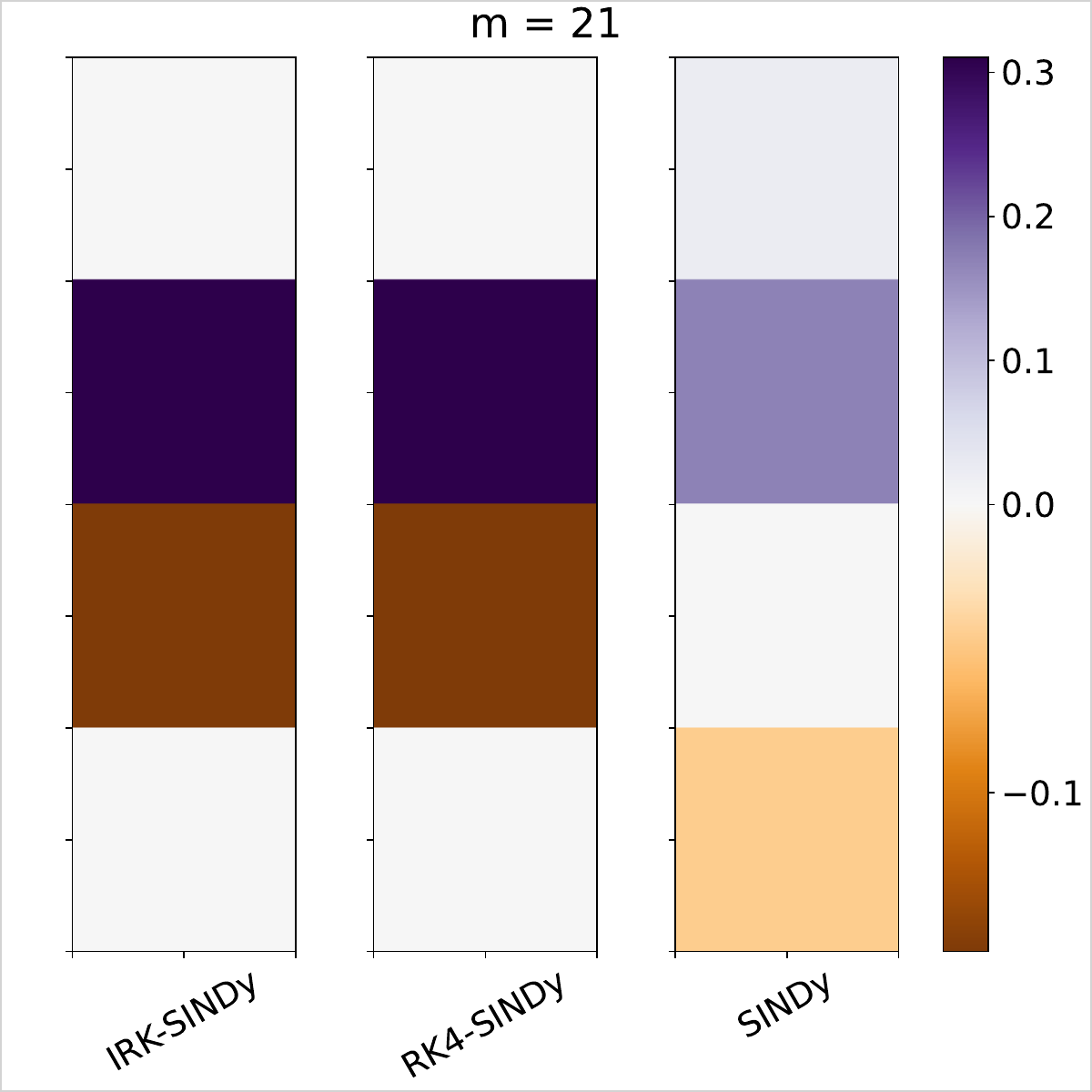}};
\node[right=of img1] (img3) {\includegraphics[width=0.45\textwidth]{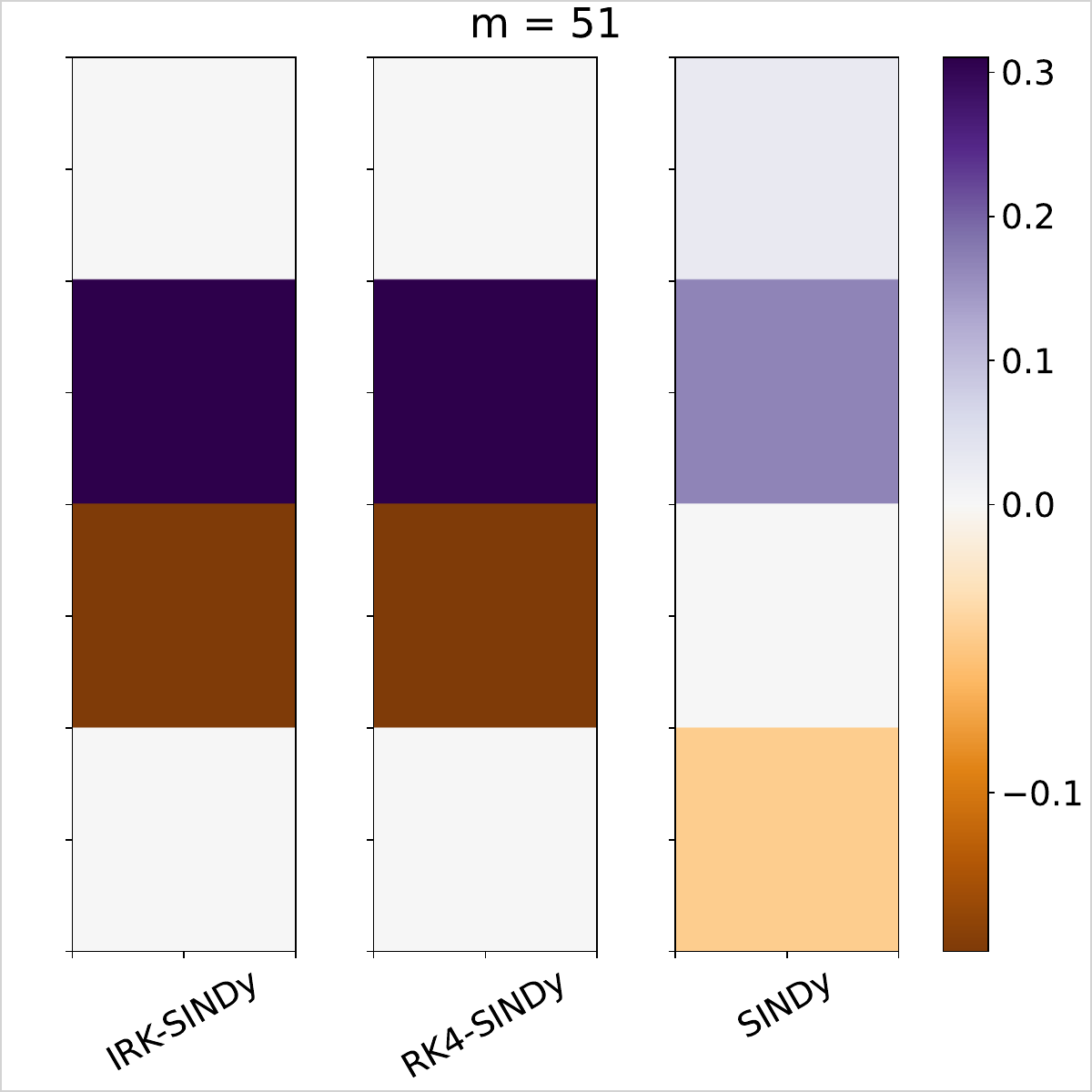}};
\node[below=of img3] (img4) {\includegraphics[width=0.45\textwidth]{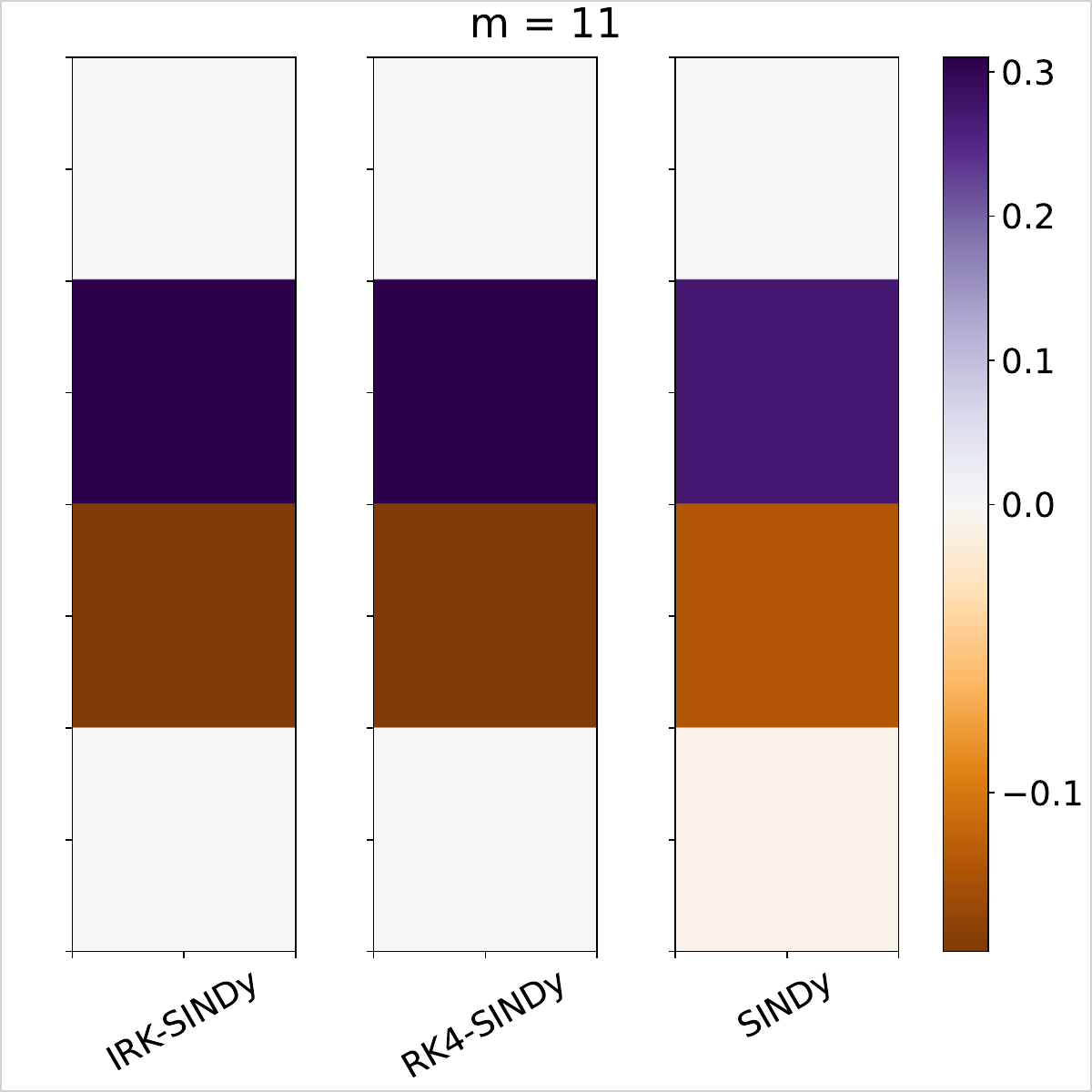}};
\end{scope}
\end{tikzpicture}
\end{subfigure}
\caption{Logistic growth model: a comparison of the reference model and recovered models using data collected at constant time stepsize. a. Phase portraits, b. coefficient matrices.}\label{fig8}
\end{figure}

In the last numerical experiment, we study the discovery of the governing equations for tumor growth. Despite the extensive advancement of various effective mathematical models within the realm of mathematical oncology, in this context, we exclude the high dimensionality and the consideration of complex tumor-immune interactions, directing our attention exclusively towards the well-established logistic growth model. This nonlinear model is a modified exponential growth model by taking into account the carrying capacity of the system, which is particularly crucial for accurately modeling the mechanisms underlying tumor growth\cite{Foryś2003}. The general form of this nonlinear dynamical system is given by eq.\eqref{Logistic}:
\begin{equation}\label{Logistic}
\dot{T}(t) = r T(t) (1 - {T(t) \over K}) = a T(t) - b T^{2}(t),
\end{equation}
where $T(t)$ signifies the temporal evolution of tumor concentration, while $r$ and $K$ represent the growth rate and carrying capacity, respectively. To generate data, we assign tumor-specific parameters of $r=0.31$ and $K=2$, conducting measurements under the initial condition $T(0)=0.1$ over the time interval $t \in [0, 50]$. We set the thresholding value to $\lambda=0.025$ and consider the polynomial space of up to order $5$ to serve as our nonlinear feature library. Throughout four successive thresholding iterations, with $20,000$ epochs allocated to the first iteration and $5,000$ epochs to each of the subsequent iterations, we simultaneously identify the active terms in the library while training the neural network. The DNN employed in this experiment comprises $3$ hidden layers, each containing $32$ neurons, utilizing the $\tanh$ activation function. We use a learning rate of $10^{-3}$ in learning $\xi$ and a learning rate of $10^{-4}$ in learning $\theta$ to discover the governing equations with varying values of $m$. Figure\ref{fig8} depicts the superior efficacy of deep IRK-SINDy in comparison to the RK4-SINDy and Conv-SINDy methodologies.

\section{Conclusion}
We have proposed an implicit Runge-Kutta based sparse identification of nonlinear dynamics, representing a new class of data-driven methods for discovering the governing equations of nonlinear dynamical systems from sparse and noisy datasets. Our innovative methodology, by integrating Gauss methods as a subclass of A-stable IRK methods by the highest accuracy with sparse regression, has demonstrated an impressive capability to directly encode the physical laws and biological mechanisms governing a specified dataset into a system of differential equations. In this work, we develop data-driven algorithms that are independent of derivative information and exhibit high robustness to data scarcity. The major challenge of these algorithms pertains to the computation of the stage values of IRKs, which has led to two general approaches: (a) iterative schemes for solving systems of nonlinear algebraic equations, including fixed-point and Newton's iterations, and (b) deep neural networks. The resultant algorithms have evidenced promising outcomes across a diverse family of benchmark problems in the field of data-driven discovery of differential equations, particularly those with biological motivation. This framework opens a new path for applying the integration of three pivotal techniques—numerical methods of solving differential equations, sparse regression, and deep learning—to directly model physical and biological phenomena from datasets.

Nevertheless, it is essential to acknowledge that the application of these algorithms to benchmark problems has revealed impracticality of approach (a)  due to the complexity of the calculations and the execution time of the algorithm. This is why approach (b) was designed and has successfully demonstrated its efficacy through the appropriate selection of neural network architecture. Deep IRK-SINDy, when applied to benchmark problems such as the two-dimensional damped harmonic oscillatory systems, FitzHugh-Nagumo model, Lorenz attractor, predator–prey dynamics, and logistic growth, has outperformed RK4-SINDy, which was similarly introduced by integrating the fourth-order Runge–Kutta method with sparse regression, as well as the conventional SINDy in the case of data scarcity. Throughout this study, it was revealed that these algorithms are resistant to noise. In this work, the Savitzky–Golay filter was employed to reduce noise and enhance the fidelity of the discovered equations.

Although this work has produced highly promising results, it is likely that the reader would concur that the number of questions engendered by this investigation significantly exceeds the answers it provides.  Which neural network architecture is optimally suited for a particular dataset? How can parametric dynamical systems, dynamical systems incorporating control terms, and equations including rational terms be effectively identified using the IRK-SINDy framework? Is the mean-squared error the most appropriate choice for the loss function? How can we develop algorithms that maintain robustness in the face of high noise levels within scenarios characterized by data scarcity, particularly considering the recent advancements in neural network architectures?

Certainly, in light of the numerous challenges present, further research is needed to establish a robust foundation in this field. Finally, the principal factor contributing to the robustness of IRK-SINDy n the context of data scarcity is related to the reduced stepsize constraints in A-stable methods and the high accuracy of IRK methods. In the future, we would like to extend the proposed framework through employing alternative high-order implicit methods that possess lower computational cost and data-independent implementations, thereby facilitating integration with sparse identification in such a way that approach (a) becomes practical and approach (b) more efficacious.

\subsection*{Data accessibility}
Our code will be attached after acceptance.

\bibliographystyle{unsrtnat}
\bibliography{IRK-SINDy}


\end{document}